\documentclass[a4paper, twoside, 11pt]{article}

\usepackage[backend=bibtex, sorting=nyt, isbn=false, url=false, doi=false, style=ieee]{biblatex}
\usepackage{amssymb, amsfonts, amsmath, amsthm}
\usepackage{hyperref} 
\usepackage{bm} 
\usepackage{graphicx} 
\usepackage[ruled, linesnumbered]{algorithm2e}
\usepackage{subcaption} 
\usepackage{multicol}
\usepackage{makecell} 
\usepackage{enumerate} 
\usepackage{indentfirst} 

\newtheorem{lemma}{Lemma}[section]
\newtheorem{theorem}{Theorem}[section]
\newtheorem{definition}{Definition}[section]
\newtheorem{corollary}{Corollary}[section]
\newtheorem{remark}{Remark}[section]
\newtheorem{example}{Example}[section]
\newtheorem{proposition}{Proposition}[section]
\numberwithin{equation}{section}
\numberwithin{figure}{section}
\numberwithin{table}{section}

\def\R{\mathbb{R}}
\def\N{\mathbb{N}}
\def\d{\mathrm{d}}
\def\norm#1{\left\lVert#1\right\rVert}
\def\B#1{\left\{#1\right\}} 
\def\abs#1{\left\lvert#1\right\rvert}
\def\jump#1{\left[\!\left[{#1}\right]\!\right]} 
\def\avg#1{\left\{\!\left\{{#1}\right\}\!\right\}} 

\allowdisplaybreaks
\title{A High-Order Local Discontinuous Galerkin Method for the \texorpdfstring{\(p\)}{p}-Laplace Equation}
\author{
	Yue~Wu\footnote{Division of Applied Mathematics, Brown University, Providence, RI 02912, USA (yue\_wu3@brown.edu).} \ 
	and Yan~Xu\footnote{Corresponding author. School of Mathematical Sciences, University of Science and Technology of China, Hefei, Anhui 230026, China (yxu@ustc.edu.cn). The research of Yan Xu was supported by NSFC grant No. 12071455.}
}
\date{} 

\begin{document}

\maketitle

\centerline{\textit{Dedicated to Prof. Chi-Wang Shu's 65th birthday}}


\begin{abstract}
	We study the high-order local discontinuous Galerkin (LDG) method for the \(p\)-Laplace equation. We reformulate our spatial discretization as an equivalent convex minimization problem and use a preconditioned gradient descent method as the nonlinear solver. For the first time, a weighted preconditioner that provides \(hk\)-independent convergence is applied in the LDG setting. For polynomial order \(k \geqslant 1\), we rigorously establish the solvability of our scheme and provide \emph{a priori} error estimates in a mesh-dependent energy norm. Our error estimates are under a different and non-equivalent distance from existing LDG results. For arbitrarily high-order polynomials under the assumption that the exact solution has enough regularity, the error estimates demonstrate the potential for high-order accuracy. Our numerical results exhibit the desired convergence speed facilitated by the preconditioner, and we observe best convergence rates in gradient variables in alignment with linear LDG, and \emph{optimal} rates in the primal variable when \(1 < p \leqslant 2\). 
\end{abstract}


\section{Introduction}

The \(p\)-Laplace equation represents a classic example of the calculus of variations in Sobolev spaces and has applications in various nonlinear physical problems. Therefore, it has attracted researchers' interest over the past decades from different disciplines. In fluid dynamics, it describes the sheer-thinning or shear-thickening effect of quasi-Newtonian fluids \cite{Lady1967New} and the nonlinear Darcy law in porous media flows \cite{originpLaplace}. In game theory, it models the tug-of-war game with noise \cite{10.1215/00127094-2008-048}. In image processing, it is used for total variation denoising \cite{10.1016/0167-2789_92_90242-F} when \(p = 1\) and optimal Lipschitz extensions \cite{10.1109/83.661188} when \(p = +\infty\). In solid mechanics, it can describe the elastic-plastic torsional creep \cite{10.1515/crll.1990.410.1}. In this study, we exclusively focus on the case \(p \in (1, +\infty)\), of which the weak form is well-posed.

Due to the nonlinear nature of the \(p\)-Laplace equation, numerical methods are challenging to design. The main reason is the singularity when \(p \in (1, 2)\) or degeneracy when \(p \in (2, +\infty)\) near stationary points where the gradient vanishes, which poses a significant challenge for linearization techniques to work well and even causes the widely used Newton-Raphson iteration to fail \cite{10.1051/m2an/197509R200411}. Therefore, numerical schemes based on minimization forms are seemingly more stable than those based on weak forms, and convex optimization algorithms are often employed for solving the nonlinear discrete system. 

Among various numerical schemes for PDEs, the traditional finite element methods (FEM) are the most widely used spatial discretization for the \(p\)-Laplace equation in the literature. The earliest one is the \(\mathcal{P}^{1}\) continuous Galerkin (CG) \cite{10.1016/S0168-2024_08_70184-X, 10.1051/m2an/197509R200411} based on a discrete minimization form. After decades, Barrett and Liu \cite{10.1090/S0025-5718-1993-1192966-4} proved that the \(\mathcal{P}^{1}\) CG converges optimally, given that the exact solution is sufficiently smooth and that the source term satisfies some constraints. Later, other FEMs for the \(p\)-Laplace equation were developed, including those belonging to nonconforming \cite{10.4208/aamm.OA-2018-0117} and mixed categories \cite{10.1093/imanum/18.1.121, 10.1216/camq/1008957338, 10.1016/j.cma.2006.11.023, 10.4208/ijnam2023-1012}. In 1996, Barrett and Liu introduced a quasi-norm technique \cite{10.1137/0733006}, and sharper error estimates \cite{10.1137/S0036142901393589, 10.1007/s00211-005-0594-5, 10.1137/070681508, 10.1093/imanum/drr016, 10.1016/j.apnum.2019.11.018} have been obtained since then. The above-mentioned researches primarily focus on the lowest-order FEMs. Nevertheless, it is worth noting that high-order and \(hp\)-FEMs offer superior approximation properties \cite{10.1007/s002110050423, AINSWORTH2000329} and thus should also be taken into consideration. 

Besides spatial discretization, another important issue is how to solve the nonlinear discrete system stably and efficiently. For primal FEMs, if the discrete equation has an equivalent minimization form, which might be a discrete analog of the variational minimization form of the \(p\)-Laplace equation, gradient descent methods could have good performance. Barrett and Liu \cite{10.1090/S0025-5718-1993-1192966-4} and Bermejo and Infante \cite{10.1137/S1064827598339098} employed the Polak-Ribi\`{e}re nonlinear conjugate gradient method without mentioning any preconditioner. In the papers \cite{IJNAM-2-123, 10.1007/s10915-007-9134-z, 10.1016/j.amc.2013.07.096}, the researchers proposed a weighted preconditioner that can provide mesh-independent (i.e. \(h\)-independent) convergence performance. However, for mixed FEMs, it's more complicated because a nonlinear saddle point system has to be solved \cite{10.1051/m2an/197509R200411, 10.1016/j.cma.2006.11.023, 10.1137/16M1067792}.

The discontinuous Galerkin (DG) methods \cite{10.1007/978-3-642-59721-3_1} allow for arbitrarily shaped mesh elements, non-conforming meshes, and totally independent elemental function spaces. During computation, inter-element communication is solely through numerical fluxes on faces, meaning that the template is compact, thus providing high parallel efficiency. Additionally, they can handle strong discontinuities like shock waves, while achieving high-order accuracy in smooth regions. The earliest DG method was proposed by Reed and Hill \cite{reed1973triangular} in 1973 to solve the neutron transport equation. In the late 1980s and 1990s, Cockburn and Shu \cite{10.1090/S0025-5718-1989-0983311-4, 0021999189901836, 10.2307/2008501, 10.1006/jcph.1998.5892} developed the RKDG method for nonlinear time-dependent hyperbolic conservation systems. Independent of DG methods for hyperbolic PDEs, those for elliptic PDEs \cite{10.1007/978-3-642-59721-3_5} were proposed in the meantime. The interior penalty DG (IPDG) method \cite{10.1137/0710071, 10.1007/BFb0120591, 10.1137/0715010} was proposed to solve linear elliptic PDEs, featuring penalty terms at faces that enforce inter-element continuity asymptotically like the Nitsche \cite{10.1007/BF02995904} method. In the 1990s, Bassi and Rebay proposed the BR scheme \cite{10.1006/jcph.1996.5572} for the diffusive part of the compressible Navier-Stokes equations. Later, Cockburn and Shu generalized the BR scheme to the LDG scheme \cite{10.1137/S0036142997316712}, using the idea that second-order PDEs can be rewritten into equivalent first-order systems. In comparison to the IPDG method, the LDG method has advantages in easy choice of penalty coefficients \cite{10.1137/S0036142900371003} and superconvergence properties \cite{10.1137/S0036142900371544, 10.1090/S0025-5718-01-01317-5}. One major drawback is the relatively high storage requirement for matrices \cite{10.1137/S1064827501388339}, but it can be alleviated by a special choice of parameters \cite{10.1137/070685518}. For a more detailed comparison of DG methods for elliptic PDEs, readers can refer to \cite{10.1137/S0036142901384162} for a theoretical perspective and \cite{10.1137/S1064827501388339} for a computational point of view. 

Besides the desire for higher accuracy and parallel efficiency, the conservation property of numerical schemes is of great importance in computational fluid dynamics. To numerically solve more sophisticated and practical PDEs like the \(p\)-Navier-Stokes equations, the \(p\)-Laplace subpart should be numerically conservative. To this end, various studies apply the DG methods to the \(p\)-Laplace equation, including LDG methods \cite{10.1016/j.crma.2008.07.005, 10.1093/imanum/drt040, 10.48550/arXiv.1709.04297, IJNAM-19-315}, IPDG methods \cite{10.1137/110820324, 10.1093/imanum/drx040}, hybridizable DG (HDG) methods \cite{10.1137/15M1008014}, and hybrid high-order (HHO) methods \cite{10.1090/mcom/3180}. In the meanwhile, theoretical tools \cite{10.1093/imanum/drn038, 10.1090/S0025-5718-10-02333-1} for broken Sobolev space analysis have been developed. However, up till now, among these methods, only the LDG \cite{10.1093/imanum/drt040}, HDG \cite{10.1007/s10915-019-00967-6} and HHO \cite{10.1142/S0218202517500191, 10.1007/s10092-021-00410-z} have \emph{a priori} error estimates for arbitrarily high-order polynomials. Nevertheless, the current LDG and HDG estimates only assume the exact solution to have the lowest regularity, so they do not reveal the potential of high-order accuracy.

In this paper, we study the LDG method for the \(p\)-Laplace equation, from derivation to analysis and numerical tests. The novelty of our work lies in the following aspects. 
\begin{enumerate}
	\item For the first time, the weighted preconditioner that is expected to have \(hk\)-independent convergence as a nonlinear solver is employed in the DG setting, whose efficiency is demonstrated in numerical results. 
	\item We present error estimates, under a different and non-equivalent distance from the existing LDG results \cite[Theorem 4.8]{10.1093/imanum/drt040}. When the exact solution has \(W^{k+1, p}\)-regularity and at least \(k\)-th (\(k \geqslant 1\)) order polynomials are used, the energy norm of the error scales like \(\mathcal{O}(h^{(p-1)k})\) when \(1 < p\leqslant 2\) and \(\mathcal{O}(h^{\frac{k}{p-1}})\) when \(p \geqslant 2\), demonstrating the potential for high-order accuracy.
	\item Furthermore, our numerical results exhibit the best convergence rates in the gradient variables in alignment with the LDG setting, and \emph{optimal} convergence rates in the primal variable when \(p \in (1, 2]\), which coincides with the improved HHO estimates \cite[Theorem 1]{10.1007/s10092-021-00410-z} and may be attributed to good local regimes of the exact solution. 
\end{enumerate}


Our paper is organized as follows. In Section \ref{se:eq}, we introduce the \(p\)-Laplace equation and some of its theoretical properties. In Section \ref{se:ldg}, we present our numerical method, from the LDG spatial discretization to a preconditioned gradient descent method solving the resulting discrete minimization problem. Section \ref{se:th} gives theoretical analysis of our proposed numerical method, which includes convergence proof and \emph{a priori} error estimates, and several numerical results are presented in Section \ref{se:nu}. We end in Section \ref{se:con} with a conclusion.

\section{The Governing Equation}\label{se:eq}

\subsection{The \texorpdfstring{\(p\)}{p}-Laplace Equation}

The primal form of the \(p\)-Laplace equation is 
\begin{equation}\label{PDE: p-Laplace equation - primal form}
	\begin{cases}
		-\nabla \cdot \left(\mathcal{A}(\nabla u)\right) = f & \text{in \(\Omega\)}, \\
		u = g_{D} & \text{on \(\Gamma^{D}\)}, \\
		\mathcal{A}(\nabla u) \cdot \bm n = \bm g_{N} \cdot \bm n & \text{on \(\Gamma^{N}\)}, 
	\end{cases}
\end{equation}
and the equivalent mixed form is 
\begin{equation}\label{PDE: p-Laplace equation - mixed form}
	\begin{cases}
		\mathcal{A}^{-1}(\bm\sigma) - \nabla u = 0 & \text{in \(\Omega\)}, \\
		-\nabla \cdot \bm\sigma = f & \text{in \(\Omega\)}, \\
		u = g_{D} & \text{on \(\Gamma^{D}\)}, \\
		\bm\sigma \cdot \bm n = \bm g_{N} \cdot \bm n & \text{on \(\Gamma^{N}\)}, 
	\end{cases}
\end{equation}
where \(\mathcal{A}: \R^{d} \mapsto \R^{d}\) is defined as \(\mathcal{A}(\bm\tau) := \abs{\bm\tau}^{p-2} \bm\tau\) with \(\abs{\cdot}\) the standard \(\ell^{2}\) norm on \(\R^{d}\), \(\Omega \subset \R^{d}\) is a polyhedral open bounded connected set, \(\Gamma^{D} \subset \partial\Omega\) has positive measure, and \(\Gamma^{D} \sqcup \Gamma^{N} = \partial\Omega\). 

The \(p\)-Laplace equation has the following behavior as \(p\) varies. \cite{10.1007/978-3-030-14501-9}
\begin{enumerate}[1)]
	\item \(p = 1\): The solution is not unique and is determined by their level sets because \(\nabla \cdot \frac{\nabla u}{\abs{\nabla u}} = \nabla \cdot \frac{\nabla \varphi(u)}{\abs{\nabla \varphi(u)}} = -H\) (the mean curvature). 
	\item \(p \in (1, +\infty)\): It is the case that we consider in this study. \begin{enumerate}
		\item \(p \in (1, 2)\): \(\mathcal{A}(\nabla u) \to \infty\) as \(\nabla u \to \bm 0\). As an elliptic PDE, the equation \eqref{p-Laplace primal weak form} is singular. 
		\item \(p = 2\): The equation \eqref{p-Laplace primal weak form} is the Poisson equation and is well-posed. 
		\item \(p \in (2, +\infty)\): \(\mathcal{A}(\nabla u) \to \bm 0\) as \(\nabla u \to \bm 0\). As an elliptic PDE, the equation \eqref{p-Laplace primal weak form} is degenerate. 
	\end{enumerate}
	\item \(p = +\infty\): The concept of viscosity solution is needed. 
\end{enumerate}

\subsection{Variational Forms of the \texorpdfstring{\(p\)}{p}-Laplace Equation}

\begin{definition}
	We define some notations for integrals. 
	\begin{enumerate}
		\item \(\norm{\cdot}_{L^{p}(K), w(x)}\): weighted \(L^{p}\) norm on the domain \(K \subset \R^{d}\) with weight function \(w(x)\). 
		\item \(\left(\cdot, \cdot\right)_{K}\) and \(\left(\cdot; w(x); \cdot\right)_{K}\): the standard and weighted dual pair on the domain \(K \subset \R^{d}\). 
		\item \(\left\langle \cdot, \cdot \right\rangle_{e}\) and \(\left\langle \cdot; w(x); \cdot \right\rangle_{e}\): the standard and weighted dual pair on the \((d-1)\)-dimensional manifold \(e\). 
	\end{enumerate}
\end{definition}

\begin{definition}[Sobolev norms]
	Let \(\Omega \subset \R^{d}\) be an open set, \(p \in [1, +\infty]\), and \(k \in \N\). The Sobolev norms for \(u \in L_{\rm loc}^{1}(\Omega)\) are defined as 
	\begin{equation}
		\norm{u}_{W^{k,p}(\Omega)} 
		:= \begin{cases}
			\left(\sum_{\abs{\bm\alpha} \leqslant k} \norm{\nabla^{\bm\alpha} u}_{L^{p}(\Omega)}^{p} \right)^{\frac{1}{p}}, & p \in [1, +\infty), \\
			\max_{\abs{\bm\alpha} \leqslant k} \norm{\nabla^{\bm\alpha} u}_{L^{\infty}(\Omega)}, & p = +\infty, 
		\end{cases}
	\end{equation}
	where \(\bm\alpha \in \N^{d}\) is a multi-index, \(\abs{\bm\alpha} = \abs{\bm\alpha}_{\ell^{1}}\), and \(\nabla\) is the weak gradient operator. 
\end{definition}

\begin{definition}[Sobolev spaces]
	Let \(\Omega \subset \R^{d}\) be an open set, \(p \in [1, +\infty]\), and \(k \in \N\). The Sobolev spaces are defined as 
	\begin{equation}
		W^{k,p}(\Omega) 
		:= \B{u \in L^{1}_{\rm loc}(\Omega) \bigg| \norm{u}_{W^{k,p}(\Omega)} < +\infty}. 
	\end{equation}
\end{definition}

\begin{definition}[Conjugate index]
	For \(p \in [1, +\infty]\), its conjugate index \(p' \in [1, +\infty]\) is defined to satisfy \(\frac{1}{p} + \frac{1}{p'} = 1\). 
\end{definition}

To introduce the variational forms of the \(p\)-Laplace equation, we choose the following spaces 
\begin{subequations}
	\begin{eqnarray}
		&& V = \B{v \in W^{1,p}(\Omega) \bigg| v\big|_{\Gamma^{D}} = g_D}, \\
		&& \widetilde{V} = \B{v \in W^{1,p}(\Omega) \bigg| v\big|_{\Gamma^{D}} = 0}. 
	\end{eqnarray}
\end{subequations}
The primal weak form of the \(p\)-Laplace equation \eqref{PDE: p-Laplace equation - primal form} is 
\begin{equation}\label{p-Laplace primal weak form}
	\begin{aligned}
		& \text{find \(u \in V\) such that \(\forall v \in \widetilde{V}\):} \\
		& \left(\mathcal{A}(\nabla u), \nabla v\right)_{\Omega} = \left(f, v\right)_{\Omega} + \left\langle \bm g_{N} \cdot \bm n, v \right\rangle_{\Gamma^{N}}, 
	\end{aligned}
\end{equation}
and the primal minimization form is 
\begin{subequations}\label{p-Laplace primal minimization form}
	\begin{equation}
		\begin{aligned}
			& \text{find \(u \in V\) such that:} \\
			& J(u) = \min_{v \in V} J(v), 
		\end{aligned}
	\end{equation}
	with the energy functional \(J: V \mapsto \R\) given as 
	\begin{equation}
		J(u) := \frac{1}{p}\norm{\nabla u}_{L^p(\Omega)}^{p} - \left(f, u\right)_{\Omega} - \left\langle \bm g_{N} \cdot \bm n, u \right\rangle_{\Gamma^{N}}. 
	\end{equation}
\end{subequations}
Here, we assume \(f \in L^{p'}(\Omega) \cap W^{-1, p'}(\Omega)\), \(g_{D} \in W^{\frac{1}{p'}, p}(\Gamma^{D})\), and \(\bm g_{N} \cdot \bm n \in L^{p'}(\Gamma^{N}) \cap W^{-\frac{1}{p'}, p'}(\Gamma^{N})\). 

\begin{proposition}[Existence, uniqueness and equivalence]
	When \(p \in (1, +\infty)\), the primal weak form \eqref{p-Laplace primal weak form} and the minimization form \eqref{p-Laplace primal minimization form} are equivalent, and their solution exists and is unique. 

	\begin{proof}
		\cite[Section 6.2]{10.1007/978-3-319-78390-1} provides a proof when \(\abs{\Gamma^{N}} = 0\) and \(g_{D} = 0\). The general proof is similar without these two conditions. 
	\end{proof}
\end{proposition}

\begin{lemma}\label{lemma: Discrete Holder inequality}
	Let \(\abs{\cdot}\) satisfy the triangle inequality (for example, a semi-norm), then 
	\begin{equation}
		\forall \alpha \in [0, +\infty), \quad 
		\abs{f + g}^{\alpha} \leqslant 2^{\max\B{\alpha-1, 0}} \left(\abs{f}^{\alpha} + \abs{g}^{\alpha}\right). 
	\end{equation}
\end{lemma}

\begin{proposition}[Energy estimates]\label{p-Laplace energy estimate}
	Let \(u \in V\) be the solution to the primal weak form \eqref{p-Laplace primal weak form} and choose \(u_{D} = E_{\gamma} E_{\Gamma^{D}}^{\partial\Omega} g_{D}\), where \(E_{\Gamma^{D}}^{\partial\Omega}: W^{\frac{1}{p'}, p}(\Gamma^{D}) \mapsto W^{\frac{1}{p'}, p}(\partial\Omega)\) is the extension operator, and \(E_{\gamma}: W^{\frac{1}{p'}, p}(\partial\Omega) \mapsto W^{1,p}(\Omega)\) is the bounded right inverse of the trace operator \cite[Theorem 3.10]{10.1007/978-3-030-56341-7}. Then we have estimates 
	\begin{equation}
		\norm{u_{D}}_{W^{1,p}(\Omega)} 
		\lesssim \norm{E_{\Gamma^{D}}^{\partial\Omega} g_{D}}_{W^{\frac{1}{p'}, p}(\partial\Omega)} 
		\lesssim \norm{g_{D}}_{W^{\frac{1}{p'}, p}(\Gamma^{D})}, 
	\end{equation}
	and 
	\begin{equation}
		\norm{u - u_{D}}_{W^{1,p}(\Omega)} 
		\lesssim \norm{\nabla u_{D}}_{L^{p}(\Omega)} + \norm{f}_{L^{p'}(\Omega)}^{\frac{1}{p-1}} + \norm{\bm g_{N} \cdot \bm n}_{L^{p'}(\Gamma^{N})}^{\frac{1}{p-1}}. 
	\end{equation}

	\begin{proof}
		The first estimate follows directly from the boundedness of these operators, so we only prove the second one. Choose \(w = u - u_{D}\) in \eqref{p-Laplace primal weak form}, then 
		\begin{align*}
			& \norm{\nabla w}_{L^{p}(\Omega)}^{p} \\
			= & \norm{\nabla u - \nabla u_{D}}_{L^{p}(\Omega)}^{p} \\
			\lesssim & \norm{\nabla u}_{L^{p}(\Omega)}^{p} + \norm{\nabla u_{D}}_{L^{p}(\Omega)}^{p} \tag{by Lemma \ref{lemma: Discrete Holder inequality}} \\
			= & \left(\mathcal{A}(\nabla u), \nabla u_{D}\right)_{\Omega} + \left(f, w\right)_{\Omega} + \left\langle \bm g_{N} \cdot \bm n, w \right\rangle_{\Gamma^{N}} + \norm{\nabla u_{D}}_{L^{p}(\Omega)}^{p} \\
			\leqslant & \norm{\nabla u}_{L^{p}(\Omega)}^{p-1} \norm{\nabla u_{D}}_{L^{p}(\Omega)} + \norm{f}_{L^{p'}(\Omega)} \norm{w}_{L^{p}(\Omega)} + \norm{\bm g_{N} \cdot \bm n}_{W^{-\frac{1}{p'}, p'}(\Gamma^{N})} \norm{w}_{W^{\frac{1}{p'}, p}(\Gamma^{N})} + \norm{\nabla u_{D}}_{L^{p}(\Omega)}^{p} \\
			\leqslant & \norm{\nabla u}_{L^{p}(\Omega)}^{p-1} \norm{\nabla u_{D}}_{L^{p}(\Omega)} + \norm{f}_{L^{p'}(\Omega)} \norm{w}_{L^{p}(\Omega)} + \norm{\bm g_{N} \cdot \bm n}_{W^{-\frac{1}{p'}, p'}(\Gamma^{N})} \norm{w}_{W^{1,p}(\Omega)} + \norm{\nabla u_{D}}_{L^{p}(\Omega)}^{p} \tag{by the trace theorem} \\
			\lesssim & \norm{\nabla w}_{L^{p}(\Omega)}^{p-1} \norm{\nabla u_{D}}_{L^{p}(\Omega)} + \norm{f}_{L^{p'}(\Omega)} \norm{w}_{L^{p}(\Omega)} + \norm{\bm g_{N} \cdot \bm n}_{W^{-\frac{1}{p'}, p'}(\Gamma^{N})} \norm{w}_{W^{1,p}(\Omega)} + \norm{\nabla u_{D}}_{L^{p}(\Omega)}^{p} \tag{by Lemma \ref{lemma: Discrete Holder inequality}}. 
		\end{align*}
		By the Poincar\'{e}-Friedrichs inequality, 
		\begin{align*}
			& \norm{w}_{W^{1,p}(\Omega)}^{p} \\
			\lesssim & \norm{w}_{W^{1,p}(\Omega)}^{p-1} \norm{\nabla u_{D}}_{L^{p}(\Omega)} + \left(\norm{f}_{L^{p'}(\Omega)} + \norm{\bm g_{N} \cdot \bm n}_{W^{-\frac{1}{p'}, p'}(\Gamma^{N})}\right) \norm{w}_{W^{1,p}(\Omega)} + \norm{\nabla u_{D}}_{L^{p}(\Omega)}^{p}. 
		\end{align*}
		Finally, the proof can be finished by a contradiction argument. 
	\end{proof}
\end{proposition}

\section{Numerical Methods}\label{se:ldg}

\subsection{Derivation of the LDG Schemes}

\subsubsection{Notations, Meshes, and Function Spaces}

In order to describe the schemes, we now state the requirements for the mesh and introduce some notations. 

Let \(\Omega \subset \R^{d}\) be a polyhedral open bounded connected set. We define \(\mathcal{T}_{h}\) the mesh, which is collection of elements \(K \subset \Omega\) such that \(\Omega = \bigsqcup_{K \in \mathcal{T}_{h}} K\). Here, each \(K \in \mathcal{T}_{h}\) is a simplex, and any two distinct \(K_1\) and \(K_2\) must either be disjoint or only share part of one common face. Any face \(e \subset \partial K\) is either an interior face, i.e. both sides are distinct elements, or a boundary face. When \(e\) is a boundary face, it must be associated with only one type of boundary condition. We use \(h\) to denote the largest diameter of elements, use \(\bm n\) to denote the unit outward normal vector of an element's face, use \(\mathcal{E}_{h}^{K}\) to denote the collection of faces of an element \(K\), use \(\mathcal{E}_{h}^{o}\), \(\mathcal{E}_{h}^{D}\), \(\mathcal{E}_{h}^{N}\), \(\mathcal{E}_{h}^{\partial}\) and \(\mathcal{E}_{h}\) to denote the collection of interior, Dirichlet, Neumann, boundary and all faces respectively, and use \(\Gamma\) to denote the union of all faces. For each face \(e\), we use \(h_{e}\) to denote the shortest normal characteristic length associated with \(e\) (could be the minimum height of neighboring elements based on this face). We require that the aspect ratios of elements are uniformly bounded, and that \(h_{e} \abs{e} \sim \abs{K}\) and \(h_{e} \sim h_{K}\) for any \(e \in \mathcal{E}_{h}^{K}\), which can be satisfied if the mesh is quasi-uniform. Lastly, we require that \(0\) is a limit point of the bounded index set \(\mathcal{H} \subset (0, h_{\rm max}]\). 

On each mesh \(\mathcal{T}_{h}\), we use the following DG function spaces 
\begin{subequations}
	\begin{eqnarray}
		&& Q_{h} := \bigoplus_{K\in\mathcal{T}_{h}} Q(K) = \bigoplus_{K\in\mathcal{T}_{h}} \left(\mathcal{P}^{k}(K)\right)^d, \\
		&& \Sigma_{h} := \bigoplus_{K\in\mathcal{T}_{h}} \Sigma(K) = \bigoplus_{K\in\mathcal{T}_{h}} \left(\mathcal{P}^{k}(K)\right)^d, \\
		&& V_{h} := \bigoplus_{K\in\mathcal{T}_{h}} V(K) = \bigoplus_{K\in\mathcal{T}_{h}} \mathcal{P}^{k}(K), 
	\end{eqnarray}
\end{subequations}
where \(\mathcal{P}^{k}(K)\) is the space of polynomials of total degree at most \(k \in \N\) supported on element \(K\), and \(\bigoplus\) denotes the direct sum of linear spaces. We use \(\nabla_{h}\) to denote the piecewise weak gradient operator, i.e. \(\nabla_{h} := \bigoplus_{K\in\mathcal{T}_{h}} \nabla\big|_{K}\). In DG schemes, boundary conditions and inter-element continuity constraints are integrated into the numerical fluxes \(\widehat{u}_{h}\) and \(\widehat{\bm\sigma}_{h}\), which are functions defined on each face \(e\) depending on \(u_{h}\) and \(\bm\sigma_{h}\) from neighboring elements. We introduce the jump and average operators that will be needed to define such fluxes. 
\begin{enumerate}
	\item \(\jump{\cdot}\): \(\forall e\in\mathcal{E}_{h}^{o}\), \(\jump{q}\big|_{e} := q_1 \bm n_1 + q_2 \bm n_2\), and \(\jump{\bm\varphi}\big|_{e} := \bm\varphi_1 \cdot \bm n_1 + \bm\varphi_2 \cdot \bm n_2\). \(\forall e\in\mathcal{E}_{h}^{\partial}\), \(\jump{q}\big|_{e} := q^{-} \bm n\), while \(\jump{\cdot}\) is undefined for vector-valued functions. 
	\item \(\avg{\cdot}\): \(\forall e\in\mathcal{E}_{h}^{o}\), \(\avg{q}\big|_{e} := \frac{q_1 + q_2}{2}\), and \(\avg{\bm\varphi}\big|_{e} := \frac{\bm\varphi_1 + \bm\varphi_2}{2}\). \(\forall e\in\mathcal{E}_{h}^{\partial}\), \(\avg{\bm\varphi}\big|_{e} := \bm\varphi^{-}\), while \(\avg{\cdot}\) is undefined for scalar functions. 
\end{enumerate}
Here, we use \(q_{j}\) and \(\bm\varphi_{j}\) with \(j \in \B{1,2}\) for traces from two neighboring elements \(K_{j}\) sharing the face \(e\) and \(\bm n_{j}\) for the outward unit normal vector of element \(K_{j}\) at face \(e\). In particular, \(\bm n_{1} + \bm n_{2} = \bm 0\).

\subsubsection{Derivation of the Weak-Form LDG Scheme}

We add a gradient variable \(\bm q\) to the mixed form \eqref{PDE: p-Laplace equation - mixed form} and use this modified version for our scheme. The equation now becomes 
\begin{equation}\label{LDG: p-Laplace equation - LDG use}
	\begin{cases}
		\bm q - \nabla u = \bm 0 & \text{in \(\Omega\)}, \\
		\bm\sigma - \mathcal{A}(\bm q) = \bm 0 & \text{in \(\Omega\)}, \\
		-\nabla \cdot \bm\sigma = f & \text{in \(\Omega\)}, \\
		u = g_{D} & \text{on \(\Gamma^{D}\)}, \\
		\bm\sigma \cdot \bm n = \bm g_{N} \cdot \bm n & \text{on \(\Gamma^{N}\)}. 
	\end{cases}
\end{equation}
Multiply the first three equations by test functions, integrate by parts, and replace the trace of the unknowns by the consistent and conservative LDG numerical fluxes \cite{10.1093/imanum/drt040, IJNAM-19-315} 
\begin{subequations}\label{p-Laplace LDG numerical flux 2}
	\begin{eqnarray}
		&& \widehat{\bm\sigma}_{h}\big|_{e} = \begin{cases}
			\avg{\bm\sigma_{h}} - \bm C_{12} \jump{\bm\sigma_{h}} - \eta \mathcal{A} \left(h_{e}^{-1} \jump{u_{h}}\right), & e \in \mathcal{E}_{h}^{o}, \\
			\bm\sigma_{h} - \eta \mathcal{A}\left(h_{e}^{-1} \left(u_{h} - g_{D}\right) \bm{n} \right), & e \in \mathcal{E}_{h}^{D}, \\
			\bm g_{N}, & e \in \mathcal{E}_{h}^{N}, 
		\end{cases} \\
		&& \widehat{u}_{h}\big|_{e} = \begin{cases}
			\avg{u_{h}} + \bm C_{12} \cdot \jump{u_{h}}, & e \in \mathcal{E}_{h}^{o}, \\
			g_{D}, & e \in \mathcal{E}_{h}^{D}, \\
			u_{h}, & e \in \mathcal{E}_{h}^{N}, 
		\end{cases}
	\end{eqnarray}
\end{subequations}
then we obtain our weak-form LDG scheme 
\begin{equation}\label{p-Laplace LDG weak form 2}
	\begin{aligned}
		& \text{find \(\left(\bm q_{h}, \bm\sigma_{h}, u_{h}\right) \in Q_{h} \times \Sigma_{h} \times V_{h}\) such that \(\forall \left(\bm\zeta_{h}, \bm\tau_{h}, v_{h}\right) \in \Sigma_{h} \times Q_{h} \times V_{h}\):} \\
		& \begin{cases}
			\left(\bm q_{h}, \bm\zeta_{h}\right)_{\Omega}
			= \left(\nabla_{h} u_{h}, \bm\zeta_{h}\right)_{\Omega} - \left\langle \jump{u_{h}}, \avg{\bm\zeta_{h}} - \bm C_{12} \jump{\bm\zeta_{h}} \right\rangle_{\Gamma^{o}} - \left\langle u_{h} - g_{D}, \bm\zeta_{h} \cdot \bm n \right\rangle_{\Gamma^{D}}, \\
			\left(\bm\sigma_{h}, \bm\tau_{h}\right)_{\Omega} = \left(\mathcal{A}(\bm q_{h}), \bm\tau_{h}\right)_{\Omega}, \\
			\begin{aligned}
				\left(\bm\sigma_{h}, \nabla_{h} v_{h}\right)_{\Omega} 
				= & \left(f, v_{h}\right)_{\Omega} + \left\langle \bm g_{N}\cdot\bm n, v_{h}\right\rangle_{\Gamma^{N}} \\
				& + \left\langle \jump{v_{h}}, \avg{\bm\sigma_{h}} - \bm C_{12} \jump{\bm\sigma_{h}}\right\rangle_{\Gamma^{o}} + \left\langle v_{h}, \bm\sigma_{h} \cdot \bm n\right\rangle_{\Gamma^{D}} \\
				& - \left\langle \eta\mathcal{A}\left(h_{e}^{-1}\jump{u_{h}}\right), \jump{v_{h}}\right\rangle_{\Gamma^{o}} - \left\langle \eta\mathcal{A}\left(h^{-1}\left(u_{h} - g_{D}\right) \bm n\right), v_{h} \bm n \right\rangle_{\Gamma^{D}}. 
			\end{aligned}
		\end{cases}
	\end{aligned}
\end{equation}
Here, \(\eta\) and \(\bm C_{12}\) are constants along each face \(e\). We assume \(\eta = \Theta(1)\) and \(\bm C_{12} = \mathcal{O}(1)\) uniformly for all \(h \in \mathcal{H}\).

\subsubsection{Conversion to the Minimization-Form LDG Scheme}

Notice that the first and third equations in the weak-form LDG scheme \eqref{p-Laplace LDG weak form 2} share some similar terms. To simplify our scheme, we introduce the following spaces and operators. 

\begin{definition}[Broken Sobolev spaces]
	Suppose \(\mathcal{T}_{h}\) is a mesh discretization of \(\Omega\). For \(k \in \N\), \(p \in [1, +\infty]\), define the broken Sobolev space \(W^{k,p}(\mathcal{T}_{h})\) to be 
	\begin{equation}
		W^{k,p}(\mathcal{T}_{h}) 
		:= \B{v \in L^{1}_{\rm loc}(\Omega) \bigg| v\big|_{K} \in W^{k,p}(K), \forall K \in \mathcal{T}_{h}}. 
	\end{equation}
\end{definition}

\begin{definition}[Interpolation operators]
	Suppose \(\mathcal{T}_{h}\) is a mesh discretization of \(\Omega\). Define the interpolation operator \(\Pi_{V_{h}}: L^{1}(\Omega) \mapsto V_{h}\) for \(u \in L^{1}(\Omega)\) as follows. 
	\begin{equation}
		\forall v_{h} \in V_{h}, \quad 
		\left(\Pi_{V_{h}} u, v_{h}\right)_{\Omega} = \left(u, v_{h}\right)_{\Omega}. 
	\end{equation}
	Besides, we also define the interpolation to \(Q_{h}\) and \(\Sigma_{h}\) similarly. 
\end{definition}

\begin{definition}[DG discrete weak gradient operators]
	Suppose \(\mathcal{T}_{h}\) is a mesh discretization of \(\Omega\). Define the bilinear operator \(D_{DG}: W^{1,1}(\mathcal{T}_{h}) \times L^{1}(\Gamma^{D}) \mapsto Q_{h}\) as follows. 
	\begin{equation}
		\begin{aligned}
			\forall \bm\zeta_{h} \in \Sigma_{h}, \quad 
			& \left(D_{DG}\left(v; g\right), \bm\zeta_{h}\right)_{\Omega} \\
			= & \left(\nabla_{h} v, \bm\zeta_{h}\right)_{\Omega} - \left\langle \jump{v}, \avg{\bm\zeta_{h}} - \bm C_{12} \jump{\bm\zeta_{h}} \right\rangle_{\Gamma^{o}} - \left\langle v - g, \bm\zeta_{h} \cdot \bm n \right\rangle_{\Gamma^{D}} \\
			= & \sum_{K \in \mathcal{T}_{h}} \left( \left( \nabla v, \bm\zeta_{h} \right)_{K} + \left\langle \widehat{v} - v, \bm\zeta_{h} \cdot \bm n \right\rangle_{\partial K} \right) \\
			= & \sum_{K \in \mathcal{T}_{h}} \left(-\left(v, \nabla \cdot \bm\zeta_{h}\right)_{K} + \left\langle \widehat{v}, \bm\zeta_{h} \cdot \bm n \right\rangle_{\partial K}\right) \\
			= & -\left(v, \nabla_{h} \cdot \bm \zeta_{h}\right)_{\Omega} + \left\langle \avg{v} + \bm C_{12} \cdot \jump{v}, \jump{\bm\zeta_{h}} \right\rangle_{\Gamma^{o}} + \left\langle v, \bm\zeta_{h} \cdot \bm n \right\rangle_{\Gamma^{N}} + \left\langle g, \bm\zeta_{h} \cdot \bm n \right\rangle_{\Gamma^{D}}. 
		\end{aligned}
	\end{equation}
	Here, \(\widehat{v}\) is the scalar numerical flux of the same form as \eqref{p-Laplace LDG numerical flux 2}, whose related Dirichlet boundary function on \(\Gamma^{D}\) is taken to be \(g\). 
\end{definition}

\begin{remark}
	\(D_{DG}(v; g)\) is in fact a discrete analogue of \(\nabla v\) with Dirichlet boundary condition \(g\), and it has the decomposition \(D_{DG}(v; g) = D_{DG}(v; 0) + D_{DG}(0; g)\). 
\end{remark}

Using the above operators, the weak-form LDG scheme \eqref{p-Laplace LDG weak form 2} can be rewritten into one equation solely involving the primal variable, 
\begin{equation}\label{p-Laplace LDG primal weak form 2}
	\begin{aligned}
		& \text{find \(u_{h} \in V_{h}\) such that \(\forall v_{h} \in V_{h}\):} \\
		& \begin{aligned}
			& \left(\mathcal{A}\left(D_{DG}\left(u_{h}; g_{D}\right)\right), D_{DG}\left(v_{h}; 0\right)\right)_{\Omega} \\& + \left\langle \eta\mathcal{A}\left(h_{e}^{-1} \jump{u_{h}}\right), \jump{v_{h}}\right\rangle_{\Gamma^{o}} + \left\langle \eta\mathcal{A}\left(h_{e}^{-1} \left(u_{h} - g_{D}\right) \bm n\right), v_{h} \bm n \right\rangle_{\Gamma^{D}} \\
			= & \left(f, v_{h}\right)_{\Omega} + \left\langle \bm g_{N} \cdot \bm n, v_{h}\right\rangle_{\Gamma^{N}}, 
		\end{aligned}
	\end{aligned}
\end{equation}
and then the gradient variables are given explicitly by 
\begin{subequations}
	\begin{eqnarray}
		&& \bm q_{h} = D_{DG}\left(u_{h}; g_{D}\right), \\
		&& \bm\sigma_{h} = \Pi_{\Sigma_{h}} \mathcal{A}\left(\bm q_{h}\right). 
	\end{eqnarray}
\end{subequations}
Its corresponding minimization form is 
\begin{subequations}\label{p-Laplace LDG minimization form 2}
	\begin{equation}
		\begin{aligned}
			& \text{find \(u_{h} \in V_{h}\) such that:} \\
			& J_{h}(u_{h}) = \min_{v_{h} \in V_{h}} J_{h}(v_{h}), 
		\end{aligned}
	\end{equation}
	where the discrete energy functional \(J_{h}: V_{h} \mapsto \R\) is 
	\begin{equation}
		\begin{aligned}
			J_{h}(v_{h}) 
			:= & \frac{1}{p}\norm{D_{DG}\left(v_{h}; g_{D}\right)}_{L^p(\Omega)}^{p} + \frac{1}{p}\norm{\jump{v_{h}}}_{L^{p}(\Gamma^{o}, \eta h_{e}^{1-p})}^{p} \\& + \frac{1}{p}\norm{v_{h} - g_{D}}_{L^{p}(\Gamma^{D}, \eta h_{e}^{1-p})}^{p} - \left(f, v_{h}\right)_{\Omega} - \left\langle \bm g_{N} \cdot \bm n, v_{h}\right\rangle_{\Gamma^{N}}. 
		\end{aligned}
	\end{equation}
\end{subequations}
The equivalence between the weak-form LDG scheme \eqref{p-Laplace LDG primal weak form 2} and minimization-form LDG scheme \eqref{p-Laplace LDG minimization form 2} will be proved rigorously.

\subsection{Solving the Minimization Problem}

While we have two forms, \eqref{p-Laplace LDG primal weak form 2} and \eqref{p-Laplace LDG minimization form 2}, for the LDG scheme, we opt for the minimization form \eqref{p-Laplace LDG minimization form 2} in practical computation due to stability concerns. Therefore, we now introduce our algorithm for solving the latter. 

\subsubsection{Preconditioned Gradient Descent}

The \(p\)-Laplace-type energy minimization problem \eqref{p-Laplace LDG minimization form 2} is either singular or degenerate, although it can always be categorized as a convex minimization problem. We should mention here that a correct choice of the preconditioner is crucial. If we adopt a na\"{\i}ve preconditioner, such as the Poisson preconditioner \cite{10.1007/s10915-007-9134-z}, or if we do not use any preconditioner \cite{10.48550/arXiv.1709.04297}, the nonlinear solver could be extremely slow. To achieve \(hk\)-independent convergence, we generalize the preconditioned gradient descent method for CG \cite{10.1007/s10915-007-9134-z} to our LDG scheme. To begin with, we calculate the G\^{a}teaux derivatives of the discrete energy functional. 

\begin{lemma}[G\^{a}teaux derivative of \(J_{h}\)]
	The first-order G\^{a}teaux derivative of \(J_{h}\) is 
	\begin{equation}\label{p-Laplace LDG2 1st order Gateaux derivative of discrete J}
		\begin{aligned}
			J_{h}'(u_{h})(v_{h}) 
			= & \left(\mathcal{A}\left(D_{DG}\left(u_{h}; g_{D}\right)\right), D_{DG}\left(v_{h}; 0\right)\right)_{\Omega} + \left\langle \eta\mathcal{A}\left(h^{-1}\jump{u_{h}}\right), \jump{v_{h}}\right\rangle_{\Gamma^{o}} \\& + \left\langle \eta\mathcal{A}\left(h^{-1}\left(u_{h} - g_{D}\right) \bm n\right), v_{h} \bm n \right\rangle_{\Gamma^{D}} - \left(f, v_{h}\right)_{\Omega} - \left\langle \bm g_{N}\cdot\bm n, v_{h}\right\rangle_{\Gamma^{N}}, 
		\end{aligned}
	\end{equation}
	
	\begin{proof}
		It follows from direct calculation. 
	\end{proof}
\end{lemma}

By the characterization of steepest descent, let \(\norm{\cdot}\) be any norm on \(V_{h}\), then the unnormalized steepest descent direction \(w_{h}\) at point \(u_{h}\) for \(J_{h}\) satisfies 
\begin{equation}\label{characterization of steepest descent}
	J_{h}'(u_{h})(w_{h}) 
	= -\norm{J'(u_{h})} \norm{w_{h}} 
	= -\sup_{v_{h} \in V_{h}}\frac{J_{h}'(u_{h})(v_{h})}{\norm{v_{h}}} \norm{w_{h}}. 
\end{equation}
As a generalization from \cite{10.1007/s10915-007-9134-z}, we suggest the following linearized energy norm for the descent direction, which can also be called a weighted norm, 
\begin{equation}\label{p-Laplace LDG minimization descent norm 2}
	\norm{w_{h}}^{2} 
	:= \begin{cases}
		\begin{aligned}
			& \norm{D_{DG}\left(w_{h}; 0\right)}_{L^2(\Omega, (\epsilon + \abs{D_{DG}\left(u_{h}; g_{D}\right)})^{p-2})}^{2} \\& + \norm{\jump{w_{h}}}_{L^2(\Gamma^{o}, \eta h_{e}^{-1}(\epsilon + \abs{h_{e}^{-1} \jump{u_{h}}})^{p-2})}^{2} \\& + \norm{w_{h}}_{L^2(\Gamma^{D}, \eta h_{e}^{-1}(\epsilon + \abs{h_{e}^{-1}(u_{h} - g_{D})})^{p-2})}^{2}, 
		\end{aligned} & p \in (1, 2), \\
		\begin{aligned}
			\norm{D_{DG}\left(w_{h}; 0\right)}_{L^2(\Omega)}^{2} + \norm{\jump{w_{h}}}_{L^2(\Gamma^{o}, \eta h_{e}^{-1})}^{2} + \norm{w_{h}}_{L^2(\Gamma^{D}, \eta h_{e}^{-1})}^{2}, 
		\end{aligned} & p = 2, \\
		\begin{aligned}
			& \norm{D_{DG}\left(w_{h}; 0\right)}_{L^2(\Omega, \epsilon + \abs{D_{DG}\left(u_{h}; g_{D}\right)}^{p-2})}^{2} \\& + \norm{\jump{w_{h}}}_{L^2(\Gamma^{o}, \eta h_{e}^{-1}(\epsilon + \abs{h_{e}^{-1} \jump{u_{h}}}^{p-2}))}^{2} \\& + \norm{w_{h}}_{L^2(\Gamma^{D}, \eta h_{e}^{-1}(\epsilon + \abs{h_{e}^{-1}(u_{h} - g_{D})}^{p-2}))}^{2}, 
		\end{aligned} & p \in (2, +\infty). 
	\end{cases}
\end{equation}
Combining \eqref{characterization of steepest descent} and \eqref{p-Laplace LDG minimization descent norm 2}, the resulting scheme to solve for the descent direction \(w_{h}\) is then: 
\begin{equation}\label{p-Laplace LDG minimization descent direction 2}
	\begin{aligned}
		& \text{find \(w_{h} \in V_{h}\) such that \(\forall v_{h} \in V_{h}\):} \\
		& -J_{h}'(u_{h})(v_{h}) = \begin{cases}
			\begin{aligned}
				& \left( D_{DG}\left(w_{h}; 0\right); (\epsilon + \abs{D_{DG}\left(u_{h}; g_{D}\right)})^{p-2}; D_{DG}\left(v_{h}; 0\right)\right)_{\Omega} \\
				& + \left\langle \eta h^{-1} \jump{w_{h}}; (\epsilon + \abs{h^{-1} \jump{u_{h}}})^{p-2}; \jump{v_{h}} \right\rangle_{\Gamma^{o}} \\
				& + \left\langle \eta h^{-1} w_{h}; (\epsilon + \abs{h^{-1} (u_{h} - g_{D})})^{p-2}; v_{h}\right\rangle_{\Gamma^{D}}, 
			\end{aligned} & p \in (1, 2), \\
			\begin{aligned}
				& \left( D_{DG}\left(w_{h}; 0\right), D_{DG}\left(v_{h}; 0\right)\right)_{\Omega} \\
				& + \left\langle \eta h^{-1} \jump{w_{h}}, \jump{v_{h}} \right\rangle_{\Gamma^{o}} + \left\langle \eta h^{-1} w_{h}, v_{h}\right\rangle_{\Gamma^{D}}, 
			\end{aligned} & p = 2, \\
			\begin{aligned}
				& \left( D_{DG}\left(w_{h}; 0\right); \epsilon + \abs{D_{DG}\left(u_{h}; g_{D}\right)}^{p-2}; D_{DG}\left(v_{h}; 0\right)\right)_{\Omega} \\
				& + \left\langle \eta h^{-1} \jump{w_{h}}; \epsilon + \abs{h^{-1} \jump{u_{h}}}^{p-2}; \jump{v_{h}} \right\rangle_{\Gamma^{o}} \\
				& + \left\langle \eta h^{-1} w_{h}; \epsilon + \abs{h^{-1} (u_{h} - g_{D})}^{p-2}; v_{h}\right\rangle_{\Gamma^{D}}, 
			\end{aligned} & p \in (2, +\infty). 
		\end{cases}
	\end{aligned}
\end{equation}
We refer to the scheme \eqref{p-Laplace LDG minimization descent direction 2} as the \emph{weighted preconditioner}, and we hope it provides \(hk\)-independent convergence, just as demonstrated numerically in \cite{10.1007/s10915-007-9134-z} for CG. 

\begin{remark}
	The weighted preconditioner scheme \eqref{p-Laplace LDG minimization descent direction 2} is actually an LDG scheme for a linear elliptic PDE. 
\end{remark}

\begin{remark}
	The stabilization term \(\epsilon > 0\) is added to avoid singularity when \(p \in (1, 2)\) and degeneracy when \(p \in (2, +\infty)\). We have to emphasize that, in the mathematical sense, the value of \(\epsilon\) does not influence the convergence and the limit \(u_{h}\). However, in practical computation, if \(\epsilon\) is too large, the convergence may become slow and dependent on the mesh. If it is too small, the condition number of the linear system solved at each iteration can be huge. 
\end{remark}

Our preconditioned gradient descent algorithm is described in Algorithm \ref{algo: Preconditioned Steepest Descent Algorithm}. Although the initial guess \(u_{h}^{(0)}\) has no impact on the convergence limit both theoretically and computationally, we opt to use the numerical solution of the Poisson equation as a reasonable one. 

\begin{algorithm}[htbp]
    \caption{Preconditioned Steepest Descent}
	\label{algo: Preconditioned Steepest Descent Algorithm}
	\SetAlgoLined

	\KwIn{\(u_{h}^{(0)}\), \(\epsilon\), \(\delta_{\rho}\), \(\delta_{w_{h}}\), \(N_{\rm it}\)}

	\KwOut{\(u_{h}^{(k)}\)}

	\tcc{main loop for gradient descent}

	\(\rho^{(0)} \leftarrow 1\) \tcp{initial guess for step size (optimal for the Poisson equation)}

	\For{\(k \leftarrow 0\) \KwTo \(N_{\rm it}-1\)}
	{   
		\(w_{h} \leftarrow\) PrecondGrad(\(u_{h}^{(k)}\), \(\epsilon\))\;

		\lIf{\(\norm{w_{h}} < \delta_{w_{h}}\)}{break} 

		\(\rho^{(k+1)} \leftarrow\) LineSearch(fun=\(J_{h}(u_{h}^{(k)} - \cdot w_{h})\), start=\(\rho^{(k)}\))\;

		\lIf{\(\rho^{(k+1)} < \delta_{\rho}\)}{break}

		\(u_{h}^{(k+1)} \leftarrow u_{h}^{(k)} - \rho^{(k+1)} w_{h}\)\;
	}
\end{algorithm}

\begin{remark}\label{remark: optimal step size}
	In Algorithm \ref{algo: Preconditioned Steepest Descent Algorithm}, we use \(\rho^{(0)} = 1\) as the initial guess for the step size because it's optimal for the case \(p = 2\) (Poisson), where convergence is reached within one iteration theoretically. 
\end{remark}

\subsubsection{Line Search}

In each iteration in Algorithm \ref{algo: Preconditioned Steepest Descent Algorithm}, the number of evaluations of the discrete energy functional \(J_{h}\) is proportional to the iteration number in the line search subroutine. Hence, a key to acceleration is to employ an efficient line search algorithm. Here, we simply choose the golden-section line search method as described in Algorithm \ref{algo: Golden Section Line Search (single-side) Algorithm}, which requires only \emph{one} functional evaluation per iteration and guarantees both monotonic descent and linear convergence rate for convex functionals. 

\begin{algorithm}[htb]
    \caption{Golden-Section Line Search (One-Sided Version)}
	\label{algo: Golden Section Line Search (single-side) Algorithm}
	\SetAlgoLined

    \KwIn{\(f: [0, +\infty) \mapsto \R\) (convex), \(x_{\rm guess} \geqslant 0\), \(\delta > 0\)}

    \KwOut{\(x_{\rm out}\), \(y_{\rm out}\)}

	\begin{multicols}{2}
    
        \(\lambda \leftarrow \frac{\sqrt{5} - 1}{2}\)\;

        \(x_{1} \leftarrow 0\), \(y_{1} \leftarrow f(x_{0})\)\;

        \eIf{\(y_{0} \leqslant y_{4}\)}
        {
            \(x_{4} \leftarrow x_{\rm guess}\), \(y_{4} \leftarrow f(x_{4})\)\;
            \(x_{2} \leftarrow (1-\lambda) x_{4}\), \(y_{2} \leftarrow f(x_{2})\)\;

            \While{\(x_{4} > \delta\)}
            {
                \eIf{\(y_{1} \leqslant y_{2}\)}
                {
                    \(x_{4} \leftarrow x_{2}\), \(y_{4} \leftarrow y_{2}\)\;
                    \(x_{2} \leftarrow (1-\lambda) x_{4}\), \(y_{2} \leftarrow f(x_{2})\)\;
                }
                {
                    break\;
                }
            }
        }
        {
            \(x_{2} \leftarrow x_{\rm guess}\), \(y_{2} \leftarrow f(x_{2})\)\;
            \(x_{4} \leftarrow \lambda^{-1} x_{2}\), \(y_{4} \leftarrow f(x_{4})\)\;

            \While{\(y_{2} > y_{4}\)}
            {
                \(x_{1} \leftarrow x_{2}\), \(y_{1} \leftarrow y_{2}\)\;
                \(x_{2} \leftarrow x_{4}\), \(y_{2} \leftarrow y_{4}\)\;
                \(x_{4} \leftarrow \lambda^{-1} x_{2}\), \(y_{4} \leftarrow f(x_{4})\)\;
            }
        }

        \(x_{3} \leftarrow \lambda x_{2} + (1-\lambda) x_{4}\), \(y_{3} \leftarrow f(x_{3})\)\;

        \While{\(x_{4} - x_{1} > \delta\)}
        {
            \eIf{\(y_{2} > y_{3}\)}
            {
                \(x_{1} \leftarrow x_{2}\), \(y_{1} \leftarrow y_{2}\)\;
                \(x_{2} \leftarrow x_{3}\), \(y_{2} \leftarrow y_{3}\)\;
                \(x_{3} \leftarrow \lambda x_{3} + (1-\lambda) x_{4}\), \(y_{3} \leftarrow f(x_{3})\)\;
            }
            {
                \(x_{4} \leftarrow x_{3}\), \(y_{4} \leftarrow y_{3}\)\;
                \(x_{3} \leftarrow x_{2}\), \(y_{3} \leftarrow y_{2}\)\;
                \(x_{2} \leftarrow (1-\lambda) x_{1} + \lambda x_{2}\), \(y_{2} \leftarrow f(x_{2})\)\;
            }
        }

        \(i \leftarrow \arg\min_{k \in \B{1,2,3,4}} y_{k}\)\;
        \(x_{\rm out} \leftarrow x_{i}\), \(y_{\rm out} \leftarrow y_{i}\)\;
	\end{multicols}
\end{algorithm}

\begin{remark}
	Our Algorithm \ref{algo: Golden Section Line Search (single-side) Algorithm} is a modified version because we have \emph{a priori} knowledge that \(f'(0) < 0\) in our context. 
\end{remark}

\section{Theoretical Analysis}\label{se:th}

In this section, we provide analysis for our LDG schemes, \eqref{p-Laplace LDG primal weak form 2} and \eqref{p-Laplace LDG minimization form 2}, including equivalence of these two schemes, solvability, energy estimates, and \emph{a priori} error estimates.

\subsection{Discrete Functional Analysis and FEM Tools}

To begin with, we present some analysis tools that are introduced specifically for our LDG schemes. 

\begin{lemma}[Basic properties for \(\mathcal{A}\)]
	Suppose \(p \in [1, +\infty)\). The operator \(\mathcal{A}: \R^d \mapsto \R^d, \bm\tau \mapsto \mathcal{A}(\bm\tau) = \abs{\bm\tau}^{p-2} \bm\tau\) has the following properties. 
	\begin{enumerate}
		\item Strict monotonicity: \(\forall \bm\tau_2 \neq \bm\tau_1\), \(\left(\bm\tau_{2} - \bm\tau_{1}\right) \cdot \left(\mathcal{A}(\bm\tau_{2}) - \mathcal{A}(\bm\tau_{1})\right) > 0\). 
		\item \((p-1)\)-homogeneity: \(\forall \lambda \in \R\), \(\mathcal{A}(\lambda \bm\tau) = \lambda^{p-1} \mathcal{A}(\bm\tau)\). 
		\item Invertibility: Iff \(p \in (1, +\infty)\), \(\mathcal{A}\) is well-defined at \(\bm\tau = \bm 0\) and is a continuous bijection on \(\R^{d}\), with \(\mathcal{A}^{-1}(\bm\tau) = \abs{\bm\tau}^{p' - 2} \bm\tau\). 
	\end{enumerate}
\end{lemma}

\begin{lemma}[Continuity and coerciveness estimates for \(\mathcal{A}\) {\cite[Lemma 2.1]{10.1090/S0025-5718-1993-1192966-4}}]\label{estimate for the A-operator}
	Suppose \(p \in (1, +\infty)\) and \(\delta \geqslant 0\). There exists constants \(C_1, C_2 > 0\) such that \(\forall \bm\tau_1 \neq \bm\tau_2\): 
	\begin{enumerate}
		\item \(\abs{\mathcal{A}(\bm\tau_1) - \mathcal{A}(\bm\tau_2)} \leqslant C_1 \abs{\bm\tau_1 - \bm\tau_2}^{1-\delta} \left(\abs{\bm\tau_1} + \abs{\bm\tau_2}\right)^{p-2+\delta}\), 
		\item \(\left(\mathcal{A}(\bm\tau_1) - \mathcal{A}(\bm\tau_2)\right) \cdot \left(\bm\tau_1 - \bm\tau_2\right) \geqslant C_2 \abs{\bm\tau_1 - \bm\tau_2}^{2 + \delta} \left(\abs{\bm\tau_1} + \abs{\bm\tau_2}\right)^{p-2-\delta}\). 
	\end{enumerate}
\end{lemma}

\begin{definition}
	For \(p \in [1, +\infty)\), define the functionals \(\abs{\cdot}_{J, p}\) and \(\norm{\cdot}_{J, p}\) on \(V_{h}\) as 
	\begin{subequations}
		\begin{eqnarray}
			&& \abs{v}_{J, p} 
			:= \left( \norm{D_{DG}\left(v; v\right)}_{L^p(\Omega)}^{p} + \norm{\jump{v}}_{L^{p}(\Gamma^{o}, \eta h_{e}^{1-p})}^{p} \right)^{\frac{1}{p}}, \\
			&& \norm{v}_{J, p} 
			:= \left(\norm{D_{DG}\left(v; 0\right)}_{L^p(\Omega)}^{p} + \norm{\jump{v}}_{L^{p}(\Gamma^{o} \cup \Gamma^{D}, \eta h_{e}^{1-p})}^{p} \right)^{\frac{1}{p}}. 
		\end{eqnarray}
	\end{subequations}
	Since \(\abs{\Gamma^{D}} > 0\), \(\abs{\cdot}_{J, p}\) is a semi-norm, and \(\norm{\cdot}_{J, p}\) is a norm. Besides, define \(E_{h}: V_{h} \times L^{1}(\Gamma^{D}) \mapsto \R\) as 
	\begin{equation}
		E_{h}(u_{h}; g_{D}) 
		:= \left(\norm{D_{DG}(u_{h}; g_{D})}_{L^{p}(\Omega)}^{p} + \norm{\jump{u_{h}}}_{L^{p}(\Gamma^{o}, \eta h_{e}^{1-p})}^{p} + \norm{u_{h} - g_{D}}_{L^{p}(\Gamma^{D}, \eta h_{e}^{1-p})}^{p}\right)^{\frac{1}{p}}. 
	\end{equation}
\end{definition}

\begin{definition}[Broken Sobolev norms {\cite{10.1093/imanum/drn038}}]
	For \(p \in [1, +\infty)\), define a semi-norm \(\abs{\cdot}_{W^{1,p}(\mathcal{T}_{h})}\) and a norm \(\norm{\cdot}_{W^{1,p}_{D}(\mathcal{T}_{h})}\) on \(W^{1,p}(\mathcal{T}_{h})\) as 
	\begin{subequations}
		\begin{eqnarray}
			&& \abs{v}_{W^{1,p}(\mathcal{T}_{h})} := \left( \norm{\nabla_{h} v}_{L^p(\Omega)}^{p} + \norm{\jump{v}}_{L^{p}(\Gamma^{o}, h_{e}^{1-p})}^{p} \right)^{\frac{1}{p}}, \\
			&& \norm{v}_{W^{1,p}_{D}(\mathcal{T}_{h})} := \left(\abs{v}_{W^{1,p}(\mathcal{T}_{h})}^{p} + \norm{v}_{L^p(\Gamma^{D}, h_{e}^{1-p})}^{p}\right)^{\frac{1}{p}}. 
		\end{eqnarray}
	\end{subequations}
\end{definition}

\begin{lemma}[Consistency error of the DG discrete weak gradient operator {\cite[Lemma 7]{10.1093/imanum/drn038}}]\label{lemma: Consistency of the DG-gradient operator}
	Suppose \(p \in [1, +\infty)\). There exists a constant \(C > 0\) independent of \(h\) such that 
	\begin{equation}
		\begin{aligned}
			& \forall v \in W^{1,p}(\mathcal{T}_{h}), \forall g \in L^{p}(\Gamma^{D}), \quad \\
			& \norm{D_{DG}(v; g) - \nabla_{h} v}_{L^{p}(\Omega)} 
			\leqslant C \left(\norm{\jump{v}}_{L^{p}(\Gamma^{o}, h_{e}^{1-p})} + \norm{v - g}_{L^{p}(\Gamma^{D}, h_{e}^{1-p})}\right). 
		\end{aligned}
	\end{equation}
\end{lemma}

\begin{lemma}[Uniform equivalence of norms]\label{p-Laplace LDG2 broken norm equivalence}
	Suppose \(p \in [1, +\infty)\). There exists a constant \(C \geqslant 1\) independent of \(h\) such that 
	\begin{equation}
		\forall v_{h} \in V_{h}, \quad 
		C^{-1} \norm{v_{h}}_{W^{1,p}_{D}(\mathcal{T}_{h})} 
		\leqslant \norm{v_{h}}_{J, p} 
		\leqslant C \norm{v_{h}}_{W^{1,p}_{D}(\mathcal{T}_{h})}. 
	\end{equation}

	\begin{proof}
		On the one hand, 
		\begin{align*}
			\norm{v_{h}}_{W^{1,p}_{D}(\mathcal{T}_{h})} 
			& \leqslant \norm{\nabla_{h} v_{h}}_{L^p(\Omega)} + \norm{\jump{v_{h}}}_{L^p(\Gamma^{o} \cap \Gamma^{D}, h_{e}^{1-p})} \\
			& \leqslant \norm{D_{DG}(v_{h}; 0)}_{L^{p}(\Omega)} + (1 + \widetilde{C}) \norm{\jump{v_{h}}}_{L^p(\Gamma^{o} \cap \Gamma^{D}, h_{e}^{1-p})} \tag{by Lemma \ref{lemma: Consistency of the DG-gradient operator}} \\
			& \leqslant 2^{\frac{1}{p'}} \max\B{1, (1 + \widetilde{C}) \eta_{\rm min}^{-\frac{1}{p}}} \norm{v_{h}}_{J, p}. \tag{by the discrete H\"{o}lder's inequality}
		\end{align*}
		On the other hand, 
		\begin{align*}
			\norm{v_{h}}_{J, p} 
			& \leqslant \norm{D_{DG}(v_{h}; 0)}_{L^{p}(\Omega)} + \norm{\jump{v_{h}}}_{L^p(\Gamma^{o} \cap \Gamma^{D}, \eta h_{e}^{1-p})} \\
			& \leqslant \norm{\nabla_{h} v_{h}}_{L^p(\Omega)} + (\eta_{\rm min}^{-\frac{1}{p}} + \widetilde{C}) \norm{\jump{v_{h}}}_{L^p(\Gamma^{o} \cap \Gamma^{D}, h_{e}^{1-p})} \tag{by Lemma \ref{lemma: Consistency of the DG-gradient operator}} \\
			& \leqslant 2^{\frac{1}{p'}} \max\B{1, \eta_{\rm min}^{-\frac{1}{p}} + \widetilde{C}} \norm{v_{h}}_{W^{1,p}_{D}(\mathcal{T}_{h})}. \tag{by the discrete H\"{o}lder's inequality}
		\end{align*}
		Therefore, \(\norm{\cdot}_{W^{1,p}_{D}(\mathcal{T}_{h})}\) and \(\norm{\cdot}_{J, p}\) are uniformly equivalent. 
	\end{proof}
\end{lemma}

\begin{lemma}[Uniform equivalence of semi-norms]\label{p-Laplace LDG2 broken semi-norm equivalence}
	Suppose \(p \in [1, +\infty)\). There exists a constant \(C \geqslant 1\) independent of \(h\) such that 
	\begin{equation}
		\forall v_{h} \in V_{h}, \quad 
		C^{-1} \abs{v_{h}}_{W_{D}^{1,p}(\mathcal{T}_{h})} 
		\leqslant \abs{v_{h}}_{J, p} 
		\leqslant C \abs{v_{h}}_{W_{D}^{1,p}(\mathcal{T}_{h})}. 
	\end{equation}
	
	\begin{proof}
		We omit the proof because it is similar to that of Lemma \ref{p-Laplace LDG2 broken norm equivalence}. 
	\end{proof}
\end{lemma}

\begin{lemma}[H\"{o}lder-type energy norm inequality]\label{p-Laplace LDG2 Holder-type energy norm inequality}
	Let \(1 \leqslant q \leqslant p < +\infty\). There exists a constant \(C > 0\) independent of \(h\) and \(\abs{\Omega}\) (may depend on the shape of \(\Omega\)) such that 
	\begin{equation}
		\forall v \in W^{1,p}(\mathcal{T}_{h}), \quad 
		\norm{v}_{J, q} \leqslant C \abs{\Omega}^{\frac{1}{q} - \frac{1}{p}} \norm{v}_{J, p}. 
	\end{equation}

	\begin{proof}
		\begin{align*}
			\norm{v}_{J, q}^{q}
			& = \norm{D_{DG}\left(v; 0\right)}_{L^q(\Omega)}^{q} + \norm{\jump{v}}_{L^{q}(\Gamma^{o} \cup \Gamma^{D}, \eta h_{e}^{1-q})}^{q} \\
			& \leqslant \abs{\Omega}^{1 - \frac{q}{p}} \norm{D_{DG}\left(v; 0\right)}_{L^p(\Omega)}^{q} + \sum_{e \in \mathcal{E}_{h}^{o} \cup \mathcal{E}_{h}^{D}} \norm{\jump{v}}_{L^{p}(e, \eta h_{e}^{1-p})}^{q} \left(\eta h_{e} \abs{e}\right)^{1-\frac{q}{p}} \tag{by H\"{o}lder's inequality} \\
			& \leqslant \abs{\Omega}^{1 - \frac{q}{p}} \norm{D_{DG}\left(v; 0\right)}_{L^p(\Omega)}^{q} + \frac{1}{2} \sum_{K \in \mathcal{T}_{h}} \sum_{e \in \mathcal{E}_{h}^{K} \backslash \mathcal{E}_{h}^{N}} \norm{\jump{v}}_{L^{p}(e, \eta h_{e}^{1-p})}^{q} \left(\eta h_{e} \abs{e}\right)^{1-\frac{q}{p}} \\
			& \leqslant \abs{\Omega}^{1 - \frac{q}{p}} \norm{D_{DG}\left(v; 0\right)}_{L^p(\Omega)}^{q} + \frac{1}{2} \sum_{K \in \mathcal{T}_{h}} \sum_{e \in \mathcal{E}_{h}^{K} \backslash \mathcal{E}_{h}^{N}} \norm{\jump{v}}_{L^{p}(e, \eta h_{e}^{1-p})}^{q} \left(c_2 \eta_{\rm max} \abs{K}\right)^{1-\frac{q}{p}} \tag{by mesh assumption} \\
			& \leqslant \abs{\Omega}^{1 - \frac{q}{p}} \left(\norm{D_{DG}\left(v; 0\right)}_{L^p(\Omega)}^{q} + \left( (d+1) c_{2} \eta_{\rm max} \right)^{1-\frac{q}{p}} \norm{\jump{v}}_{L^{p}(\Gamma^{o}\cup\Gamma^{D}, \eta h_{e}^{1-p})}^{p}\right) \tag{by H\"{o}lder's inequality} \\
			& \leqslant C \abs{\Omega}^{1 - \frac{q}{p}} \norm{v}_{J, p}^{q}. \tag{by the discrete H\"{o}lder's inequality} 
		\end{align*}
	\end{proof}
\end{lemma}

\begin{definition}[Sobolev conjugate index]
	In \(\R^{d}\), the Sobolev conjugate index for \(p \in [1, d)\) is defined to be \(p^{\ast}\) such that \(\frac{1}{p^{\ast}} = \frac{1}{p} - \frac{1}{d}\). 
\end{definition}

\begin{lemma}[Broken Poincar\'{e}-Friedrichs inequality]\label{p-Laplace LDG2 broken Poincare-Friedrichs inequality}
	If \(d > 1\), suppose \(p \in [1, d)\) and \(q \in [1, p^{\ast}]\), or \(p \in [d, +\infty)\) and \(q \in [1, +\infty)\). If \(d = 1\), suppose \(p \in [1, +\infty)\) and \(q \in [1, +\infty]\). There exists a constant \(C > 0\) independent of \(h\) such that 
	\begin{equation}
		\forall v_{h} \in V_{h}, \quad 
		\norm{v_{h}}_{L^{q}(\Omega)} \leqslant C \norm{v_{h}}_{J, p}. 
	\end{equation}

	\begin{proof}
		Case \(d > 1\). When \(p \in [1, d)\) and \(q \in [1, p^{\ast}]\), 
		\begin{align*}
			\norm{v_{h}}_{L^{q}(\Omega)}
			& \leqslant \abs{\Omega}^{\frac{1}{q} - \frac{1}{p^{\ast}}} \norm{v_{h}}_{L^{p^{\ast}}(\Omega)} \tag{by H\"{o}lder's inequality} \\
			& \lesssim \left(\norm{v_{h}}_{L^{p}(\Gamma^{D})} + \abs{v_{h}}_{W^{1,p}(\mathcal{T}_{h})}\right) \tag{by {\cite[Corollary 4.3]{10.1093/imanum/drn038}}} \\
			& \lesssim \left(\norm{v_{h}}_{L^{p}(\Gamma^{D})} + \norm{v_{h}}_{W^{1,p}_{D}(\mathcal{T}_{h})}\right) \\
			& \lesssim \left(\norm{v_{h}}_{L^{p}(\Gamma^{D})} + C_{1} \norm{v_{h}}_{J, p}\right) \tag{by Lemma \ref{p-Laplace LDG2 broken norm equivalence}} \\
			& \lesssim \norm{v_{h}}_{J, p}. 
		\end{align*}
		When \(p \in [d, +\infty)\) and \(q \in [1, +\infty)\), choose \(r \in (1, d)\) such that \(r^{\ast} \geqslant q\), 
		\begin{align*}
			\norm{v_{h}}_{L^{q}(\Omega)}
			& \leqslant \abs{\Omega}^{\frac{1}{q} - \frac{1}{r^{\ast}}} \norm{v_{h}}_{L^{r^{\ast}}(\Omega)} \tag{by H\"{o}lder's inequality} \\
			& \lesssim \norm{v_{h}}_{J, r} \tag{by Lemma \ref{p-Laplace LDG2 Holder-type energy norm inequality}} \\
			& \lesssim \norm{v_{h}}_{J, p}. \tag{by the last sub-case}
		\end{align*}

		Case \(d = 1\). 
		\begin{align*}
			\norm{v_{h}}_{L^{q}(\Omega)}
			& \leqslant \abs{\Omega}^{\frac{1}{q}} \norm{v_{h}}_{L^{\infty}(\Omega)} \tag{by H\"{o}lder's inequality} \\
			& \leqslant \abs{\Omega}^{\frac{1}{q}} \left( \norm{\nabla_{h} v_{h}}_{L^{1}(\Omega)} + \norm{\jump{v_{h}}}_{L^{1}(\Gamma^{o}\cup\Gamma^{D})} \right) \tag{by Newton-Leibniz formula} \\
			& \lesssim \norm{v_{h}}_{J, 1} \\
			& \lesssim \norm{v_{h}}_{J, p}. 
		\end{align*}
	\end{proof}
\end{lemma}

\begin{lemma}[Broken trace inequality]\label{p-Laplace LDG2 broken trace inequality}
	If \(d > 1\), suppose \(p \in [1, d]\) and \(q \in [1, \frac{p(d-1)}{d-p}]\), or \(p \in (d, +\infty)\) and \(q \in [1, +\infty]\). If \(d = 1\), suppose \(p \in [1, +\infty)\) and \(q \in [1, +\infty]\). There exists a constant \(C > 0\) independent of \(h\) such that 
	\begin{equation}
		\forall v_{h} \in V_{h}, \quad 
		\norm{v_{h}}_{L^{q}(\partial\Omega)} \leqslant C \norm{v_{h}}_{J, p}. 
	\end{equation}
	
	\begin{proof}
		Case \(d > 1\). When \(p \in [1, d]\) and \(q \in [1, \frac{p(d-1)}{d-p}]\), let \(\widetilde{p} := \frac{p(d-1)}{d-p}\), 
		\begin{align*}
			\norm{v_{h}}_{L^{q}(\partial\Omega)} 
			& \leqslant \abs{\partial\Omega}^{\frac{1}{q} - \frac{1}{\widetilde{p}}} \norm{v_{h}}_{L^{\widetilde{p}}(\partial\Omega)} \tag{by H\"{o}lder's inequality} \\
			& \lesssim \norm{v_{h}}_{L^{1}(\Omega)} + \abs{v_{h}}_{W^{1,p}(\mathcal{T}_{h})} \tag{by {\cite[Theorem 4.4]{10.1093/imanum/drn038}}} \\
			& \lesssim \norm{v_{h}}_{J, p}. \tag{by Lemma \ref{p-Laplace LDG2 Holder-type energy norm inequality} and \ref{p-Laplace LDG2 broken Poincare-Friedrichs inequality}} 
		\end{align*}
		When \(p \in (d, +\infty)\) and \(q \in [1, +\infty]\), 
		\begin{align*}
			\norm{v_{h}}_{L^{q}(\partial\Omega)} 
			& \leqslant \abs{\partial \Omega}^{\frac{1}{q}} \norm{v_{h}}_{L^\infty(\partial\Omega)} \tag{by H\"{o}lder's inequality} \\
			& \lesssim \norm{v_{h}}_{J, d} \tag{by {\cite[Theorem 4.4]{10.1093/imanum/drn038}}} \\
			& \lesssim \widehat{C}_{2} \norm{v_{h}}_{J, p}. \tag{by Lemma \ref{p-Laplace LDG2 Holder-type energy norm inequality}} 
		\end{align*}

		Case \(d = 1\) directly follows from Lemma \ref{p-Laplace LDG2 broken Poincare-Friedrichs inequality}. 
	\end{proof}
\end{lemma}

\begin{remark}
	These results are similar to that of \cite{10.1093/imanum/drn038, 10.1090/S0025-5718-10-02333-1} but are in a different norm. The estimates in our norm are essential to our energy norm based error estimates presented later because the norm coincides with the nonlinear part of the discrete energy functional. 
\end{remark}

Besides the discrete functional analysis tools, we also introduce some results from FEMs before presenting our analysis. 

\begin{lemma}[Estimates for polynomial interpolation on affine elements {\cite[Theorem 5]{10.1007/BF00252458}}]\label{thm: polynomial approximation}
	Let \(K\) be any affine element whose reference element is \(\widehat{K}\). Suppose \(p, q \in [1, +\infty]\) and \(k, m \in \N\). Assume a linear operator \(\widehat{\Pi} \in L(W^{k+1, p}(\widehat{K}), W^{m, q}(\widehat{K}))\) satisfies \(\mathcal{P}^{k}(\widehat{K}) \subset \ker (I - \widehat{\Pi})\). Let \(\Pi := \Phi_{\ast} \circ \widehat{\Pi} \circ \Phi^{\ast}\), where \(\Phi^{\ast}\) satisfies \((\Phi^{\ast} u)(\bm\lambda) = u(\Phi(\bm\lambda))\), \(\Phi_{\ast} = (\Phi^{\ast})^{-1}\), and \(\Phi\) is the affine mapping from \(\widehat{K}\) to \(K\). Then there exists a constant \(C > 0\) such that 
	\begin{equation}
		\forall v \in W^{k+1, p}(K) \cap W^{m, q}(K), \quad 
		\abs{v - \Pi v}_{W^{m, q}(K)} 
		\leqslant C \abs{K}^{\frac{1}{q} - \frac{1}{p}} \frac{h_{K}^{k+1}}{\rho_{K}^{m}} \abs{v}_{W^{k+1, p}(K)}. 
	\end{equation}
	Here, \(h_{K}\) is the minimal diameter of any ball being a superset of \(K\), and \(\rho_{K}\) is the maximal diameter of any ball being a subset of \(K\). 
\end{lemma}

\begin{lemma}[Multiplicative trace inequality for simplices]\label{thm: Multiplicative trace inequality}
	Suppose \(K \subset \R^{d}\) is a simplex. There exists a constant \(c_{d} > 0\) that only depends on \(d\) such that \(\forall p \in [1, +\infty]\), \(\forall v \in W^{1,p}(K)\), \(\forall e\) being a face of \(K\), 
	\begin{equation}
		\norm{v}_{L^{p}(e)} 
		\leqslant c_{d} \norm{v}_{L^{p}(K)}^{1 - \frac{1}{p}} \left(\left(l_{e}^{\perp}\right)^{-\frac{1}{p}} \norm{v}_{L^{p}(K)}^{\frac{1}{p}} + \norm{\nabla v}_{L^{p}(K)}^{\frac{1}{p}}\right), 
	\end{equation}
	where \(l_{e}^{\perp}\) is the height on \(e\). 

	\begin{proof}
		It follows directly from the proof of \cite[Lemma 12.15]{10.1007/978-3-030-56341-7}. 
	\end{proof}
\end{lemma}

\subsection{Equivalence, Uniqueness and Existence}

\subsubsection{The Main Results}

With the above tools, we analyze the solvability of our LDG schemes. Our main results are stated as follows. 

\begin{theorem}[Existence of the minimizer]\label{thm: existence of the minimizer}
	For \(p \in [1, +\infty)\), the solution to the minimization-form LDG scheme \eqref{p-Laplace LDG minimization form 2} exists. 
\end{theorem}

\begin{theorem}[Uniqueness of the minimizer]\label{thm: uniqueness of the minimizer}
	For \(p \in (1, +\infty)\), the solution to the minimization-form LDG scheme \eqref{p-Laplace LDG minimization form 2} is unique. 
\end{theorem}

\begin{theorem}[Equivalence of the two primal forms]\label{thm: equivalence of the two forms}
	For \(p \in (1, +\infty)\), the solutions of the two primal forms of LDG schemes, \eqref{p-Laplace LDG primal weak form 2} and \eqref{p-Laplace LDG minimization form 2}, are the same. 
\end{theorem}

\subsubsection{Proofs}

Next, we are going to provide our proof. 

\begin{lemma}[Coerciveness]\label{p-Laplace LDG2 discrete coerciveness}
	Suppose \(p \in [1, +\infty)\). For any fixed \(\mathcal{T}_{h}\), the discrete energy functional \(J_{h}\) is coercive with respect to the norm \(\norm{\cdot}_{J, p}\), i.e. 
	\begin{equation}
		\lim_{\norm{v_{h}}_{J, p} \to +\infty} J_{h}(v_{h}) = +\infty. 
	\end{equation}

	\begin{proof}
		Fix any \(u_{h,0} \in V_{h}\), and let \(w_{h} := v_{h} - u_{h,0}\), then 
		\begin{align*}
			J_{h}(v_{h}) 
			= & \frac{1}{p} E_{h}(v_{h}; g_{D})^{p} - \left(f, v_{h}\right)_{\Omega} - \left\langle \bm g_{N} \cdot \bm n, v_{h} \right\rangle_{\Gamma^{N}} \\
			\geqslant & \frac{2^{1-p}}{p} \norm{w_{h}}_{J, p}^{p} - \left(f, w_{h}\right)_{\Omega} - \left\langle \bm g_{N} \cdot \bm n, w_{h} \right\rangle_{\Gamma^{N}} \\& - \frac{1}{p} E_{h}(u_{h,0}; g_{D})^{p} - \left(f, u_{h,0}\right)_{\Omega} - \left\langle \bm g_{N} \cdot \bm n, u_{h,0} \right\rangle_{\Gamma^{N}} \tag{by Lemma \ref{lemma: Discrete Holder inequality}} \\
			\geqslant & \frac{2^{1-p}}{p} \norm{w_{h}}_{J, p}^{p} - \norm{f}_{L^{p'}(\Omega)} \norm{w_{h}}_{L^{p}(\Omega)} - \norm{\bm g_{N} \cdot \bm n}_{L^{p'}(\Gamma^{N})} \norm{w_{h}}_{L^{p}(\Gamma^{N})} \\& - \frac{1}{p} E_{h}(u_{h,0}; g_{D})^{p} - \left(f, u_{h,0}\right)_{\Omega} - \left\langle \bm g_{N} \cdot \bm n, u_{h,0} \right\rangle_{\Gamma^{N}} \tag{by H\"{o}lder's inequality} \\
			\geqslant & \frac{2^{1-p}}{p} \norm{w_{h}}_{J, p}^{p} - C_{1} \norm{f}_{L^{p'}(\Omega)} \norm{w_{h}}_{J, p} - C_{2} \norm{\bm g_{N} \cdot \bm n}_{L^{p'}(\Gamma^{N})} \norm{w_{h}}_{J, p} \\& - \frac{1}{p} E_{h}(u_{h,0}; g_{D})^{p} - \left(f, u_{h,0}\right)_{\Omega} - \left\langle \bm g_{N} \cdot \bm n, u_{h,0} \right\rangle_{\Gamma^{N}}. \tag{by Lemma \ref{p-Laplace LDG2 broken Poincare-Friedrichs inequality} and \ref{p-Laplace LDG2 broken trace inequality}} 
		\end{align*}
		Here, the constants \(C_{1}\) and \(C_{2}\) are independent of \(h\), so \(J_{h}(v_{h}) \to +\infty\) as \(\norm{v_{h}}_{J, p} \to +\infty\). 
	\end{proof}
\end{lemma}

\begin{lemma}[Weak lower semi-continuity]
	The discrete energy functional \(J_{h}\) is weakly lower semi-continuous. 
\end{lemma}

\begin{proof}[Proof of Theorem \ref{thm: existence of the minimizer}]
	With the coerciveness and the weak lower semi-continuity, Theorem \ref{thm: existence of the minimizer} follows directly from classical results in calculus of variations. 
\end{proof}

\begin{lemma}[Strict convexity]
	Suppose \(p \in (1, +\infty)\). The discrete energy functional \(J_{h}\) is strictly convex. 

	\begin{proof}
		Assume there exists \(\lambda \in (0,1)\) such that 
		\begin{align*}
			\lambda J_{h}(v_{h}^{(1)}) + (1-\lambda)J_{h}(v_{h}^{(2)}) = J_{h}(\lambda v_{h}^{(1)} + (1-\lambda) v_{h}^{(2)}). 
		\end{align*}
		The strict convexity of the mapping \(t \mapsto \abs{t}^{p}\) yields \(\left(v_{h}^{(1)} - v_{h}^{(2)}\right)\big|_{\Gamma^{D}} = 0\), \(\jump{v_{h}^{(1)} - v_{h}^{(2)}}\big|_{\Gamma^{o}} = \bm 0\), and \(D_{DG}(v_{h}^{(1)} - v_{h}^{(2)}; 0) = \bm 0\). Then \(\nabla_{h} (v_{h}^{(1)} - v_{h}^{(2)}) = \bm 0\) and thus \(v_{h}^{(1)} = v_{h}^{(2)}\). 
	\end{proof}
\end{lemma}

\begin{proof}[Proof of Theorem \ref{thm: uniqueness of the minimizer}]
	It directly follows from the strict convexity of \(J_{h}\). 
\end{proof}

\begin{proof}[Proof of Theorem \ref{thm: equivalence of the two forms}]
	We first prove that the solution of \eqref{p-Laplace LDG minimization form 2} is a solution of \eqref{p-Laplace LDG primal weak form 2}. In fact, using \eqref{p-Laplace LDG2 1st order Gateaux derivative of discrete J}, we know the 1st-order G\^{a}teaux derivative of \(J_{h}\) vanishes at \(u_{h}\), which means \eqref{p-Laplace LDG primal weak form 2} holds. 
	
	Therefore, \eqref{p-Laplace LDG primal weak form 2} has solutions. We now only need to prove that every solution \(u_{h}\) of \eqref{p-Laplace LDG primal weak form 2} is a solution of \eqref{p-Laplace LDG minimization form 2}. That \(u_{h}\) solves \eqref{p-Laplace LDG primal weak form 2} gives \(J_{h}'(u_{h}) = 0\). The convexity of \(J_{h}\) yields 
	\begin{align*}
		\forall v_{h} \in V_{h}, \forall t \in \R, \quad 
		J_{h}(u_{h} + t v_{h}) 
		\geqslant J_{h}(u_{h}) + t J_{h}'(u_{h})(v_{h}) 
		= J_{h}(u_{h}), 
	\end{align*}
	then \(u_{h}\) also solves \eqref{p-Laplace LDG minimization form 2}. 
\end{proof}

So far, we have established the relationship between the discrete solutions of the LDG schemes \eqref{p-Laplace LDG primal weak form 2} and \eqref{p-Laplace LDG minimization form 2} on a single mesh.

\subsection{Boundedness and Convergence}

On a series of meshes, we provide the following boundedness and convergence results of the numerical solutions. 

\begin{proposition}[Energy estimates for the discrete solution]\label{p-Laplace LDG2 energy-norm estimate}
	Let \(u_{h}\) be the solution of the LDG scheme \eqref{p-Laplace LDG primal weak form 2} or \eqref{p-Laplace LDG minimization form 2}. Let \(u_{h,D} := \Pi_{V_{h}} u_{D}\), where \(u_{D} = E_{\gamma} E_{\Gamma^{D}}^{\partial\Omega} g_{D}\). Then there exists a constant \(C > 0\) independent of \(h\) such that 
	\begin{subequations}
		\begin{equation}
			E_{h}(u_{h}; g_{D}) + \norm{u_{h} - u_{h,D}}_{J, p} 
			\leqslant C \left( \norm{g_{D}}_{W^{\frac{1}{p'},p}(\Gamma^{D})} + \norm{f}_{L^{p'}(\Omega)}^{\frac{1}{p-1}} + \norm{\bm g_{N} \cdot \bm n}_{L^{p'}(\Gamma^{N})}^{\frac{1}{p-1}} \right). 
		\end{equation}
	\end{subequations}

	\begin{proof}
		See {\cite[Theorem 3.2]{10.1093/imanum/drt040}} for the proof with the Dirichlet boundary condition. The generalization to our case is trivial. 
	\end{proof}
\end{proposition}

\begin{lemma}[Convergence]
	If \(d > 1\), suppose \(p \in (1, d)\) and \(q \in [1, p^{\ast})\) and \(r \in [1, \frac{(d-1)p}{d-p})\), or \(p \in [d, +\infty)\) and \(q \in [1, +\infty)\) and \(r \in [1, +\infty)\). If \(d = 1\), suppose \(p \in (1, +\infty)\) and \(q \in [1, +\infty)\) and \(r \in [1, +\infty]\). Let \(k \geqslant 1\) be the minimal polynomial degree used in the LDG scheme \eqref{p-Laplace LDG primal weak form 2} or \eqref{p-Laplace LDG minimization form 2}, and let \(u_{h}\) be the discrete solution. Let \(u \in W^{1, p}(\Omega)\) be the solution to the \(p\)-Laplace equaiton \eqref{p-Laplace primal weak form}. Then as \(h \to 0^{+}\), 
	\begin{enumerate}
		\item \(u_{h} \to u\) strongly in \(L^{q}(\Omega)\), 
		\item \(u_{h}\big|_{\partial\Omega} \to u\big|_{\partial\Omega}\) strongly in \(L^{r}(\partial\Omega)\), 
		\item \(D_{DG}(u_{h}; g_{D}) \to \nabla u\) strongly in \(L^{p}(\Omega)\), 
		\item \(\norm{u_{h} - u}_{L^{p}(\Gamma^{o} \cup \Gamma^{D}, h_{e}^{1-p})} \to 0\), 
		\item \(J_{h}(u_{h}) \to J(u)\). 
	\end{enumerate}
	
	\begin{proof}
		It follows from \cite[Lemma 8 and Theorem 6.1]{10.1093/imanum/drn038}. 
	\end{proof}
\end{lemma}

Up till now, we have only required the minimal regularity for the exact solution \(u\), i.e. \(u \in W^{1,p}(\Omega)\).

\subsection{\emph{A Priori} Error Estimates}

Finally, leveraging the aforementioned results, we give error estimates for our LDG schemes. 

\begin{definition}[Quasi-norms \cite{10.1090/S0025-5718-1993-1192966-4}]
	For functions \(f, g \in L^{p}(D)\) we define the notation 
	\begin{equation}
		\abs{f}_{g, p, \sigma, D, w} 
		:= \int_{D} \left(\abs{f} + \abs{g}\right)^{p-\sigma} \abs{f}^{\sigma} w(x) \d x. 
	\end{equation}
	Here, \(w\) can be omitted in the notation if \(w(x) \equiv 1\). 
\end{definition}

\begin{lemma}[{\cite[Lemma 2.2]{10.1090/S0025-5718-1993-1192966-4}}]\label{lemma: paper func inequality}
	For \(f, g \in L^{p}(D)\) and \(\sigma > 1\), the following estimates hold. 
	\begin{subequations}
		\begin{eqnarray}
			& \abs{f}_{g, p, \sigma, D}^{\frac{\sigma}{p}} 
			\leqslant \norm{f}_{L^{p}(D)}^{\sigma} 
			\leqslant C \left(\norm{f}_{L^{p}(D)} + \norm{g}_{L^{p}(D)}\right)^{\sigma - p} \abs{f}_{g, p, \sigma, D}, & p \in (1, \sigma], \\
			& \norm{f}_{L^{p}(D)}^{p} 
			\leqslant \abs{f}_{g, p, \sigma, D} 
			\leqslant C \left(\norm{f}_{L^{p}(D)} + \norm{g}_{L^{p}(D)}\right)^{p - \sigma} \norm{f}_{L^{p}(D)}^{\sigma}, & p \in [\sigma, +\infty). 
		\end{eqnarray}
	\end{subequations}
\end{lemma}

\begin{lemma}[The discrete equation satisfied by the exact solution \cite{10.1093/imanum/drt040}]
	Suppose \(u \in W^{1,p}(\Omega)\) is the exact solution to the weak form \eqref{p-Laplace primal weak form}. Let \(\bm q := \nabla u \in L^{p}(\Omega)\) and \(\bm\sigma := \mathcal{A}(\bm q) \in L^{p'}(\mathrm{div}, \Omega)\). Then the following equations hold. 
	\begin{equation}\label{p-Laplace discrete equation satisfied by the exact solution}
		\begin{aligned}
			& \forall \left(\bm\zeta_{h}, \bm\tau_{h}, v_{h}\right) \in \Sigma_{h} \times Q_{h} \times V_{h}: \\
			& \begin{cases}
				\left(\bm q, \bm\zeta_{h}\right)_{\Omega} = \left(\nabla_{h} u, \bm\zeta_{h}\right)_{\Omega}, \\
				\left(\bm\sigma, \bm\tau_{h}\right)_{\Omega} = \left(\mathcal{A}(\bm q), \bm\tau_{h}\right)_{\Omega}, \\
				\begin{aligned}
					\left(\bm\sigma, D_{DG}(v_{h}; 0)\right)_{\Omega} 
					= & \left(f, v_{h}\right)_{\Omega} + \left\langle \bm g_{N} \cdot \bm n, v_{h} \right\rangle_{\Gamma^{N}} \\& - \left\langle \avg{\Pi_{\Sigma_{h}} \bm\sigma} - \bm C_{12} \jump{\Pi_{\Sigma_{h}} \bm\sigma} - \bm\sigma, \jump{v_{h}} \right\rangle_{\Gamma^{o} \cup \Gamma^{D}}. 
				\end{aligned}
			\end{cases}
		\end{aligned}
	\end{equation}
\end{lemma}

\begin{lemma}[Estimates for the projected exact solution]\label{projected solution estimate}
	Let \(k \in \N\) be the minimal polynomial degree used in the LDG scheme \eqref{p-Laplace LDG primal weak form 2} or \eqref{p-Laplace LDG minimization form 2}. Let \(u\) be the solution to the \(p\)-Laplace equation \eqref{p-Laplace primal weak form}. Assume \(u \in W^{s+1, p}(\Omega)\), where WLOG we assume \(s \in \N\) satisfies \(s \leqslant k\). Let \(u_{h}^{\ast} := \Pi_{V_{h}} u\). Then \(\forall l \in \N\) such that \(l \leqslant s\), 
	\begin{subequations}
		\begin{eqnarray}
			&& \norm{D_{DG}(u_{h}^{\ast}; g_{D}) - \nabla_{h} u_{h}^{\ast}}_{L^{p}(\Omega)} 
			\lesssim \norm{\jump{u_{h}^{\ast} - u}}_{L^{p}(\Gamma^{o} \cup \Gamma^{D}, h_{e}^{1-p})}, \\
			&& \norm{\nabla_{h} u_{h}^{\ast} - \nabla u}_{L^{p}(\Omega)} 
			\lesssim h^{l} \abs{u}_{W^{l + 1,p}(\Omega)}, \\
			&& \norm{\jump{u_{h}^{\ast} - u}}_{L^{p}(\Gamma^{o} \cup \Gamma^{D}, h_{e}^{1-p})} 
			\lesssim h^{l} \abs{u}_{W^{l+1,p}(\Omega)}. 
		\end{eqnarray}
		Using these results we then have the energy estimates for \(u_{h}^{\ast}\) by 
		\begin{equation}
			E_{h}(u_{h}^{\ast}; g_{D}) 
			\leqslant \norm{\nabla u}_{L^{p}(\Omega)} + C h^{l} \abs{u}_{W^{l+1, p}(\Omega)}. 
		\end{equation}
	\end{subequations}

	\begin{proof}
		Just apply Lemma \ref{lemma: Consistency of the DG-gradient operator}, \ref{thm: polynomial approximation} and \ref{thm: Multiplicative trace inequality}. 
	\end{proof}
\end{lemma}

\begin{theorem}[\emph{A priori} error estimates for the primal variable]\label{our error estimates for primal variable}
	Let \(k \in \N\) be the minimal polynomial degree used in the LDG scheme \eqref{p-Laplace LDG primal weak form 2} or \eqref{p-Laplace LDG minimization form 2} and \(u_{h}\) be the discrete solution. Let \(u\) be the solution to the \(p\)-Laplace equation \eqref{p-Laplace primal weak form}, and let \(\bm q := \nabla u\) and \(\bm\sigma := \mathcal{A}(\bm q)\). Assume \((\bm q, \bm\sigma, u) \in W^{s, p}(\Omega) \times (W^{r, p'}(\Omega) \cap L^{p'}(\mathrm{div}, \Omega)) \times W^{s+1, p}(\Omega)\), where WLOG we assume \(s, r \in \N\) satisfy \(s \leqslant k\) and \(r \leqslant k+1\). Let \(u_{h}^{\ast} := \Pi_{V_{h}} u\), then we have the following estimates. 
	\begin{subequations}
		\begin{eqnarray}
			&& \begin{aligned}
				& \text{For \(p \in (1, 2]\):} \\
				& \norm{u_{h} - u_{h}^{\ast}}_{J, p} 
				\lesssim C_{g_{D}, f, \bm g_{N}\cdot\bm n}^{2-p} 
				\left(
					h^{s(p-1)} \abs{u}_{W^{s+1,p}(\Omega)}^{p-1} 
					+ h^{r} \abs{\bm\sigma}_{W^{r, p'}(\Omega)} 
				\right). 
			\end{aligned} \\
			&& \begin{aligned}
				& \text{For \(p \in [2, +\infty)\):} \\
				& \norm{u_{h} - u_{h}^{\ast}}_{J, p} 
				\lesssim h^{\frac{s}{p-1}} C_{g_{D}, f, \bm g_{N}\cdot\bm n}^{\frac{p-2}{p-1}} \abs{u}_{W^{s+1,p}(\Omega)}^{\frac{1}{p-1}} 
                + h^{\frac{r}{p-1}} \abs{\bm\sigma}_{W^{r, p'}(\Omega)}^{\frac{1}{p-1}} 
				+ h^{s} \abs{u}_{W^{s+1,p}(\Omega)}. 
			\end{aligned}
		\end{eqnarray}
	\end{subequations}
	Here, the hidden constant \(C > 0\) is independent of \(h\), and \(C_{g_{D}, f, \bm g_{N}\cdot\bm n}\) is the right-hand side of the energy estimates in Lemma \ref{p-Laplace energy estimate} and \ref{p-Laplace LDG2 energy-norm estimate}, i.e. 
	\begin{equation}
		C_{g_{D}, f, \bm g_{N}\cdot\bm n} 
		:= \norm{g_{D}}_{W^{\frac{1}{p'}, p}(\Gamma^{D})} + \norm{f}_{L^{p'}(\Omega)}^{\frac{1}{p-1}} + \norm{\bm g_{N} \cdot \bm n}_{L^{p'}(\Gamma^{N})}^{\frac{1}{p-1}}. 
	\end{equation}

	\begin{proof}
		Similar to the \(u_{h}^{\ast}\), we define the projected gradient variables \(\bm q_{h}^{\ast} := \Pi_{Q_{h}} \bm q\) and \(\bm\sigma_{h}^{\ast} := \Pi_{\Sigma_{h}} \bm\sigma\). Besides, we define the discrete error as \(e_{h} := u_{h} - u_{h}^{\ast}\). 

		\begin{paragraph}{The Core Inequality.}
			Since \(u_{h}\) minimizes \(J_{h}\), the Newton-Leibniz formula gives 
			\begin{align*}
				0 
				\geqslant & J_{h}(u_{h}) - J_{h}(u_{h}^{\ast}) \\
				= & \int_{0}^{1} J_{h}'(u_{h}^{\ast} + t e_{h})(e_{h}) \d t \\
				= & J_{h}'(u_{h}^{\ast})(e_{h}) + \int_{0}^{1} \left(J_{h}'(u_{h}^{\ast} + t e_{h})(e_{h}) - J_{h}'(u_{h}^{\ast})(e_{h})\right) \d t, 
			\end{align*}
			yielding our core inequality 
			\begin{equation}\label{the core inequality}\tag{\(\ast\)}
				J_{h}'(u_{h}^{\ast})(-e_{h}) 
				\geqslant \int_{0}^{1} \left(J_{h}'(u_{h}^{\ast} + t e_{h})(e_{h}) - J_{h}'(u_{h}^{\ast})(e_{h})\right) \d t. 
			\end{equation}
			In the following part, we will refer to the left- and right-hand side of the inequality \eqref{the core inequality} as LHS\eqref{the core inequality} and RHS\eqref{the core inequality} respectively. 
		\end{paragraph}

		\begin{paragraph}{RHS\eqref{the core inequality} Estimates with Quasi-Norms.}
			Notice that for \(t \in [0, 1]\), 
			\begin{align*}
				\frac{t}{t+1} \left(\abs{a} + \abs{b}\right) 
				\leqslant \abs{a} + \abs{a \pm t b} 
				\leqslant 2 \left(\abs{a} + \abs{b}\right). 
			\end{align*}
			Use this inequality and Lemma \ref{estimate for the A-operator}, for any fixed \(\delta \geqslant 0\), 
			\begin{align*}
				\text{RHS\eqref{the core inequality}} 
				= & \int_{0}^{1} \left(\mathcal{A}(D_{DG}(u_{h}^{\ast} + t e_{h}; g_{D})) - \mathcal{A}(D_{DG}(u_{h}^{\ast}; g_{D})), D_{DG}(t e_{h}; 0)\right)_{\Omega} t^{-1} \d t \\& + \int_{0}^{1} \left\langle \eta\left(\mathcal{A}(h_{e}^{-1} \jump{u_{h}^{\ast} + t e_{h} - u}) - \mathcal{A}(h_{e}^{-1} \jump{u_{h}^{\ast} - u})\right), \jump{t e_{h}} \right\rangle_{\Gamma^{o} \cup \Gamma^{D}} t^{-1} \d t \\
				\geqslant & C_{\delta} \int_{0}^{1} \left(\abs{D_{DG}(e_{h}; 0)}^{2 + \delta}, \left(\abs{D_{DG}(u_{h}^{\ast} + t e_{h}; g_{D})} + \abs{D_{DG}(u_{h}^{\ast}; g_{D})}\right)^{p - 2 - \delta} \right)_{\Omega} t^{1 + \delta} \d t \\& + C_{\delta} \int_{0}^{1} \left\langle \eta h_{e}^{1-p} \abs{\jump{e_{h}}}^{2 + \delta}, \left(\abs{\jump{u_{h}^{\ast} + t e_{h} - u}} + \abs{\jump{u_{h}^{\ast} - u}}\right)^{p - 2 - \delta} \right\rangle_{\Gamma^{o} \cup \Gamma^{D}} t^{1 + \delta} \d t \\
				\geqslant & C_{p,\delta} \left( \abs{D_{DG}(e_{h}; 0)}_{D_{DG}(u_{h}^{\ast}; g_{D}), p, \delta + 2, \Omega} + \abs{\jump{e_{h}}}_{\jump{u_{h}^{\ast} - u}, p, \delta + 2, \Gamma^{o} \cup \Gamma^{D}, \eta h_{e}^{1-p}} \right). 
			\end{align*}
		\end{paragraph}

		\begin{paragraph}{Partial Estimates for LHS\eqref{the core inequality}.}
			Substitute the equation that \(\bm\sigma\) satisfies in \eqref{p-Laplace discrete equation satisfied by the exact solution} into LHS\eqref{the core inequality} to obtain 
			\begin{align*}
				\text{LHS\eqref{the core inequality}} 
				= & -\left(\mathcal{A}(D_{DG}(u_{h}^{\ast}; g_{D})), D_{DG}(e_{h}; 0)\right)_{\Omega} - \left\langle \eta \mathcal{A}(h_{e}^{-1} \jump{u_{h}^{\ast}}), \jump{e_{h}} \right\rangle_{\Gamma^{o}} \\& - \left\langle \eta \mathcal{A}(h_{e}^{-1} (u_{h}^{\ast} - g_{D}) \bm n), e_{h} \bm n \right\rangle_{\Gamma^{D}} + \left(f, e_{h}\right)_{\Omega} + \left\langle \bm g_{N} \cdot \bm n, e_{h} \right\rangle_{\Gamma^{N}} \\
				= & \left(\bm\sigma - \mathcal{A}(D_{DG}(u_{h}^{\ast}; g_{D})), D_{DG}(e_{h}; 0)\right)_{\Omega} \\& + \left\langle \avg{\bm\sigma_{h}^{\ast}} - \bm C_{12} \jump{\bm\sigma_{h}^{\ast}} - \bm\sigma, \jump{e_{h}} \right\rangle_{\Gamma^{o} \cup \Gamma^{D}} \\& - \left\langle \eta \mathcal{A}(h_{e}^{-1} \jump{u_{h}^{\ast} - u}), \jump{e_{h}} \right\rangle_{\Gamma^{o} \cup \Gamma^{D}} \\
				:= & T_{1} + T_{2} + T_{3}. 
			\end{align*}
			For \(T_{2}\) and \(T_{3}\), we have the following estimates for \(p \in (1, +\infty)\). 
			\begin{align*}
				T_{2} 
				= & \left\langle \avg{\bm\sigma_{h}^{\ast}} - \bm C_{12} \jump{\bm\sigma_{h}^{\ast}} - \bm\sigma, \jump{e_{h}} \right\rangle_{\Gamma^{o} \cup \Gamma^{D}} \\
				\lesssim & \norm{(\bm\sigma_{h}^{\ast} - \bm\sigma) \cdot \bm n}_{L^{p'}(\Gamma^{o} \cup \Gamma^{D})} \norm{\jump{e_{h}}}_{L^{p}(\Gamma^{o} \cup \Gamma^{D}, \eta)} \\
				= & \norm{h_{e}^{1-\frac{1}{p}} (\bm\sigma_{h}^{\ast} - \bm\sigma) \cdot \bm n}_{L^{p'}(\Gamma^{o} \cup \Gamma^{D})} \norm{\jump{e_{h}}}_{L^{p}(\Gamma^{o} \cup \Gamma^{D}, \eta h_{e}^{1-p})} \\
				\lesssim & h^{r} \abs{\bm\sigma}_{W^{r, p'}(\Omega)} \norm{\jump{e_{h}}}_{L^{p}(\Gamma^{o} \cup \Gamma^{D}, \eta h_{e}^{1-p})}, \tag{by Lemma \ref{thm: polynomial approximation} and \ref{thm: Multiplicative trace inequality}} 
			\end{align*}
			and 
			\begin{align*}
				T_{3} 
				= & -\left\langle \eta \mathcal{A}(h_{e}^{-1} \jump{u_{h}^{\ast} - u}), \jump{e_{h}} \right\rangle_{\Gamma^{o} \cup \Gamma^{D}} \\
				\leqslant & \norm{ \eta^{-\frac{1}{p}} \mathcal{A}( h_{e}^{-1} (u_{h}^{\ast} - u) \bm n) }_{L^{p'}(\Gamma^{o} \cup \Gamma^{D})} \norm{\jump{e_{h}}}_{L^{p}(\Gamma^{o} \cup \Gamma^{D}, \eta)} \\
				\lesssim & \norm{ \mathcal{A}( h_{e}^{\frac{1}{p} - 1} (u_{h}^{\ast} - u) \bm n) }_{L^{p'}(\Gamma^{o} \cup \Gamma^{D})} \norm{\jump{e_{h}}}_{L^{p}(\Gamma^{o} \cup \Gamma^{D}, \eta h_{e}^{1-p})} \\
				= & \norm{ h_{e}^{\frac{1}{p} - 1} (u_{h}^{\ast} - u) \bm n}_{L^{p}(\Gamma^{o} \cup \Gamma^{D})}^{p-1} \norm{\jump{e_{h}}}_{L^{p}(\Gamma^{o} \cup \Gamma^{D}, \eta h_{e}^{1-p})} \\
				\lesssim & h^{s(p-1)} \abs{u}_{W^{s+1,p}(\Omega)}^{p-1} \norm{\jump{e_{h}}}_{L^{p}(\Gamma^{o} \cup \Gamma^{D}, \eta h_{e}^{1-p})}. \tag{by Lemma \ref{thm: polynomial approximation} and \ref{thm: Multiplicative trace inequality}} 
			\end{align*}
		\end{paragraph}

		Now, we provide the remaining estimates for cases \(p \in (1, 2]\) and \(p \in [2, +\infty)\) separately. 

		\begin{paragraph}{Case 1: \(p \in (1, 2]\).}
			We choose \(\delta = 0\) in RHS\eqref{the core inequality} in this case. Pick \(\sigma = \delta + 2 = 2\) in Lemma \ref{lemma: paper func inequality} to get 
			\begin{eqnarray*}
				&& \begin{aligned}
					& \norm{D_{DG}(e_{h}; 0)}_{L^{p}(\Omega)}^{2} \\
					\lesssim & \left(\norm{D_{DG}(e_{h}; 0)}_{L^{p}(\Omega)} + \norm{D_{DG}(u_{h}^{\ast}; g_{D})}_{L^{p}(\Omega)}\right)^{2-p} \abs{D_{DG}(e_{h}; 0)}_{D_{DG}(u_{h}^{\ast}; g_{D}), p, 2, \Omega}, 
				\end{aligned} \\
				&& \begin{aligned}
					& \norm{\jump{e_{h}}}_{L^{p}(\Gamma^{o} \cup \Gamma^{D}, \eta h_{e}^{1-p})}^{2} \\
					\lesssim & \left(\norm{\jump{e_{h}}}_{L^{p}(\Gamma^{o} \cup \Gamma^{D}, \eta h_{e}^{1-p})} + \norm{\jump{u_{h}^{\ast} - u}}_{L^{p}(\Gamma^{o} \cup \Gamma^{D}, \eta h_{e}^{1-p})}\right)^{2-p} \abs{\jump{e_{h}}}_{\jump{u_{h}^{\ast} - u}, p, 2, \Gamma^{o} \cup \Gamma^{D}, \eta h_{e}^{1-p}}, 
				\end{aligned}
			\end{eqnarray*}
			and 
			\begin{align*}
				& \norm{D_{DG}(e_{h}; 0)}_{L^{p}(\Omega)} + \norm{\jump{e_{h}}}_{L^{p}(\Gamma^{o} \cup \Gamma^{D}, \eta h_{e}^{1-p})} \\
				\leqslant & \abs{\norm{D_{DG}(u_{h}; g_{D})}_{L^{p}(\Omega)} - \norm{D_{DG}(u_{h}^{\ast}; g_{D})}_{L^{p}(\Omega)}} \\& + \abs{\norm{\jump{u_{h} - u}}_{L^{p}(\Gamma^{o} \cup \Gamma^{D}, \eta h_{e}^{1-p})} - \norm{\jump{u_{h}^{\ast} - u}}_{L^{p}(\Gamma^{o} \cup \Gamma^{D}, \eta h_{e}^{1-p})}} \\
				\leqslant & \norm{D_{DG}(u_{h}; g_{D})}_{L^{p}(\Omega)} + \norm{D_{DG}(u_{h}^{\ast}; g_{D})}_{L^{p}(\Omega)} \\& + \norm{\jump{u_{h} - u}}_{L^{p}(\Gamma^{o} \cup \Gamma^{D}, \eta h_{e}^{1-p})} + \norm{\jump{u_{h}^{\ast} - u}}_{L^{p}(\Gamma^{o} \cup \Gamma^{D}, \eta h_{e}^{1-p})} \\ 
				\leqslant & 2^{1-\frac{1}{p}} \left( E_{h}(u_{h}; g_{D}) + E_{h}(u_{h}^{\ast}; g_{D}) \right) \\
				\lesssim & C_{g_{D}, f, \bm g_{N}\cdot\bm n}, \tag{by Proposition \ref{p-Laplace LDG2 energy-norm estimate}, Lemma \ref{projected solution estimate}, and Proposition \ref{p-Laplace energy estimate}} 
			\end{align*}
			hence 
			\begin{eqnarray*}
				&& \norm{D_{DG}(e_{h}; 0)}_{L^{p}(\Omega)}^{2} 
				\lesssim C_{g_{D}, f, \bm g_{N}\cdot\bm n}^{2-p} \abs{D_{DG}(e_{h}; 0)}_{D_{DG}(u_{h}^{\ast}; g_{D}), p, 2, \Omega}, \\
				&& \norm{\jump{e_{h}}}_{L^{p}(\Gamma^{o} \cup \Gamma^{D}, \eta h_{e}^{1-p})}^{2} 
				\lesssim C_{g_{D}, f, \bm g_{N}\cdot\bm n}^{2-p} \abs{\jump{e_{h}}}_{\jump{u_{h}^{\ast} - u}, p, 2, \Gamma^{o} \cup \Gamma^{D}, \eta h_{e}^{1-p}}. 
			\end{eqnarray*}
			For \(T_{1}\), we have 
			\begin{align*}
				T_{1} 
				= & \left(\bm\sigma - \mathcal{A}(D_{DG}(u_{h}^{\ast}; g_{D})), D_{DG}(e_{h}; 0)\right)_{\Omega} \\
				\leqslant & \norm{\mathcal{A}(\bm q) - \mathcal{A}(D_{DG}(u_{h}^{\ast}; g_{D}))}_{L^{p'}(\Omega)} \norm{D_{DG}(e_{h}; 0)}_{L^{p}(\Omega)} \\
				\lesssim & \norm{ \abs{\bm q - D_{DG}(u_{h}^{\ast}; g_{D})}^{p-1} }_{L^{p'}(\Omega)} \norm{D_{DG}(e_{h}; 0)}_{L^{p}(\Omega)} \tag{by Lemma \ref{estimate for the A-operator} with \(\delta = 2-p\)} \\
				= & \norm{\nabla u - D_{DG}(u_{h}^{\ast}; g_{D})}_{L^{p}(\Omega)}^{p-1} \norm{D_{DG}(e_{h}; 0)}_{L^{p}(\Omega)} \\
				\lesssim & h^{s(p-1)} \abs{u}_{W^{s+1,p}(\Omega)}^{p-1} \norm{D_{DG}(e_{h}; 0)}_{L^{p}(\Omega)}. \tag{by Lemma \ref{projected solution estimate} and \ref{thm: polynomial approximation}} 
			\end{align*}
			Substitute all these estimates into the core inequality \eqref{the core inequality} to obtain 
			\begin{align*}
				\norm{e_{h}}_{J, p}^{2} 
				\leqslant & 2^{\frac{2}{p} - 1} \left( \norm{D_{DG}(e_{h}; 0)}_{L^{p}(\Omega)}^{2} + \norm{\jump{e_{h}}}_{L^{p}(\Gamma^{o} \cup \Gamma^{D}, \eta h_{e}^{1-p})}^{2} \right) \\
				\lesssim & C_{g_{D}, f, \bm g_{N}\cdot\bm n}^{2-p} \text{RHS\eqref{the core inequality}} \\
				\leqslant & C_{g_{D}, f, \bm g_{N}\cdot\bm n}^{2-p} \text{LHS\eqref{the core inequality}} \\
				\lesssim & C_{g_{D}, f, \bm g_{N}\cdot\bm n}^{2-p} h^{s(p-1)} \abs{u}_{W^{s+1,p}(\Omega)}^{p-1} \norm{D_{DG}(e_{h}; 0)}_{L^{p}(\Omega)} \\& + C_{g_{D}, f, \bm g_{N}\cdot\bm n}^{2-p} h^{r} \abs{\bm\sigma}_{W^{r, p'}(\Omega)} \norm{\jump{e_{h}}}_{L^{p}(\Gamma^{o} \cup \Gamma^{D}, \eta h_{e}^{1-p})} \\& + C_{g_{D}, f, \bm g_{N}\cdot\bm n}^{2-p} h^{s(p-1)} \abs{u}_{W^{s+1,p}(\Omega)}^{p-1} \norm{\jump{e_{h}}}_{L^{p}(\Gamma^{o} \cup \Gamma^{D}, \eta h_{e}^{1-p})} \\
				\lesssim & C_{g_{D}, f, \bm g_{N}\cdot\bm n}^{2-p} 
				\left(
					h^{s(p-1)} \abs{u}_{W^{s+1,p}(\Omega)}^{p-1} 
					+ h^{r} \abs{\bm\sigma}_{W^{r, p'}(\Omega)} 
				\right) 
				\norm{e_{h}}_{J, p}, 
			\end{align*}
			hence 
			\begin{align*}
				\norm{e_{h}}_{J, p} 
				\lesssim C_{g_{D}, f, \bm g_{N}\cdot\bm n}^{2-p} 
				\left(
					h^{s(p-1)} \abs{u}_{W^{s+1,p}(\Omega)}^{p-1} 
					+ h^{r} \abs{\bm\sigma}_{W^{r, p'}(\Omega)} 
				\right). 
			\end{align*}
		\end{paragraph}

		\begin{paragraph}{Case 2: \(p \in [2, +\infty)\).}
			We choose \(\delta = p-2\) in RHS\eqref{the core inequality} in this case. Pick \(\sigma = \delta + 2 = p\) in Lemma \ref{lemma: paper func inequality} to get 
			\begin{eqnarray*}
				&& \norm{D_{DG}(e_{h}; 0)}_{L^{p}(\Omega)}^{p} 
				\lesssim \abs{D_{DG}(e_{h}; 0)}_{D_{DG}(u_{h}^{\ast}; g_{D}), p, p, \Omega}, \\
				&& \norm{\jump{e_{h}}}_{L^{p}(\Gamma^{o}, \eta h_{e}^{1-p})}^{p} 
				\lesssim \abs{\jump{e_{h}}}_{\jump{u_{h}^{\ast} - u}, p, p, \Gamma^{o} \cup \Gamma^{D}, \eta h_{e}^{1-p}}. 
			\end{eqnarray*}
			For \(T_{1}\), we have 
			\begin{align*}
				T_{1} 
				= & \left(\bm\sigma - \mathcal{A}(D_{DG}(u_{h}^{\ast}; g_{D})), D_{DG}(e_{h}; 0)\right)_{\Omega} \\
				\leqslant & \norm{\mathcal{A}(\bm q) - \mathcal{A}(D_{DG}(u_{h}^{\ast}; g_{D}))}_{L^{p'}(\Omega)} \norm{D_{DG}(e_{h}; 0)}_{L^{p}(\Omega)} \\
				\lesssim & \norm{ \abs{\bm q - D_{DG}(u_{h}^{\ast}; g_{D})} \left(\abs{\bm q} + \abs{D_{DG}(u_{h}^{\ast}; g_{D})}\right)^{p-2} }_{L^{p'}(\Omega)} \norm{D_{DG}(e_{h}; 0)}_{L^{p}(\Omega)} \tag{by Lemma \ref{estimate for the A-operator} with \(\delta = 0\)} \\
				\leqslant & \norm{\nabla u - D_{DG}(u_{h}^{\ast}; g_{D})}_{L^{p}(\Omega)} \norm{\abs{\nabla u} + \abs{D_{DG}(u_{h}^{\ast}; g_{D})}}_{L^{p}(\Omega)}^{p-2} \norm{D_{DG}(e_{h}; 0)}_{L^{p}(\Omega)} \\
				\lesssim & \norm{\nabla u - D_{DG}(u_{h}^{\ast}; g_{D})}_{L^{p}(\Omega)} \norm{ 2 \abs{\nabla u} + \abs{\nabla u - D_{DG}(u_{h}^{\ast}; g_{D})}}_{L^{p}(\Omega)}^{p-2} \norm{D_{DG}(e_{h}; 0)}_{L^{p}(\Omega)} \\
				\lesssim & \left(h^{s} \abs{u}_{W^{s+1,p}(\Omega)}\right) \left(C_{g_{D}, f, \bm g_{N}\cdot\bm n} + \abs{u}_{W^{1,p}(\Omega)}\right)^{p-2} \norm{D_{DG}(e_{h}; 0)}_{L^{p}(\Omega)} \tag{by Lemma \ref{projected solution estimate} and \ref{thm: polynomial approximation} and Proposition \ref{p-Laplace energy estimate}} \\
				\lesssim & C_{g_{D}, f, \bm g_{N}\cdot\bm n}^{p-2} h^{s} \abs{u}_{W^{s+1,p}(\Omega)} \norm{D_{DG}(e_{h}; 0)}_{L^{p}(\Omega)}. \tag{by Proposition \ref{p-Laplace energy estimate}} 
			\end{align*}
			Substitute all these estimates into the core inequality \eqref{the core inequality} to obtain 
			\begin{align*}
				\norm{e_{h}}_{J, p}^{p} 
				= & \norm{D_{DG}(e_{h}; 0)}_{L^{p}(\Omega)}^{p} + \norm{\jump{e_{h}}}_{L^{p}(\Gamma^{o} \cup \Gamma^{D}, \eta h_{e}^{1-p})}^{p} \\
				\lesssim & \text{RHS\eqref{the core inequality}} \\
				\leqslant & \text{LHS\eqref{the core inequality}} \\
				\lesssim & C_{g_{D}, f, \bm g_{N}\cdot\bm n}^{p-2} h^{s} \abs{u}_{W^{s+1,p}(\Omega)} \norm{D_{DG}(e_{h}; 0)}_{L^{p}(\Omega)} \\& + h^{r} \abs{\bm\sigma}_{W^{r, p'}(\Omega)} \norm{\jump{e_{h}}}_{L^{p}(\Gamma^{o} \cup \Gamma^{D}, \eta h_{e}^{1-p})} \\& + h^{s(p-1)} \abs{u}_{W^{s+1,p}(\Omega)}^{p-1} \norm{\jump{e_{h}}}_{L^{p}(\Gamma^{o} \cup \Gamma^{D}, \eta h_{e}^{1-p})} \\
				\lesssim & \left( 
					C_{g_{D}, f, \bm g_{N}\cdot\bm n}^{p-2} h^{s} \abs{u}_{W^{s+1,p}(\Omega)} 
					+ h^{r} \abs{\bm\sigma}_{W^{r, p'}(\Omega)} 
					+ h^{s(p-1)} \abs{u}_{W^{s+1,p}(\Omega)}^{p-1} 
				\right) 
				\norm{e_{h}}_{J, p}, 
			\end{align*}
			hence 
			\begin{align*}
				\norm{e_{h}}_{J, p} 
				\lesssim C_{g_{D}, f, \bm g_{N}\cdot\bm n}^{\frac{p-2}{p-1}} h^{\frac{s}{p-1}} \abs{u}_{W^{s+1,p}(\Omega)}^{\frac{1}{p-1}} 
				+ h^{\frac{r}{p-1}} \abs{\bm\sigma}_{W^{r, p'}(\Omega)}^{\frac{1}{p-1}} 
				+ h^{s} \abs{u}_{W^{s+1,p}(\Omega)}. 
			\end{align*}
		\end{paragraph}

		So far, we have completed the proof. 
	\end{proof}
\end{theorem}

\begin{remark}[Comparison to existing HHO and HDG error estimates] 
	Our \emph{a priori} error estimates for primal variable in the mesh-dependent energy norm in Theorem \ref{our error estimates for primal variable} have similar regularity requirements for the exact solution as the HHO estimates \cite[Theorem 3.2]{10.1142/S0218202517500191} and similar results, which is summed up in Table \ref{table: error estimates comparison}, where \(s\) and \(r\) are the regularity parameters in the assumptions of Theorem \ref{our error estimates for primal variable}. Because the \(\mathcal{P}^{k}(K)\otimes \mathcal{P}^{k}(F)\)-HHO methods have more degrees of freedom, higher convergence rate comparable to that of \(\mathcal{P}^{k+1}(K)\)-LDG should be expected, which is also true in the case of linear elliptic PDEs. If we assume the same regularity for \(u\), the estimated convergence rates will be exactly the same if the polynomial orders are high enough. We remark that under additional local regime requirements, the HHO method can have \emph{optimal} convergence rate in the primal variable when \(p \in (1, 2]\) as analyzed in \cite[Theorem 1]{10.1007/s10092-021-00410-z}, and our estimates have not taken these requirements into consideration yet. 
	\begin{table}[htbp]
		\centering
		\caption{Comparison between our LDG estimates and HHO estimates.}
		\label{table: error estimates comparison}
		\begin{tabular}{c | c c c c}
			\hline
			method & polynomial & regularity & error (\(p\in(1, 2]\)) & error (\(p \in [2, +\infty)\)) \\
			\hline
			LDG & \(\mathcal{P}^{k}(K)\) & \(s=r=k\) & \(\mathcal{O}(h^{(p-1)k})\) & \(\mathcal{O}(h^{\frac{k}{p-1}})\) \\
			HHO \cite{10.1090/mcom/3180} & \(\mathcal{P}^{k}(K) \otimes \mathcal{P}^{k}(F)\) & \(s=r=k+1\) & \(\mathcal{O}(h^{(p-1)(k+1)})\) \cite{10.1142/S0218202517500191} & \(\mathcal{O}(h^{\frac{k+1}{p-1}})\) \cite{10.1142/S0218202517500191} \\
			\hline
		\end{tabular}
	\end{table}

	Besides, our estimates are also similar to those of the HDG scheme \cite[Theorem 3.2]{10.1007/s10915-019-00967-6} under the lowest regularity assumption, except for a minor difference caused by the choice of the norm of the penalty terms on faces. 
\end{remark}

\begin{remark}[Comparison to existing LDG error estimates]
	Our estimates are in the \(\norm{\cdot}_{J, p}\)-norm (similar to the \(W_{0}^{1,p}\)-norm), which is not only different but also non-equivalent distance functional from that used in previous LDG estimates \cite[Theorem 4.8]{10.1093/imanum/drt040}. Moreover, only our estimates assume that the exact solution has sufficiently high regularity, and we reveal the potential for high-order accuracy with high-order polynomials. 
\end{remark}

\begin{corollary}[\emph{A priori} error estimates for all variables]\label{our error estimates for all variables}
	Let \(k \in \N\) be the minimal polynomial degree used in the LDG scheme \eqref{p-Laplace LDG weak form 2} and \(\left(\bm q_{h}, \bm\sigma_{h}, u_{h}\right)\) be the discrete solution. Let \(u\) be the solution to the \(p\)-Laplace equation \eqref{p-Laplace primal weak form}, and let \(\bm q := \nabla u\) and \(\bm\sigma := \mathcal{A}(\bm q)\). Assume \((\bm q, \bm\sigma, u) \in W^{s, p}(\Omega) \times (W^{r, p'}(\Omega) \cap L^{p'}(\mathrm{div}, \Omega)) \times W^{s+1, p}(\Omega)\), where WLOG we assume \(s, r \in \N\) satisfy \(s \leqslant k\) and \(r \leqslant k+1\). Let \(u_{h}^{\ast} := \Pi_{V_{h}} u\), then we have the following estimates. 
	\begin{subequations}
		\begin{eqnarray}
			&& \begin{aligned}
				& \text{For \(p \in (1, 2]\):} \\
				& \begin{cases}
					\norm{u_{h} - u}_{L^{p}(\Omega)} 
					\lesssim \norm{u_{h} - u_{h}^{\ast}}_{J, p} + h^{s + 1} \abs{u}_{W^{s+1, p}(\Omega)}, \\
					\norm{\bm q_{h} - \bm q}_{L^{p}(\Omega)} + \norm{\jump{u_{h} - u}}_{L^{p}(\Gamma^{o} \cup \Gamma^{D}, h^{1-p})} 
					\lesssim \norm{u_{h} - u_{h}^{\ast}}_{J, p} + h^{s} \abs{u}_{W^{s+1, p}(\Omega)}, \\
					\norm{\bm\sigma_{h} - \bm\sigma}_{L^{p'}(\Omega)} 
					\lesssim \norm{\bm q_{h} - \bm q}_{L^{p}(\Omega)}^{p-1}, \\
					\norm{u_{h} - u_{h}^{\ast}}_{J, p} 
					\lesssim C_{g_{D}, f, \bm g_{N}\cdot\bm n}^{2-p} 
					\left(
						h^{s(p-1)} \abs{u}_{W^{s+1,p}(\Omega)}^{p-1} 
						+ h^{r} \abs{\bm\sigma}_{W^{r, p'}(\Omega)} 
					\right). 
				\end{cases}
			\end{aligned} \\
			&& \begin{aligned}
				& \text{For \(p \in [2, +\infty)\):} \\
				& \begin{cases}
					\norm{u_{h} - u}_{L^{p}(\Omega)} 
					\lesssim \norm{u_{h} - u_{h}^{\ast}}_{J, p} + h^{s + 1} \abs{u}_{W^{s+1, p}(\Omega)}, \\
					\norm{\bm q_{h} - \bm q}_{L^{p}(\Omega)} + \norm{\jump{u_{h} - u}}_{L^{p}(\Gamma^{o} \cup \Gamma^{D}, h^{1-p})} 
					\lesssim \norm{u_{h} - u_{h}^{\ast}}_{J, p} + h^{s} \abs{u}_{W^{s+1, p}(\Omega)}, \\
					\norm{\bm\sigma_{h} - \bm\sigma}_{L^{p'}(\Omega)} 
					\lesssim C_{g_{D}, f, \bm g_{N}\cdot\bm n}^{p-2} \norm{\bm q_{h} - \bm q}_{L^{p}(\Omega)}, \\
					\norm{u_{h} - u_{h}^{\ast}}_{J, p} 
                    \lesssim h^{\frac{s}{p-1}} C_{g_{D}, f, \bm g_{N}\cdot\bm n}^{\frac{p-2}{p-1}} \abs{u}_{W^{s+1,p}(\Omega)}^{\frac{1}{p-1}} 
                    + h^{\frac{r}{p-1}} \abs{\bm\sigma}_{W^{r, p'}(\Omega)}^{\frac{1}{p-1}} 
                    + h^{s} \abs{u}_{W^{s+1,p}(\Omega)}. 
				\end{cases}
			\end{aligned}
		\end{eqnarray}
	\end{subequations}
\end{corollary}

\begin{remark}
	Our error estimates in Theorem \ref{our error estimates for primal variable} and Corollary \ref{our error estimates for all variables} are estimates of the energy-norm type. They coincide with \(H_{0}^{1}\)-estimates for CG for the Poisson equation when \(p = 2\). In order to enhance our \(L^{p}\)-estimates, it is necessary to involve dual arguments, which warrants further study. 
\end{remark}

\section{Numerical Results}\label{se:nu}


Before carrying out numerical experiments, we first mention several aspects that should be paid special attention to. 

The most critical issue is numerical quadrature rules. Our scheme involves the evaluation of integrals of nonlinear functions on elements and faces, which are discretized via quadrature. The positivity-preservation property of quadrature rules is essential to stability because the fully discretized minimization problem must be convex and bounded from below. As a result, the quadrature rules must have positive weights. For forward computation,  we use the Gauss-Legendre rules \cite{10.1137/140954969, 10.1137/S1064827500379690} with an algebraic degree of exactness \(2k+1\) in 1D. In the 2D case, we choose the fully symmetric and positive Gauss quadrature rules \cite{10.2307/43693493} provided in the open-source software PHG \cite{parallel_hierarchical_grid} with an algebraic degree of exactness at least \(2k\). For computing the errors, we also use the same rules. 

Due to the degenerate or singular nature of the \(p\)-Laplace problem, floating-point errors should also be considered. We utilize the Lobatto collocation points \cite{10.1093/imamat/hxh077} and the Bernstein basis \cite{10.1137/11082539X} in order to improve floating-point stability. 
\begin{figure}[hbp]
	\centering
	\begin{subfigure}{0.32\linewidth}
		\centering
		\includegraphics[width=\linewidth]{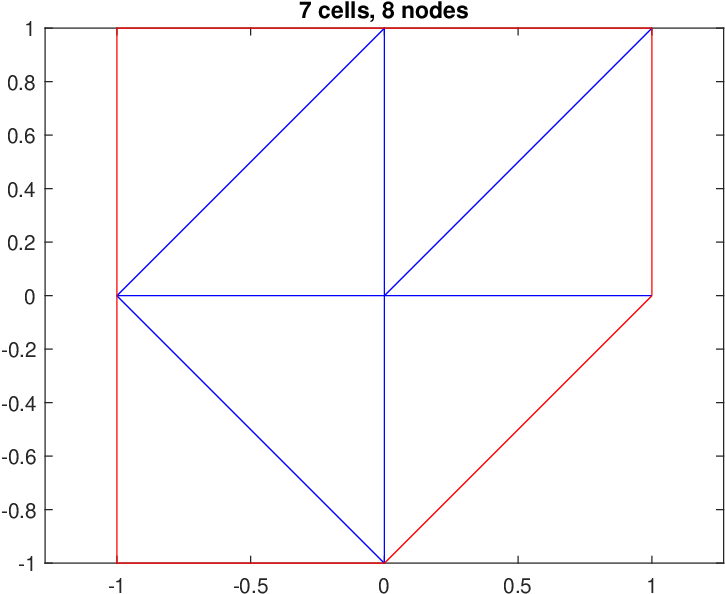}
		\caption{mesh 1}
	\end{subfigure}
	\begin{subfigure}{0.32\linewidth}
		\centering
		\includegraphics[width=\linewidth]{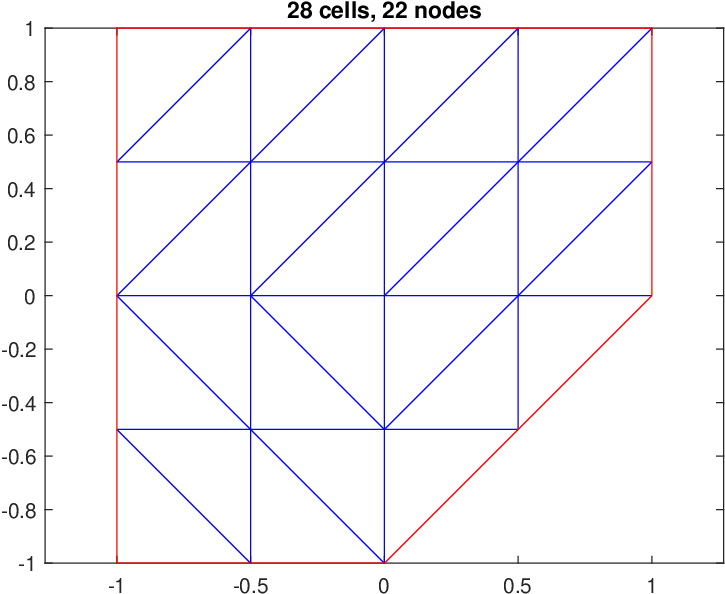}
		\caption{mesh 2}
	\end{subfigure}
	\begin{subfigure}{0.32\linewidth}
		\centering
		\includegraphics[width=\linewidth]{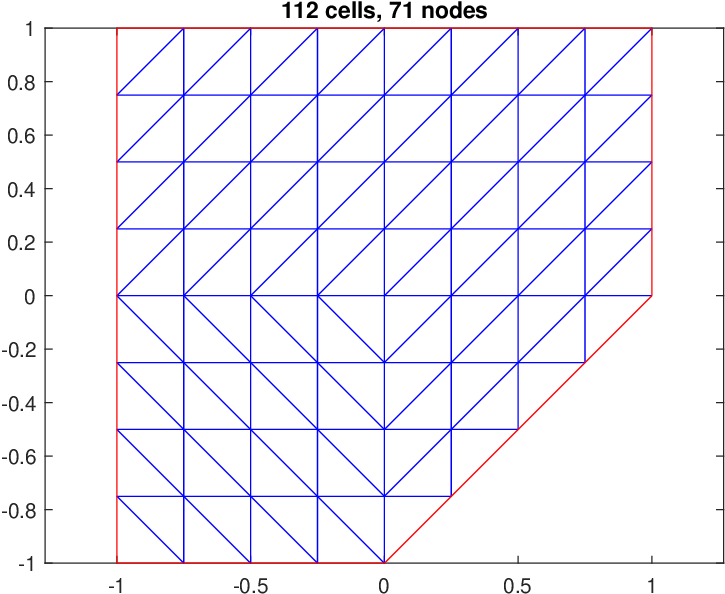}
		\caption{mesh 3}
	\end{subfigure}
	\\
	\begin{subfigure}{0.32\linewidth}
		\centering
		\includegraphics[width=\linewidth]{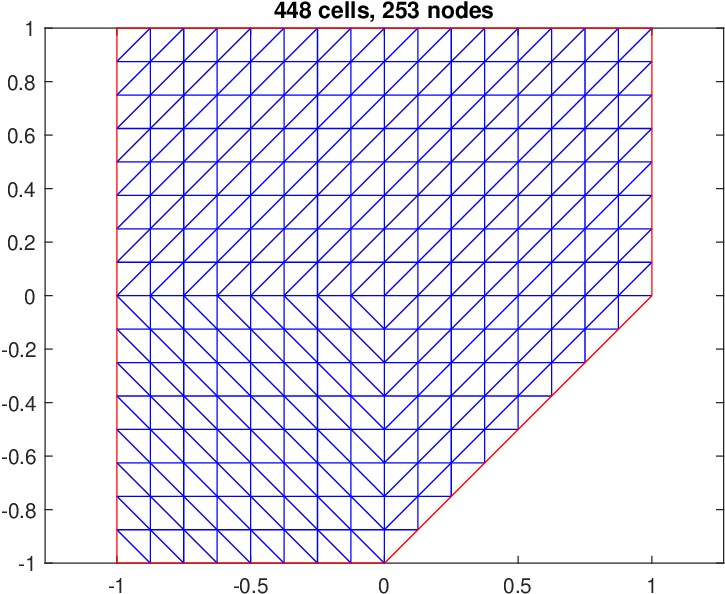}
		\caption{mesh 4}
	\end{subfigure}
	\begin{subfigure}{0.32\linewidth}
		\centering
			\includegraphics[width=\linewidth]{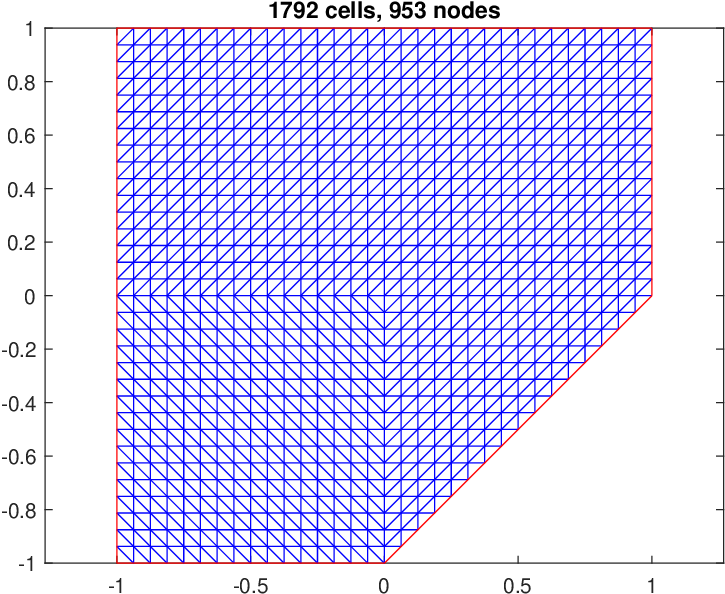}
			\caption{mesh 5}
	\end{subfigure}
	\begin{subfigure}{0.32\linewidth}
		\centering
		\includegraphics[width=\linewidth]{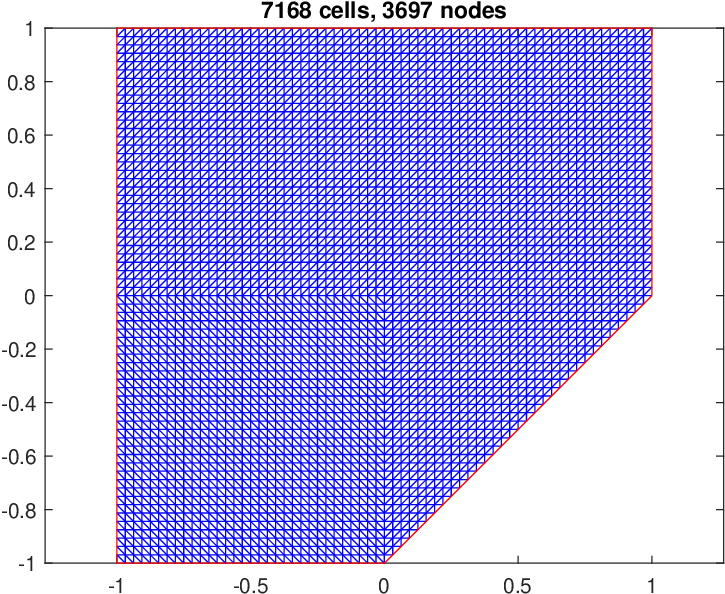}
		\caption{mesh 6}
	\end{subfigure}
	\caption{Uniformly refined conforming triangular meshes for a 2D polygonal domain.}
	\label{Mesh: 2D Dirichlet}
\end{figure}

In the following text, we are going to present our numerical results. For Example \ref{exm:lin} to \ref{exm:de}, we take the polygonal domain \(\Omega = \B{(x,y) \big| -1\leqslant x, y \leqslant 1, y-x+1\geqslant 0}\) and use the uniformly refined conforming triangular meshes as depicted in Figure \ref{Mesh: 2D Dirichlet}. For simplicity, we only use the Dirichlet boundary condition in our tests. The penalty coefficients for the numerical flux \eqref{p-Laplace LDG numerical flux 2} are taken to be \(\eta = 10.0\), and the upwind coefficients \(\bm C_{12}\) are set according to the MD-LDG scheme \cite{10.1007/s10915-007-9130-3}, with the (constant) auxiliary vector field \(\bm v_0 = (1, 0)^{T}\). In the following text, we use \(\epsilon\) to denote any fixed small positive constant.

\begin{example}{\bf A Linear Case}\label{exm:lin}
\end{example}

In our first example, we take the exact solution to be 
\begin{eqnarray*}
	&& u(x, y) = e^{\sin(\pi x)} \cos(\pi(x + y)), \\
	&& \bm q(x, y) = \begin{pmatrix}
		\pi e^{\sin(\pi x)} \left(\cos(\pi x) \cos(\pi(x+y)) - \sin(\pi(x+y))\right) \\
		-\pi e^{\sin(\pi x)} \sin(\pi(x+y)) 
	\end{pmatrix}, \\
	&& \bm \sigma(x, y) = \begin{pmatrix}
		\pi e^{\sin(\pi x)} \left(\cos(\pi x) \cos(\pi(x+y)) - \sin(\pi(x+y))\right) \\
		-\pi e^{\sin(\pi x)} \sin(\pi(x+y)) 
	\end{pmatrix}, \\
	&& f(x, y) = \pi^{2} e^{\sin(\pi x)} \left(\cos(\pi(x+y)) \left(2 - \cos^2(\pi x) + \sin(\pi x)\right) + 2\cos(\pi x) \sin(\pi(x+y))\right). 
\end{eqnarray*}
We test with \(p = 2\) to verify the consistency of our gradient descent method for the Poisson equation. The errors are shown in Table \ref{table: 2D p-Laplace Dirichlet LDG, example 1}, demonstrating that our solver provides the best convergence in this linear case, when compared to linear LDG error estimates \cite{10.1137/S0036142900371003}. We also mention that every test here converges within one iteration, which aligns with Remark \ref{remark: optimal step size} well. 

\begin{table}[htbp]
	\centering
	\caption{Error table for Example \ref{exm:lin}. Parameters: \(p = 2\), \(\epsilon = 0\), \(\delta_{w_{h}} = 10^{-16}\), \(\delta_{\rho} = 10^{-16}\).}
	\label{table: 2D p-Laplace Dirichlet LDG, example 1}
	\begin{tabular}{| c | c c | c c | c c | c c |}
		\hline

		\(k\) & \(N_{e}\) & \(N_{\rm dof}\) & \(\norm{u - u_{h}}_{L^p}\) & order & \(\norm{\bm q - \bm q_{h}}_{L^p}\) & order & \(\norm{\bm\sigma - \bm\sigma_{h}}_{L^{p'}}\) & order \\

		\hline

		1	& 7 & 21 & 1.2653e+00 & - & 1.0916e+01 & - & 1.0916e+01 & - \\
			& 28 & 84 & 7.0497e-01 & 0.8438 & 5.7008e+00 & 0.9372 & 5.7008e+00 & 0.9372 \\
			& 112 & 336 & 2.6383e-01 & 1.4180 & 3.4436e+00 & 0.7272 & 3.4436e+00 & 0.7272 \\
			& 448 & 1344 & 7.5449e-02 & 1.8060 & 1.8233e+00 & 0.9174 & 1.8233e+00 & 0.9174 \\
			& 1792 & 5376 & 1.9627e-02 & 1.9427 & 9.2622e-01 & 0.9771 & 9.2622e-01 & 0.9771 \\
			& 7168 & 21504 & 4.9640e-03 & 1.9832 & 4.6540e-01 & 0.9929 & 4.6540e-01 & 0.9929 \\

		\hline

		2	& 7 & 42 & 6.1560e-01 & - & 5.7034e+00 & - & 5.7034e+00 & - \\
			& 28 & 168 & 1.4129e-01 & 2.1233 & 2.3923e+00 & 1.2534 & 2.3923e+00 & 1.2534 \\
			& 112 & 672 & 2.0681e-02 & 2.7724 & 7.1105e-01 & 1.7504 & 7.1105e-01 & 1.7504 \\
			& 448 & 2688 & 2.7014e-03 & 2.9365 & 1.8887e-01 & 1.9126 & 1.8887e-01 & 1.9126 \\
			& 1792 & 10752 & 3.4226e-04 & 2.9805 & 4.8095e-02 & 1.9734 & 4.8095e-02 & 1.9734 \\
			& 7168 & 43008 & 4.3003e-05 & 2.9926 & 1.2087e-02 & 1.9925 & 1.2087e-02 & 1.9925 \\

		\hline

		3	& 7 & 70 & 4.4329e+00 & - & 4.8816e+00 & - & 4.8816e+00 & - \\
			& 28 & 280 & 3.5407e-02 & 3.6118 & 7.5929e-01 & 2.6846 & 7.5929e-01 & 2.6846 \\
			& 112 & 1120 & 2.7498e-03 & 3.6866 & 1.2018e-01 & 2.6594 & 1.2018e-01 & 2.6594 \\
			& 448 & 4480 & 1.7279e-04 & 3.9922 & 1.5933e-02 & 2.9151 & 1.5933e-02 & 2.9151 \\
			& 1792 & 17920 & 1.0786e-05 & 4.0017 & 2.0237e-03 & 2.9770 & 2.0237e-03 & 2.9770 \\
			& 7168 & 71680 & 6.7332e-07 & 4.0018 & 2.5385e-04 & 2.9950 & 2.5385e-04 & 2.9950 \\

		\hline

		4	& 7 & 105 & 1.2932e-01 & - & 1.5622e+00 & - & 1.5622e+00 & - \\
			& 28 & 420 & 9.2744e-03 & 3.8016 & 2.2166e-01 & 2.8171 & 2.2166e-01 & 2.8171 \\
			& 112 & 1680 & 3.6503e-04 & 4.6672 & 1.7668e-02 & 3.6491 & 1.7668e-02 & 3.6491 \\
			& 448 & 6720 & 1.2501e-05 & 4.8680 & 1.2025e-03 & 3.8770 & 1.2025e-03 & 3.8770 \\
			& 1792 & 26880 & 3.9980e-07 & 4.9666 & 7.6083e-05 & 3.9824 & 7.6083e-05 & 3.9824 \\
			& 7168 & 107520 & 1.2568e-08 & 4.9914 & 4.7602e-06 & 3.9985 & 4.7602e-06 & 3.9985 \\

		\hline

		5	& 7 & 147 & 1.1381e-01 & - & 1.0080e+00 & - & 1.0080e+00 & - \\
			& 28 & 588 & 2.3911e-03 & 5.5728 & 5.5520e-02 & 4.1824 & 5.5520e-02 & 4.1824 \\
			& 112 & 2352 & 5.2502e-05 & 5.5092 & 2.4154e-03 & 4.5227 & 2.4154e-03 & 4.5227 \\
			& 448 & 9408 & 8.6522e-07 & 5.9232 & 7.9605e-05 & 4.9232 & 7.9605e-05 & 4.9232 \\
			& 1792 & 37632 & 1.3790e-08 & 5.9714 & 2.5347e-06 & 4.9730 & 2.5347e-06 & 4.9730 \\
			& 7168 & 150528 & 2.1659e-10 & 5.9925 & 7.9403e-08 & 4.9965 & 7.9403e-08 & 4.9965 \\
		
		\hline

		6   & 7 & 196 & 1.9906e-02 & - & 4.3179e-01 & - & 4.3179e-01 & - \\
			& 28 & 784 & 4.0022e-04 & 5.6363 & 1.3970e-02 & 4.9499 & 1.3970e-02 & 4.9499 \\
			& 112 & 3136 & 4.0023e-06 & 6.6438 & 2.8650e-04 & 5.6077 & 2.8650e-04 & 5.6077 \\
			& 448 & 12544 & 3.5067e-08 & 6.8346 & 5.0342e-06 & 5.8306 & 5.0342e-06 & 5.8306 \\
			& 1792 & 50176 & 2.7932e-10 & 6.9721 & 7.9191e-08 & 5.9903 & 7.9191e-08 & 5.9903 \\
			& 7168 & 200704 & 2.1943e-12 & 6.9920 & 1.2349e-09 & 6.0028 & 1.2349e-09 & 6.0028 \\

		\hline
	\end{tabular}
\end{table}

\begin{example}{\bf A Regular Case}\label{exm:re}
\end{example}

We take the radial exact solution \cite{10.1090/S0025-5718-1993-1192966-4, 10.1137/15M1008014} as 
\begin{eqnarray*}
	&& u(x,y) = \frac{p-1}{(\sigma + 2)^{\frac{1}{p-1}}} \frac{1 - r(x,y)^{\frac{\sigma + p}{p-1}}}{\sigma + p} \in W^{\frac{\sigma+1}{p-1}+\frac{2}{p}+1 - \varepsilon, p}, \\
	&& \bm q(x,y) = \begin{pmatrix}
		-\frac{r(x,y)^{\frac{\sigma+1}{p-1}}}{(\sigma + 2)^{\frac{1}{p-1}}} \frac{x}{r(x,y)} \\
		-\frac{r(x,y)^{\frac{\sigma+1}{p-1}}}{(\sigma + 2)^{\frac{1}{p-1}}} \frac{y}{r(x,y)} 
	\end{pmatrix} \in W^{\frac{\sigma+1}{p-1}+\frac{2}{p} - \varepsilon, p}, \\
	&& \bm\sigma(x,y) = \begin{pmatrix}
		\frac{-r(x,y)^{\sigma} x}{\sigma + 2} \\
		\frac{-r(x,y)^{\sigma} y}{\sigma + 2} 
	\end{pmatrix} \in W^{\sigma + 3 - \frac{2}{p} - \varepsilon, p'}, \\
	&& f(x,y) = r(x,y)^{\sigma}, 
\end{eqnarray*}
where \(r(x,y) = \sqrt{x^2 + y^2}\) is the radial distance. Although the exact solution is not smooth, it does not have singular points or zero-gradient regions with positive measures. For tests, we choose the following two groups of parameters, \((\sigma, p) = (0, 1.5)\) and \((\sigma, p) = (7, 4)\). The errors are listed in Table \ref{table: 2D p-Laplace Dirichlet LDG, example 2, test1} and \ref{table: 2D p-Laplace Dirichlet LDG, example 2, test2} respectively, and some of the convergence history of the nonlinear solver are shown in Figure \ref{pic: 2D p-Laplace Dirichlet LDG, example 2, test1} and \ref{pic: 2D p-Laplace Dirichlet LDG, example 2, test2}. 

The convergence rates are going up only some of the way up as \(k\) increases, possibly due to limited regularity. 
\begin{enumerate}
	\item Case \((\sigma, p) = (0, 1.5)\): \((u, \bm q, \bm\sigma) \in W^{s + 1 - \varepsilon, p} \times W^{s - \varepsilon, p} \times W^{r - \varepsilon, p'}\) with \((s, r) = (\frac{10}{3}, \frac{5}{3})\). Our scheme attains the best convergence rate for all variables compared to the linear LDG scheme, which coincides with observations from the HDG scheme \cite{10.1137/15M1008014}, although in different norms. Moreover, our numerical results coincide with the improved HHO estimates \cite[Theorem 1]{10.1007/s10092-021-00410-z}, so we guess this is due to good local regimes of the exact solution. 
	\item Case \((\sigma, p) = (7, 4)\): \((u, \bm q, \bm\sigma) \in W^{s + 1 - \varepsilon, p} \times W^{s - \varepsilon, p} \times W^{r - \varepsilon, p'}\) with \((s, r) = (\frac{19}{6}, \frac{19}{2})\). The convergence rates for gradient variables \(\bm q\) and \(\bm\sigma\) are the best, but that for the primal variable \(u\) is not. 
\end{enumerate}
Besides, for any \(h\) and \(k\), both examples consistently require around 15 and 25 iterations respectively, demonstrating \(hk\)-independent convergence of our preconditioned method. 

\begin{table}[htbp]
	\centering
	\caption{Error table for Example \ref{exm:re}. Parameters: \(\sigma = 0\), \(p = 1.5\), \(\epsilon = 10^{-14}\), \(\delta_{w_{h}} = 10^{-16}\), \(\delta_{\rho} = 10^{-16}\).}
	\label{table: 2D p-Laplace Dirichlet LDG, example 2, test1}
	\begin{tabular}{| c | c c | c c | c c | c c |}
		\hline

		\(k\) & \(N_{e}\) & \(N_{\rm dof}\) & \(\norm{u - u_{h}}_{L^p}\) & order & \(\norm{\bm q - \bm q_{h}}_{L^p}\) & order & \(\norm{\bm\sigma - \bm\sigma_{h}}_{L^{p'}}\) & order \\
		
		\hline

		1	& 7 & 21 & 4.5425e-02 & - & 2.2627e-01 & - & 3.2533e-01 & - \\
			& 28 & 84 & 1.2729e-02 & 1.8354 & 1.2465e-01 & 0.8601 & 1.6710e-01 & 0.9612 \\
			& 112 & 336 & 3.2486e-03 & 1.9702 & 6.4513e-02 & 0.9503 & 8.4892e-02 & 0.9770 \\
			& 448 & 1344 & 8.1708e-04 & 1.9913 & 3.2663e-02 & 0.9820 & 4.2871e-02 & 0.9856 \\
			& 1792 & 5376 & 2.0453e-04 & 1.9982 & 1.6402e-02 & 0.9938 & 2.1543e-02 & 0.9928 \\
			& 7168 & 21504 & 5.1147e-05 & 1.9996 & 8.2140e-03 & 0.9977 & 1.0796e-02 & 0.9967 \\

		\hline

		2	& 7 & 42 & 7.7338e-03 & - & 5.2735e-02 & - & 9.2414e-02 & - \\
			& 28 & 168 & 8.9019e-04 & 3.1190 & 1.3935e-02 & 1.9200 & 2.9935e-02 & 1.6263 \\
			& 112 & 672 & 1.0777e-04 & 3.0461 & 3.6288e-03 & 1.9412 & 9.7520e-03 & 1.6181 \\
			& 448 & 2688 & 1.3320e-05 & 3.0164 & 9.2595e-04 & 1.9705 & 3.1265e-03 & 1.6412 \\
			& 1792 & 10752 & 1.6624e-06 & 3.0022 & 2.3373e-04 & 1.9861 & 9.9332e-04 & 1.6542 \\
			& 7168 & 43008 & 2.0793e-07 & 2.9990 & 5.8705e-05 & 1.9933 & 3.1419e-04 & 1.6606 \\

		\hline

		3	& 7 & 70 & 6.7287e-04 & - & 9.0177e-03 & - & 1.9052e-02 & - \\
			& 28 & 280 & 5.5204e-05 & 3.6075 & 1.5156e-03 & 2.5729 & 6.4762e-03 & 1.5567 \\
			& 112 & 1120 & 3.9437e-06 & 3.8071 & 2.1125e-04 & 2.8428 & 2.0421e-03 & 1.6651 \\
			& 448 & 4480 & 2.6731e-07 & 3.8829 & 2.8250e-05 & 2.9026 & 6.4327e-04 & 1.6665 \\
			& 1792 & 17920 & 1.7606e-08 & 3.9244 & 3.6875e-06 & 2.9375 & 2.0247e-04 & 1.6677 \\
			& 7168 & 71680 & 1.8102e-09 & 3.2818 & 4.7462e-07 & 2.9578 & 6.3822e-05 & 1.6656 \\

		\hline

		4	& 7 & 105 & 1.3377e-04 & - & 1.2196e-03 & - & 3.0487e-03 & - \\
			& 28 & 420 & 7.9208e-06 & 4.0780 & 1.6633e-04 & 2.8743 & 9.5295e-04 & 1.6777 \\
			& 112 & 1680 & 4.2657e-07 & 4.2148 & 1.8647e-05 & 3.1570 & 3.0028e-04 & 1.6661 \\
			& 448 & 6720 & 2.1984e-08 & 4.2783 & 1.9546e-06 & 3.2540 & 9.4582e-05 & 1.6667 \\
			& 1792 & 26880 & 1.2093e-09 & 4.1842 & 1.9903e-07 & 3.2958 & 2.9790e-05 & 1.6667 \\
			& 7168 & 107520 & 1.1226e-09 & 0.1074 & 2.1021e-08 & 3.2431 & 9.3838e-06 & 1.6666 \\

		\hline

		5	& 7 & 147 & 2.9285e-05 & - & 4.6112e-04 & - & 1.5947e-03 & - \\
			& 28 & 588 & 1.5832e-06 & 4.2093 & 5.1172e-05 & 3.1717 & 5.0520e-04 & 1.6584 \\
			& 112 & 2352 & 7.9726e-08 & 4.3116 & 5.1633e-06 & 3.3090 & 1.5913e-04 & 1.6666 \\
			& 448 & 9408 & 5.6283e-09 & 3.8243 & 5.1460e-07 & 3.3268 & 5.0123e-05 & 1.6667 \\
			& 1792 & 37632 & 5.3178e-09 & 0.0819 & 5.7037e-08 & 3.1735 & 1.5788e-05 & 1.6666 \\
			& 7168 & 150528 & 3.7110e-10 & 3.8410 & 5.4414e-09 & 3.3899 & 4.9727e-06 & 1.6668 \\
		
		\hline

		6   & 7 & 196 & 6.8486e-06 & - & 1.7974e-04 & - & 6.7422e-04 & - \\
			& 28 & 784 & 3.5864e-07 & 4.2552 & 1.8790e-05 & 3.2579 & 2.1222e-04 & 1.6676 \\
			& 112 & 3136 & 1.7881e-08 & 4.3260 & 1.8690e-06 & 3.3296 & 6.6846e-05 & 1.6667 \\
			& 448 & 12544 & 2.1894e-09 & 3.0299 & 1.8631e-07 & 3.3265 & 2.1056e-05 & 1.6666 \\
			& 1792 & 50176 & 2.0939e-09 & 0.0644 & 2.0968e-08 & 3.1514 & 6.6326e-06 & 1.6666 \\
			& 7168 & 200704 & 1.9309e-09 & 0.1169 & 6.8024e-09 & 1.6241 & 2.0898e-06 & 1.6662 \\

		\hline
	\end{tabular}
\end{table}

\begin{figure}[htbp]
	\centering
	\begin{subfigure}{0.32\linewidth}
		\centering
		\includegraphics[width=\linewidth]{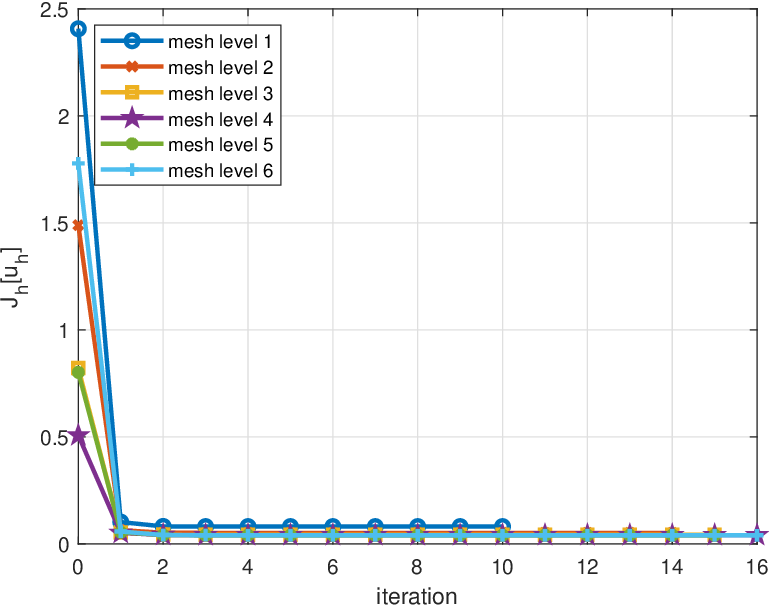}
		\caption{\(J_{h}(u_{h})\), \(k = 1\)}
	\end{subfigure}
	\begin{subfigure}{0.32\linewidth}
		\centering
		\includegraphics[width=\linewidth]{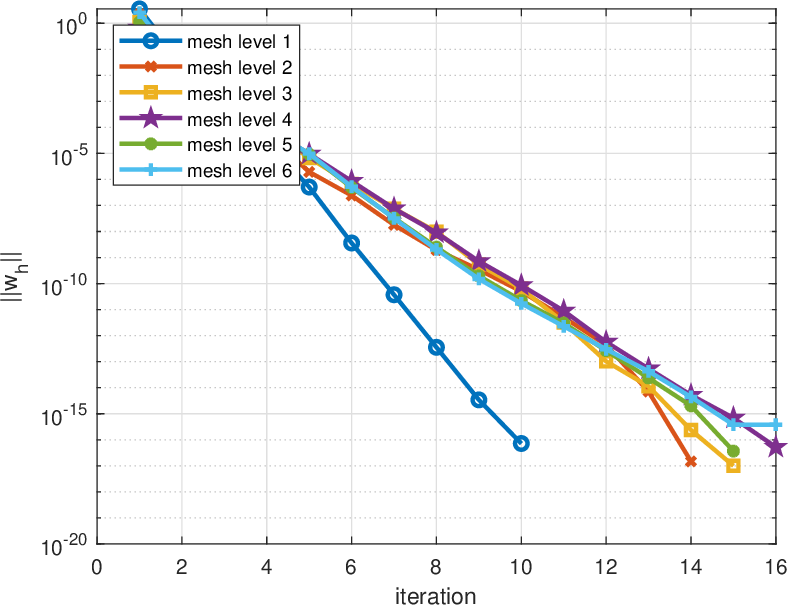}
		\caption{\(\norm{w_{h}}\), \(k = 1\)}
	\end{subfigure}
	\begin{subfigure}{0.32\linewidth}
		\centering
		\includegraphics[width=\linewidth]{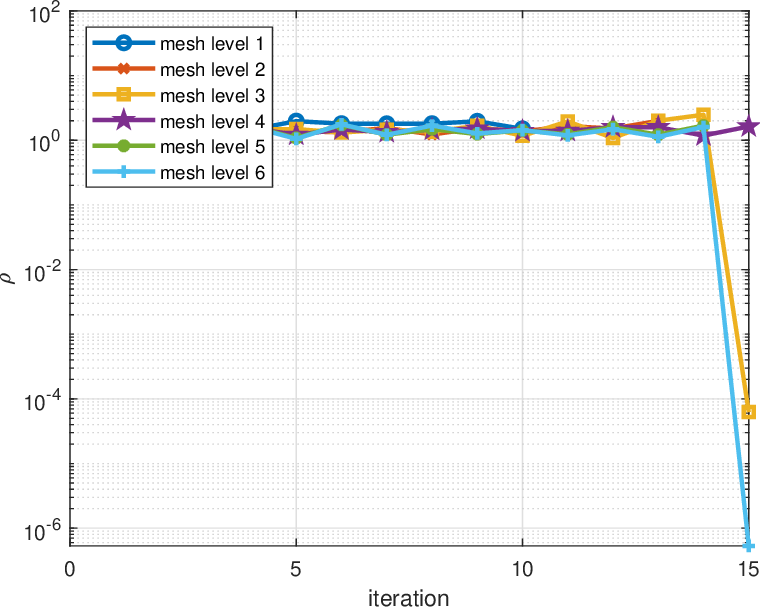}
		\caption{\(\rho\), \(k = 1\)}
	\end{subfigure}
	\\
	\begin{subfigure}{0.32\linewidth}
		\centering
		\includegraphics[width=\linewidth]{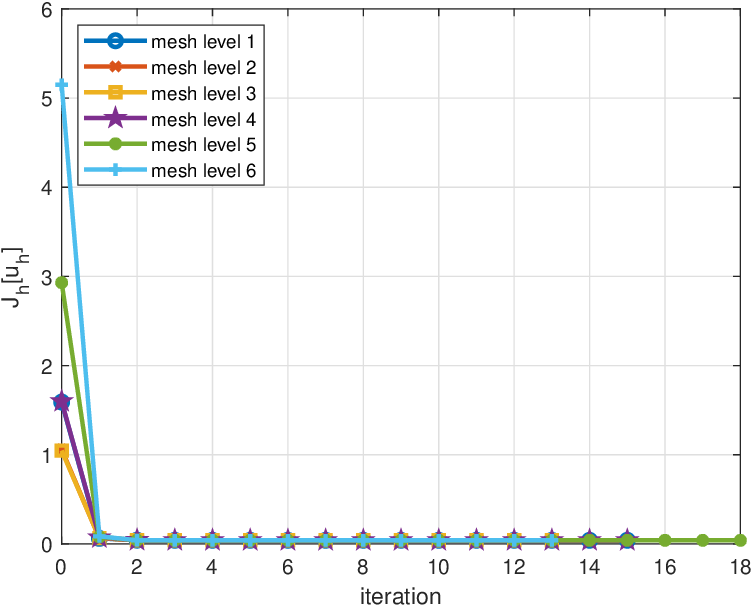}
		\caption{\(J_{h}(u_{h})\), \(k = 6\)}
	\end{subfigure}
	\begin{subfigure}{0.32\linewidth}
		\centering
		\includegraphics[width=\linewidth]{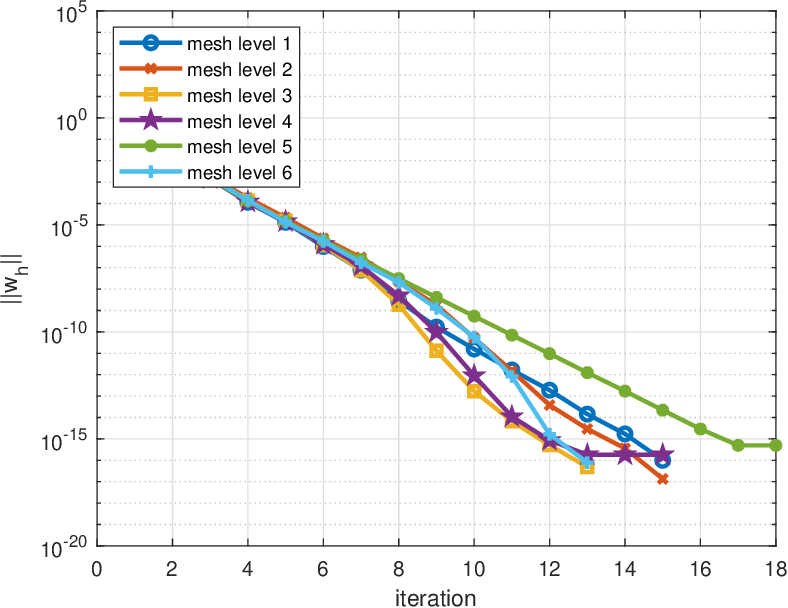}
		\caption{\(\norm{w_{h}}\), \(k = 6\)}
	\end{subfigure}
	\begin{subfigure}{0.32\linewidth}
		\centering
		\includegraphics[width=\linewidth]{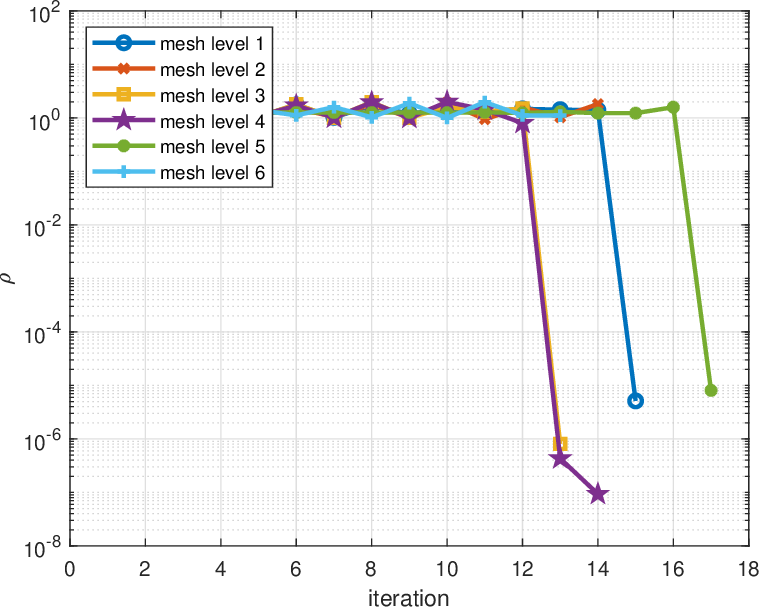}
		\caption{\(\rho\), \(k = 6\)}
	\end{subfigure}
	\caption{Convergence history for Example \ref{exm:re}, with \(\sigma = 0\) and \(p = 1.5\). }
	\label{pic: 2D p-Laplace Dirichlet LDG, example 2, test1}
\end{figure}

\begin{table}[htbp]
	\centering
	\caption{Error table for Example \ref{exm:re}. Parameters: \(\sigma = 7\), \(p = 4\), \(\epsilon = 10^{-14}\), \(\delta_{w_{h}} = 10^{-16}\), \(\delta_{\rho} = 10^{-16}\).}
	\label{table: 2D p-Laplace Dirichlet LDG, example 2, test2}
	\begin{tabular}{| c | c c | c c | c c | c c |}
		\hline

		\(k\) & \(N_{e}\) & \(N_{\rm dof}\) & \(\norm{u - u_{h}}_{L^p}\) & order & \(\norm{\bm q - \bm q_{h}}_{L^p}\) & order & \(\norm{\bm\sigma - \bm\sigma_{h}}_{L^{p'}}\) & order \\
		
		\hline

		1	& 7 & 21 & 1.1225e-01 & - & 4.5324e-01 & - & 1.9580e-01 & - \\
			& 28 & 84 & 4.8591e-02 & 1.2079 & 2.6553e-01 & 0.7714 & 7.0942e-02 & 1.4647 \\
			& 112 & 336 & 1.2312e-02 & 1.9806 & 7.8134e-02 & 1.7648 & 3.3611e-02 & 1.0777 \\
			& 448 & 1344 & 3.3495e-03 & 1.8781 & 3.8148e-02 & 1.0343 & 1.8592e-02 & 0.8543 \\
			& 1792 & 5376 & 9.1362e-04 & 1.8743 & 1.9966e-02 & 0.9341 & 1.0293e-02 & 0.8530 \\
			& 7168 & 21504 & 2.4744e-04 & 1.8845 & 1.0240e-02 & 0.9634 & 5.4531e-03 & 0.9165 \\

		\hline

		2	& 7 & 42 & 9.4888e-02 & - & 3.9841e-01 & - & 1.0562e-01 & - \\
			& 28 & 168 & 2.3261e-02 & 2.0283 & 1.0691e-01 & 1.8979 & 1.8035e-02 & 2.5499 \\
			& 112 & 672 & 5.3419e-03 & 2.1225 & 1.3977e-02 & 2.9353 & 2.9171e-03 & 2.6282 \\
			& 448 & 2688 & 1.1927e-03 & 2.1631 & 1.8476e-03 & 2.9193 & 5.7511e-04 & 2.3426 \\
			& 1792 & 10752 & 2.5934e-04 & 2.2013 & 3.7062e-04 & 2.3176 & 1.3894e-04 & 2.0494 \\
			& 7168 & 43008 & 5.5541e-05 & 2.2232 & 9.0708e-05 & 2.0306 & 3.5456e-05 & 1.9703 \\

		\hline

		3	& 7 & 70 & 3.5097e-02 & - & 2.2779e-01 & - & 4.0262e-02 & - \\
			& 28 & 280 & 1.0875e-02 & 1.6904 & 4.0079e-02 & 2.5068 & 2.9800e-03 & 3.7560 \\
			& 112 & 1120 & 1.9920e-03 & 2.4487 & 4.1837e-03 & 3.2600 & 1.9061e-04 & 3.9667 \\
			& 448 & 4480 & 3.5111e-04 & 2.5042 & 4.6552e-04 & 3.1679 & 1.4127e-05 & 3.7540 \\
			& 1792 & 17920 & 6.0286e-05 & 2.5420 & 5.1760e-05 & 3.1689 & 1.3110e-06 & 3.4298 \\
			& 7168 & 71680 & 1.0218e-05 & 2.5608 & 6.6295e-06 & 2.9649 & 1.4394e-07 & 3.1871 \\

		\hline

		4	& 7 & 105 & 2.4579e-02 & - & 1.5532e-01 & - & 1.1192e-02 & - \\
			& 28 & 420 & 4.4327e-03 & 2.4712 & 1.0572e-02 & 3.8770 & 3.1600e-04 & 5.1464 \\
			& 112 & 1680 & 6.4281e-04 & 2.7857 & 1.1827e-03 & 3.1601 & 9.0793e-06 & 5.1212 \\
			& 448 & 6720 & 8.8926e-05 & 2.8537 & 1.3170e-04 & 3.1667 & 2.8590e-07 & 4.9890 \\
			& 1792 & 26880 & 1.2062e-05 & 2.8821 & 1.4725e-05 & 3.1610 & 1.4783e-08 & 4.2735 \\
			& 7168 & 107520 & 1.6240e-06 & 2.8928 & 4.9436e-05 & -1.7473 & 1.1314e-09 & 3.7077 \\

		\hline

		5	& 7 & 147 & 1.2837e-02 & - & 6.6089e-02 & - & 2.4940e-03 & - \\
			& 28 & 588 & 2.7142e-03 & 2.2417 & 6.4898e-03 & 3.3482 & 4.2796e-05 & 5.8648 \\
			& 112 & 2352 & 3.0470e-04 & 3.1551 & 7.2250e-04 & 3.1671 & 6.5958e-07 & 6.0198 \\
			& 448 & 9408 & 3.3073e-05 & 3.2036 & 8.0512e-05 & 3.1657 & 1.0836e-08 & 5.9276 \\
			& 1792 & 37632 & 3.5545e-06 & 3.2180 & 8.9970e-06 & 3.1617 & 2.0937e-09 & 2.3718 \\
			& 7168 & 150528 & 3.7845e-07 & 3.2315 & 1.0077e-06 & 3.1584 & 4.4897e-09 & -1.1006 \\
		
		\hline

		6   & 7 & 196 & 5.6444e-03 & - & 3.3924e-02 & - & 4.2407e-04 & - \\
			& 28 & 784 & 7.2374e-04 & 2.9633 & 2.8283e-03 & 3.5843 & 2.1715e-06 & 7.6095 \\
			& 112 & 3136 & 6.2957e-05 & 3.5230 & 3.1497e-04 & 3.1666 & 1.7375e-08 & 6.9656 \\
			& 448 & 12544 & 5.3619e-06 & 3.5535 & 3.4873e-05 & 3.1751 & 5.0282e-09 & 1.7889 \\
			& 1792 & 50176 & 4.5741e-07 & 3.5512 & 3.8761e-06 & 3.1694 & 1.0701e-08 & -1.0897 \\
			& 7168 & 200704 & 8.3135e-08 & 2.4600 & 5.5058e-06 & -0.5064 & 2.3028e-09 & 2.2163 \\

		\hline
	\end{tabular}
\end{table}

\begin{figure}[htbp]
	\centering
	\begin{subfigure}{0.32\linewidth}
		\centering
		\includegraphics[width=\linewidth]{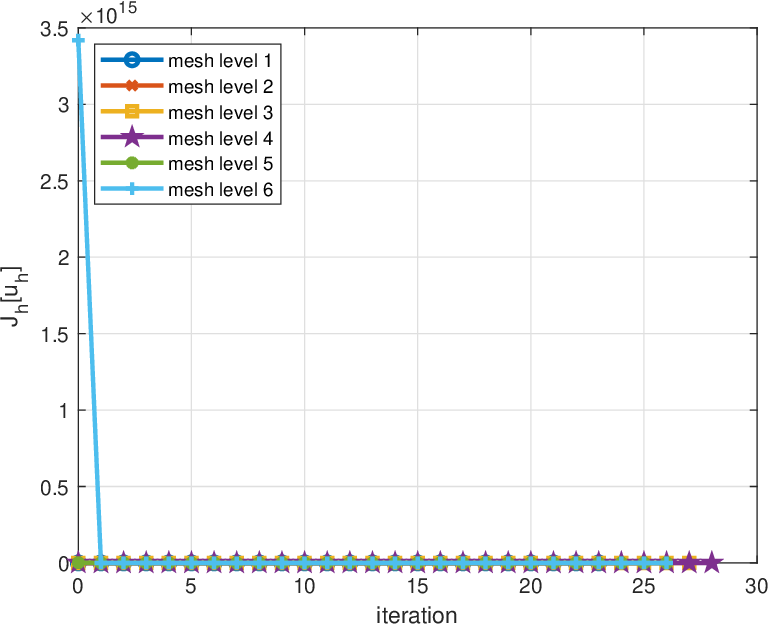}
		\caption{\(J_{h}(u_{h})\), \(k = 1\)}
	\end{subfigure}
	\begin{subfigure}{0.32\linewidth}
		\centering
		\includegraphics[width=\linewidth]{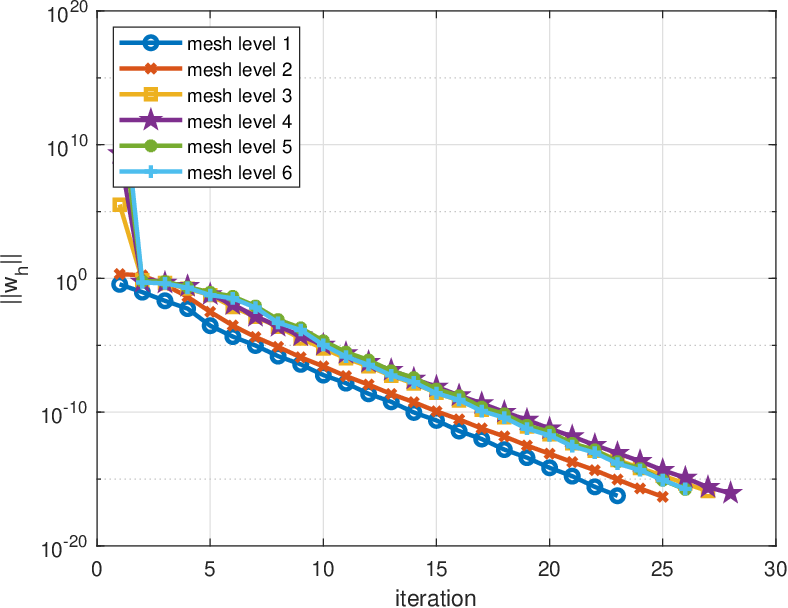}
		\caption{\(\norm{w_{h}}\), \(k = 1\)}
	\end{subfigure}
	\begin{subfigure}{0.32\linewidth}
		\centering
		\includegraphics[width=\linewidth]{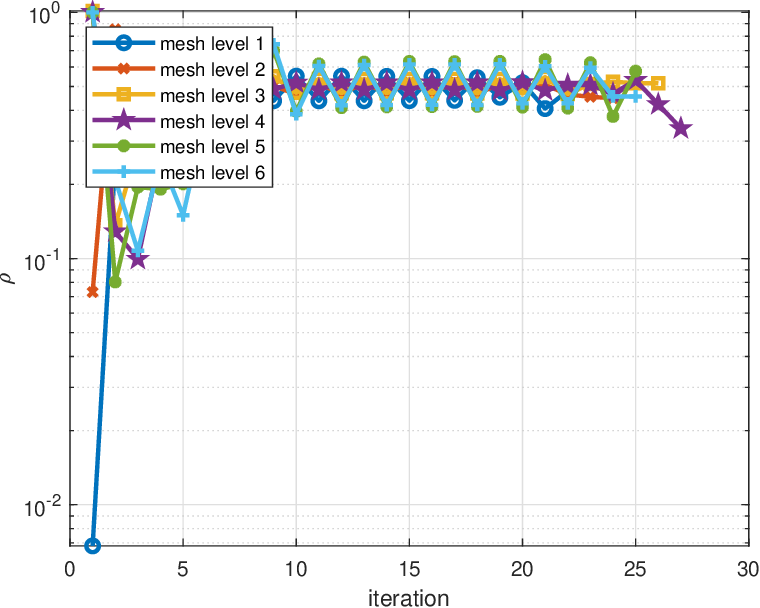}
		\caption{\(\rho\), \(k = 1\)}
	\end{subfigure}
	\\
	\begin{subfigure}{0.32\linewidth}
		\centering
		\includegraphics[width=\linewidth]{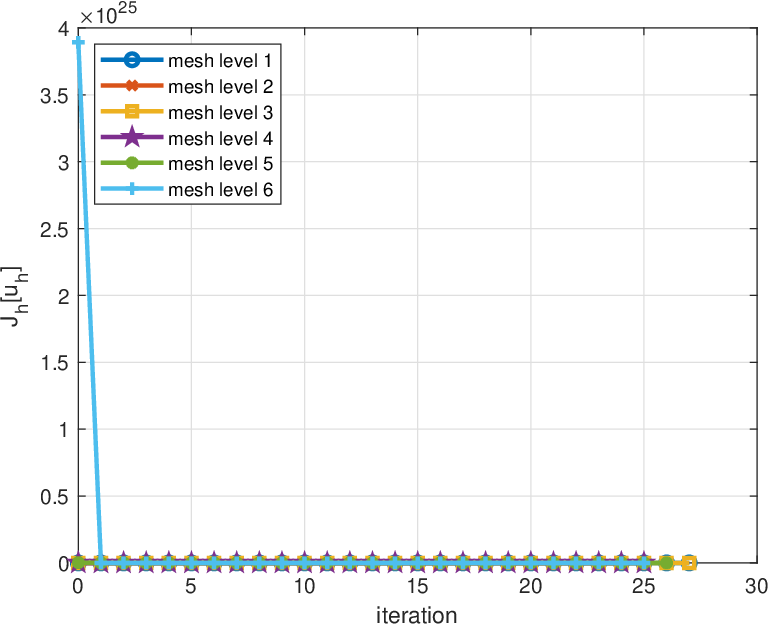}
		\caption{\(J_{h}(u_{h})\), \(k = 6\)}
	\end{subfigure}
	\begin{subfigure}{0.32\linewidth}
		\centering
		\includegraphics[width=\linewidth]{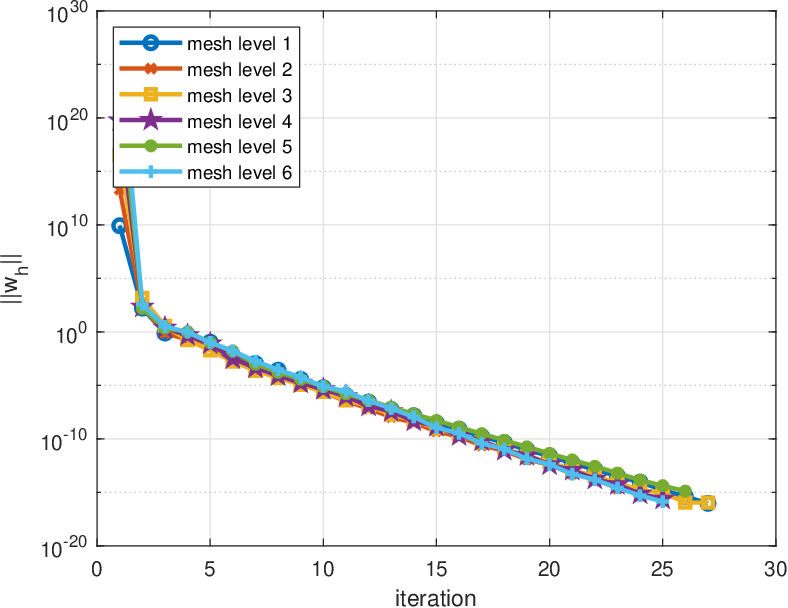}
		\caption{\(\norm{w_{h}}\), \(k = 6\)}
	\end{subfigure}
	\begin{subfigure}{0.32\linewidth}
		\centering
		\includegraphics[width=\linewidth]{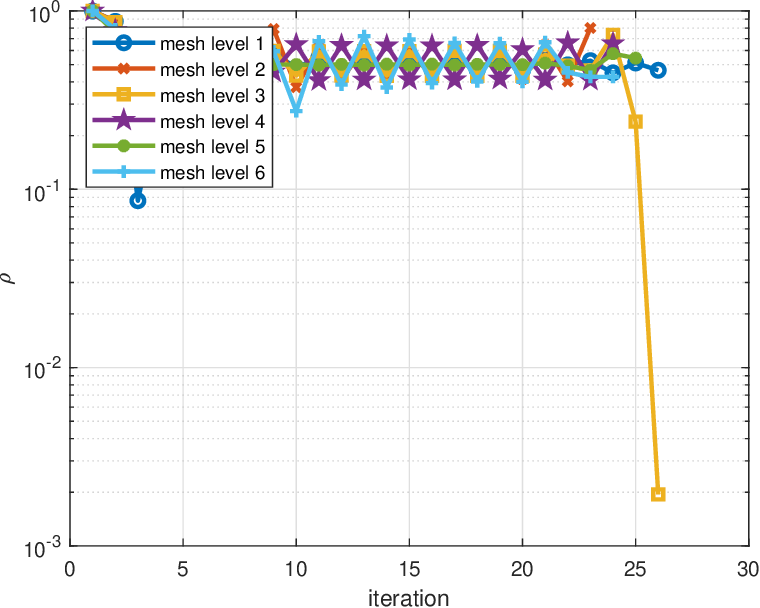}
		\caption{\(\rho\), \(k = 6\)}
	\end{subfigure}
	\caption{Convergence history for Example \ref{exm:re}, with \(\sigma = 7\) and \(p = 4\). }
	\label{pic: 2D p-Laplace Dirichlet LDG, example 2, test2}
\end{figure}

\begin{example}{\bf A Degenerate Case}\label{exm:de}
\end{example}

We take the radial exact solution \cite{10.1090/S0025-5718-1993-1192966-4} to be 
\begin{eqnarray*}
	&& u(x,y) = \begin{cases}
		0, & r(x,y) < a, \\
		(r(x,y) - a)^{4}, & r(x,y) \geqslant a, 
	\end{cases} \in W^{4, \infty} \cap W^{4 + \frac{1}{p} - \varepsilon, p}, \\
	&& \bm q(x,y) = \begin{cases}
		\bm 0, & r(x,y) < a, \\
		\begin{pmatrix}
			4 (r(x,y)-a)^{3} \frac{x}{r(x,y)} \\
			4 (r(x,y)-a)^{3} \frac{y}{r(x,y)} 
		\end{pmatrix}, & r(x,y) \geqslant a, 
	\end{cases} \in W^{3, \infty} \cap W^{3 + \frac{1}{p} - \varepsilon, p}, \\
	&& \bm\sigma(x,y) = \begin{cases}
		\bm 0, & r(x,y) < a, \\
		\begin{pmatrix}
			4^{p-1} (r(x,y)-a)^{3p-3} \frac{x}{r(x,y)} \\
			4^{p-1} (r(x,y)-a)^{3p-3} \frac{y}{r(x,y)} 
		\end{pmatrix}, & r(x,y) \geqslant a, 
	\end{cases} \in W^{3p-3, \infty} \cap W^{3p-2-\frac{1}{p} - \varepsilon, p'}, \\
	&& f(x,y) = \begin{cases}
		0, & r(x,y) < a, \\
		4^{p-1} (r(x,y) - a)^{3p-4} \left(2 - 3p + \frac{a}{r(x,y)}\right), & r(x,y) \geqslant a. 
	\end{cases}
\end{eqnarray*}
The gradient vanishes on \(B_{a}(\bm 0)\), so the problem is degenerate when \(p > 2\). Just as in the article \cite{10.1090/S0025-5718-1993-1192966-4}, we choose \(a = 0.3\) and \(p = 4\). The results are shown in Table \ref{table: 2D p-Laplace Dirichlet LDG, example 3, test1} and Figure \ref{pic: 2D p-Laplace Dirichlet LDG, example 3, test1}. The number of iterations is around 25 for most \(h\)'s and \(k\)'s. 

In this example, \((u, \bm q, \bm\sigma) \in W^{s + 1 - \varepsilon, p} \times W^{s - \varepsilon, p} \times W^{r - \varepsilon, p'}\) with \((s, r) = (3.25, 9.75)\). The convergence rates are around the best values for \(\bm q\) and \(\bm\sigma\) in comparison to the linear case. 

\begin{table}[htbp]
	\centering
	\caption{Error table for Example \ref{exm:de}. Parameters: \(p = 4\), \(\epsilon = 10^{-14}\), \(\delta_{w_{h}} = 10^{-16}\), \(\delta_{\rho} = 10^{-16}\).}
	\label{table: 2D p-Laplace Dirichlet LDG, example 3, test1}
	\begin{tabular}{| c | c c | c c | c c | c c |}
		\hline

		\(k\) & \(N_{e}\) & \(N_{\rm dof}\) & \(\norm{u - u_{h}}_{L^p}\) & order & \(\norm{\bm q - \bm q_{h}}_{L^p}\) & order & \(\norm{\bm\sigma - \bm\sigma_{h}}_{L^{p'}}\) & order \\

		\hline

		1	& 7 & 21 & 5.0087e-01 & - & 2.1984e+00 & - & 1.2909e+01 & - \\
			& 28 & 84 & 1.9990e-01 & 1.3252 & 1.1516e+00 & 0.9328 & 6.0890e+00 & 1.0841 \\
			& 112 & 336 & 6.4463e-02 & 1.6327 & 5.1304e-01 & 1.1665 & 2.3958e+00 & 1.3457 \\
			& 448 & 1344 & 1.7066e-02 & 1.9173 & 1.9450e-01 & 1.3993 & 1.2788e+00 & 0.9057 \\
			& 1792 & 5376 & 4.6698e-03 & 1.8697 & 1.0157e-01 & 0.9373 & 7.1692e-01 & 0.8349 \\
			& 7168 & 21504 & 1.2687e-03 & 1.8799 & 5.2588e-02 & 0.9497 & 3.8597e-01 & 0.8933 \\

		\hline

		2	& 7 & 42 & 3.4929e-01 & - & 1.5828e+00 & - & 1.0810e+01 & - \\
			& 28 & 168 & 1.4205e-01 & 1.2981 & 8.1954e-01 & 0.9496 & 2.0331e+00 & 2.4106 \\
			& 112 & 672 & 3.0138e-02 & 2.2367 & 1.4992e-01 & 2.4506 & 3.5512e-01 & 2.5173 \\
			& 448 & 2688 & 6.8604e-03 & 2.1352 & 1.6884e-02 & 3.1505 & 6.4558e-02 & 2.4596 \\
			& 1792 & 10752 & 1.5120e-03 & 2.1819 & 2.8383e-03 & 2.5726 & 1.4676e-02 & 2.1372 \\
			& 7168 & 43008 & 3.2590e-04 & 2.2140 & 6.7299e-04 & 2.0764 & 3.7017e-03 & 1.9872 \\

		\hline

		3	& 7 & 70 & 3.0707e-01 & - & 1.5849e+00 & - & 5.5966e+00 & - \\
			& 28 & 280 & 7.2396e-02 & 2.0846 & 4.8447e-01 & 1.7099 & 5.2846e-01 & 3.4047 \\
			& 112 & 1120 & 1.2841e-02 & 2.4951 & 4.8259e-02 & 3.3275 & 3.7767e-02 & 3.8066 \\
			& 448 & 4480 & 2.3307e-03 & 2.4620 & 4.0781e-03 & 3.5648 & 2.5268e-03 & 3.9017 \\
			& 1792 & 17920 & 4.0573e-04 & 2.5222 & 3.7353e-04 & 3.4486 & 2.0165e-04 & 3.6474 \\
			& 7168 & 71680 & 6.9167e-05 & 2.5524 & 1.6427e-04 & 1.1851 & 2.0028e-05 & 3.3318 \\

		\hline

		4	& 7 & 105 & 1.1334e-01 & - & 9.1652e-01 & - & 2.2549e+00 & - \\
			& 28 & 420 & 3.2154e-02 & 1.8176 & 1.6407e-01 & 2.4819 & 1.0525e-01 & 4.4211 \\
			& 112 & 1680 & 4.8426e-03 & 2.7312 & 1.4662e-02 & 3.4842 & 2.9308e-03 & 5.1664 \\
			& 448 & 6720 & 6.9000e-04 & 2.8111 & 1.6002e-03 & 3.1957 & 8.9106e-05 & 5.0396 \\
			& 1792 & 26880 & 9.4671e-05 & 2.8656 & 2.6576e-04 & 2.5901 & 3.0897e-06 & 4.8500 \\
			& 7168 & 107520 & 1.2770e-05 & 2.8902 & 1.0274e-04 & 1.3711 & 7.3647e-07 & 2.0688 \\

		\hline
		
		5	& 7 & 147 & 8.6038e-02 & - & 5.7931e-01 & - & 7.9631e-01 & - \\
			& 28 & 588 & 2.2907e-02 & 1.9092 & 7.2822e-02 & 2.9919 & 1.9460e-02 & 5.3547 \\
			& 112 & 2352 & 2.6796e-03 & 3.0957 & 5.0373e-03 & 3.8536 & 3.1838e-04 & 5.9336 \\
			& 448 & 9408 & 2.9864e-04 & 3.1656 & 6.1208e-04 & 3.0409 & 5.0076e-06 & 5.9905 \\
			& 1792 & 37632 & 3.2344e-05 & 3.2068 & 6.5927e-05 & 3.2148 & 2.2571e-07 & 4.4716 \\
			& 7168 & 150528 & 3.4570e-06 & 3.2259 & 8.1953e-05 & -0.3139 & 2.7001e-07 & -0.2585 \\
		
		\hline

		6   & 7 & 196 & 4.9951e-02 & - & 4.3889e-01 & - & 1.9182e-01 & - \\
			& 28 & 784 & 7.6588e-03 & 2.7053 & 2.5296e-02 & 4.1169 & 1.8994e-03 & 6.6581 \\
			& 112 & 3136 & 6.9588e-04 & 3.4602 & 2.9030e-03 & 3.1233 & 1.1439e-05 & 7.3754 \\
			& 448 & 12544 & 6.0415e-05 & 3.5259 & 2.9858e-04 & 3.2814 & 4.7298e-07 & 4.5961 \\
			& 1792 & 50176 & 5.2296e-06 & 3.5301 & 1.3389e-04 & 1.1571 & 3.3386e-07 & 0.5025 \\
			& 7168 & 200704 & 5.9421e-07 & 3.1377 & 4.6034e-05 & 1.5403 & 1.8878e-07 & 0.8225 \\

		\hline
	\end{tabular}
\end{table}

\begin{figure}[htbp]
	\centering
	\begin{subfigure}{0.32\linewidth}
		\centering
		\includegraphics[width=\linewidth]{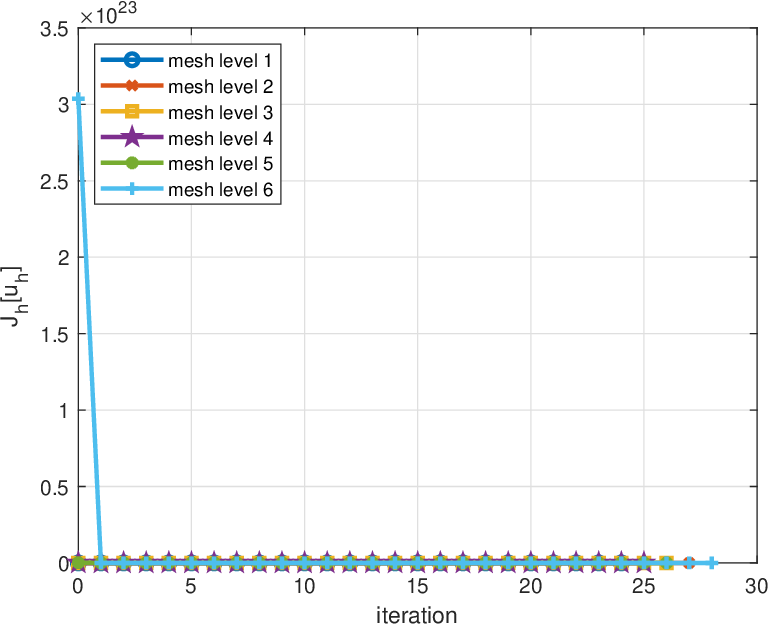}
		\caption{\(J_{h}(u_{h})\), \(k = 1\)}
	\end{subfigure}
	\begin{subfigure}{0.32\linewidth}
		\centering
		\includegraphics[width=\linewidth]{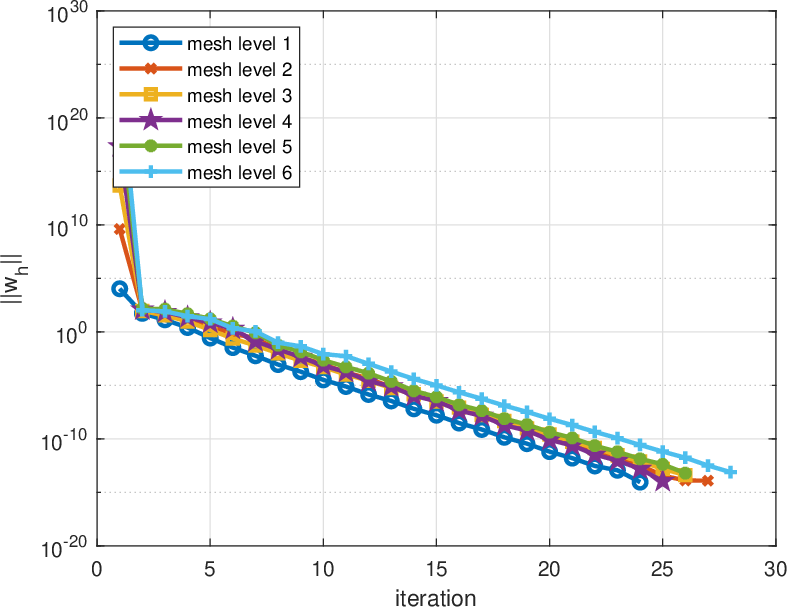}
		\caption{\(\norm{w_{h}}\), \(k = 1\)}
	\end{subfigure}
	\begin{subfigure}{0.32\linewidth}
		\centering
		\includegraphics[width=\linewidth]{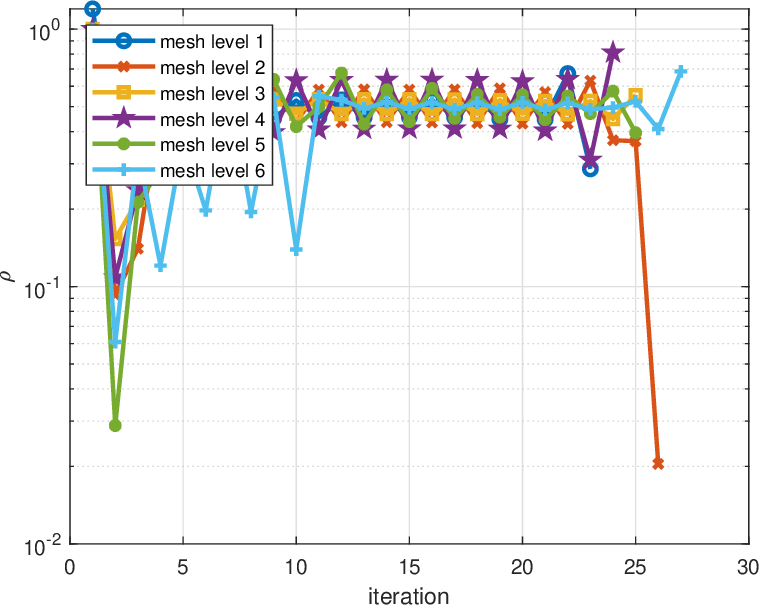}
		\caption{\(\rho\), \(k = 1\)}
	\end{subfigure}
	\\
	\begin{subfigure}{0.32\linewidth}
		\centering
		\includegraphics[width=\linewidth]{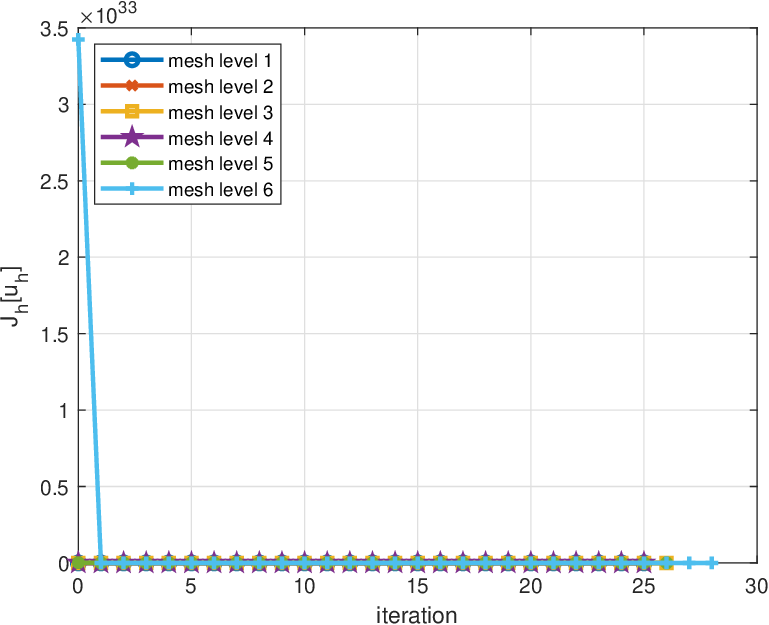}
		\caption{\(J_{h}(u_{h})\), \(k = 6\)}
	\end{subfigure}
	\begin{subfigure}{0.32\linewidth}
		\centering
		\includegraphics[width=\linewidth]{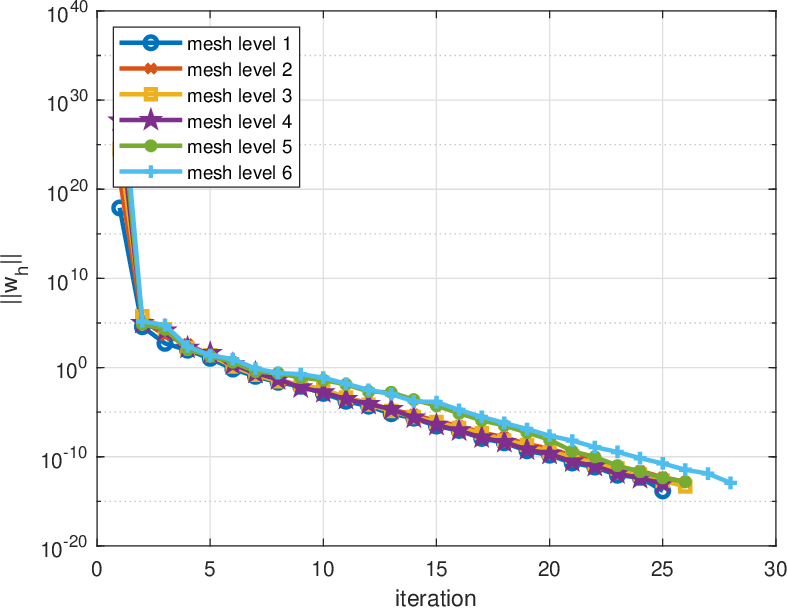}
		\caption{\(\norm{w_{h}}\), \(k = 6\)}
	\end{subfigure}
	\begin{subfigure}{0.32\linewidth}
		\centering
		\includegraphics[width=\linewidth]{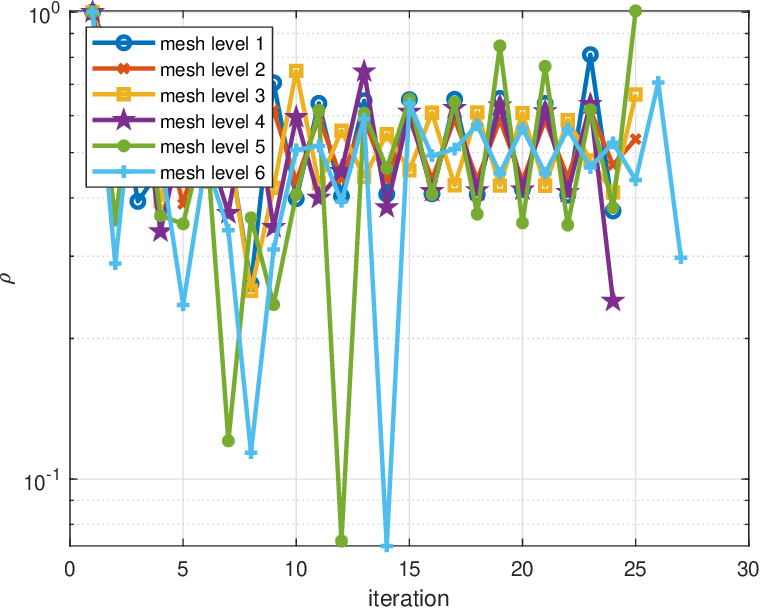}
		\caption{\(\rho\), \(k = 6\)}
	\end{subfigure}
	\caption{Convergence history for Example \ref{exm:de}, with \(p = 4\). }
	\label{pic: 2D p-Laplace Dirichlet LDG, example 3, test1}
\end{figure}

\begin{example}{\bf A Smooth Case}\label{exm:sm}
\end{example}

Take the square domain \(\Omega = [1, 2]^2\), and take the radial \(p\)-harmonic exact solution \cite{10.1007/s10915-007-9134-z, 10.1137/15M1008014} to be 
\begin{eqnarray*}
	&& u(x,y) = r(x,y)^{\frac{p-2}{p-1}}, \\
	&& \bm q(x,y) = \begin{pmatrix}
		\frac{p-2}{p-1} r(x,y)^{\frac{-p}{p-1}} x \\
		\frac{p-2}{p-1} r(x,y)^{\frac{-p}{p-1}} y 
	\end{pmatrix}, \\
	&& \bm\sigma(x,y) = \begin{pmatrix}
		\mathrm{sgn}(p-2) \abs{\frac{p-2}{p-1}}^{p-1} \frac{x}{r(x,y)^2} \\
		\mathrm{sgn}(p-2) \abs{\frac{p-2}{p-1}}^{p-1} \frac{y}{r(x,y)^2} 
	\end{pmatrix}, \\
	&& f(x,y) = 0. 
\end{eqnarray*}
This solution is smooth, and moreover its gradient is bounded from both above and below in the domain. We use the same parameters as in the article \cite{10.1137/15M1008014}, \(p = 1.5\) and \(p = 3\). For both tests, a series of uniformly refined conforming triangular meshes is used. The numerical results for \(p = 1.5\) and \(p = 3\) are reported in Table \ref{table: 2D p-Laplace Dirichlet LDG, example 4, test1} with Figure \ref{pic: 2D p-Laplace Dirichlet LDG, example 4, test1} and Table \ref{table: 2D p-Laplace Dirichlet LDG, example 4, test2} with Figure \ref{pic: 2D p-Laplace Dirichlet LDG, example 4, test2}. The number of iterations of these two cases is all around 15, which indicates independence of \(h\) and \(k\). 

The failure to reduce the error when \(k\) is big or \(h\) is small is possibly due to early termination of the gradient descent iteration and floating-point underflow. For small \(k\)'s, the results show the best convergence rates for all variables both when \(p = 1.5\) and when \(p = 3\), compared against the linear case. Similar results are also observed with the HDG method \cite{10.1137/15M1008014}. 

\begin{table}[htbp]
	\centering
	\caption{Error table for Example \ref{exm:sm}. Parameters: \(p = 1.5\), \(\epsilon = 10^{-14}\), \(\delta_{w_{h}} = 10^{-16}\), \(\delta_{\rho} = 10^{-16}\).}
	\label{table: 2D p-Laplace Dirichlet LDG, example 4, test1}
	\begin{tabular}{| c | c c | c c | c c | c c |}
		\hline

		\(k\) & \(N_{e}\) & \(N_{\rm dof}\) & \(\norm{u - u_{h}}_{L^p}\) & order & \(\norm{\bm q - \bm q_{h}}_{L^p}\) & order & \(\norm{\bm\sigma - \bm\sigma_{h}}_{L^{p'}}\) & order \\

		\hline

		1	& 16 & 48 & 1.0372e-03 & - & 2.8342e-02 & - & 4.1536e-02 & - \\
			& 64 & 192 & 2.7955e-04 & 1.8915 & 1.4451e-02 & 0.9718 & 2.2825e-02 & 0.8638 \\
			& 256 & 768 & 7.0869e-05 & 1.9799 & 7.2741e-03 & 0.9903 & 1.1631e-02 & 0.9726 \\
			& 1024 & 3072 & 1.7819e-05 & 1.9917 & 3.6433e-03 & 0.9975 & 5.8448e-03 & 0.9928 \\
			& 4096 & 12288 & 4.4670e-06 & 1.9960 & 1.8226e-03 & 0.9992 & 2.9266e-03 & 0.9979 \\
			& 16384 & 49152 & 1.1181e-06 & 1.9983 & 9.1150e-04 & 0.9997 & 1.4640e-03 & 0.9993 \\

		\hline

		2	& 16 & 96 & 8.6236e-05 & - & 2.2770e-03 & - & 3.8375e-03 & - \\
			& 64 & 384 & 1.0855e-05 & 2.9899 & 5.6497e-04 & 2.0109 & 1.0040e-03 & 1.9344 \\
			& 256 & 1536 & 1.3526e-06 & 3.0046 & 1.4138e-04 & 1.9986 & 2.5810e-04 & 1.9598 \\
			& 1024 & 6144 & 1.6836e-07 & 3.0061 & 3.5423e-05 & 1.9968 & 6.5469e-05 & 1.9790 \\
			& 4096 & 24576 & 2.0987e-08 & 3.0040 & 8.8690e-06 & 1.9979 & 1.6488e-05 & 1.9894 \\
			& 16384 & 98304 & 2.6235e-09 & 2.9999 & 2.2191e-06 & 1.9988 & 4.1371e-06 & 1.9947 \\

		\hline

		3	& 16 & 160 & 4.6017e-06 & - & 1.5681e-04 & - & 3.3194e-04 & - \\
			& 64 & 640 & 2.6963e-07 & 4.0931 & 1.9628e-05 & 2.9980 & 4.1896e-05 & 2.9860 \\
			& 256 & 2560 & 1.6227e-08 & 4.0545 & 2.4526e-06 & 3.0005 & 5.2193e-06 & 3.0049 \\
			& 1024 & 10240 & 1.0227e-09 & 3.9879 & 3.0633e-07 & 3.0012 & 6.5041e-07 & 3.0045 \\
			& 4096 & 40960 & 5.3706e-10 & 0.9292 & 3.9702e-08 & 2.9478 & 8.1706e-08 & 2.9928 \\
			& 16384 & 163840 & 2.2090e-10 & 1.2817 & 6.2691e-09 & 2.6629 & 1.1142e-08 & 2.8745 \\

		\hline

		4	& 16 & 240 & 2.5374e-07 & - & 1.0718e-05 & - & 2.4662e-05 & - \\
			& 64 & 960 & 7.9948e-09 & 4.9882 & 6.7138e-07 & 3.9967 & 1.5550e-06 & 3.9872 \\
			& 256 & 3840 & 2.6664e-10 & 4.9061 & 4.1716e-08 & 4.0085 & 9.5916e-08 & 4.0190 \\
			& 1024 & 15360 & 5.1579e-11 & 2.3700 & 2.9589e-09 & 3.8175 & 6.0748e-09 & 3.9809 \\
			& 4096 & 61440 & 3.4557e-10 & -2.7441 & 5.1145e-09 & -0.7895 & 7.6760e-09 & -0.3375 \\
			& 16384 & 245760 & 1.1062e-11 & 4.9653 & 2.1556e-10 & 4.5685 & 5.9605e-10 & 3.6868 \\

		\hline
		
		5	& 16 & 336 & 1.3848e-08 & - & 7.8500e-07 & - & 2.0262e-06 & - \\
			& 64 & 1344 & 2.3354e-10 & 5.8899 & 2.5019e-08 & 4.9716 & 6.4595e-08 & 4.9712 \\
			& 256 & 5376 & 3.9682e-10 & -0.7648 & 5.7921e-09 & 2.1109 & 8.2972e-09 & 2.9607 \\
			& 1024 & 21504 & 3.5467e-10 & 0.1620 & 8.0309e-09 & -0.4715 & 1.1609e-08 & -0.4846 \\	
			& 4096 & 86016 & 4.3664e-10 & -0.3000 & 6.1614e-09 & 0.3823 & 9.2395e-09 & 0.3294 \\
			& 16384 & 344064 & 9.1259e-11 & 2.2584 & 2.0321e-09 & 1.6003 & 3.6673e-09 & 1.3331 \\
		
		\hline
		
		6	& 16 & 448 & 8.4207e-10 & - & 5.6059e-08 & - & 1.3230e-07 & - \\
			& 64 & 1792 & 1.8888e-10 & 2.1564 & 3.0439e-09 & 4.2029 & 5.0056e-09 & 4.7241 \\
			& 256 & 7168 & 2.2991e-10 & -0.2836 & 5.5574e-09 & -0.8685 & 8.3898e-09 & -0.7451 \\
			& 1024 & 28672 & 2.8080e-10 & -0.2885 & 6.6446e-09 & -0.2578 & 9.9341e-09 & -0.2438 \\
			& 4096 & 114688 & 1.5732e-10 & 0.8359 & 3.4252e-09 & 0.9560 & 5.8425e-09 & 0.7658 \\
			& 16384 & 458752 & 2.9352e-10 & -0.8998 & 4.5425e-09 & -0.4073 & 6.7507e-09 & -0.2085 \\

		\hline
	\end{tabular}
\end{table}

\begin{figure}[htbp]
	\centering
	\begin{subfigure}{0.32\linewidth}
		\centering
		\includegraphics[width=\linewidth]{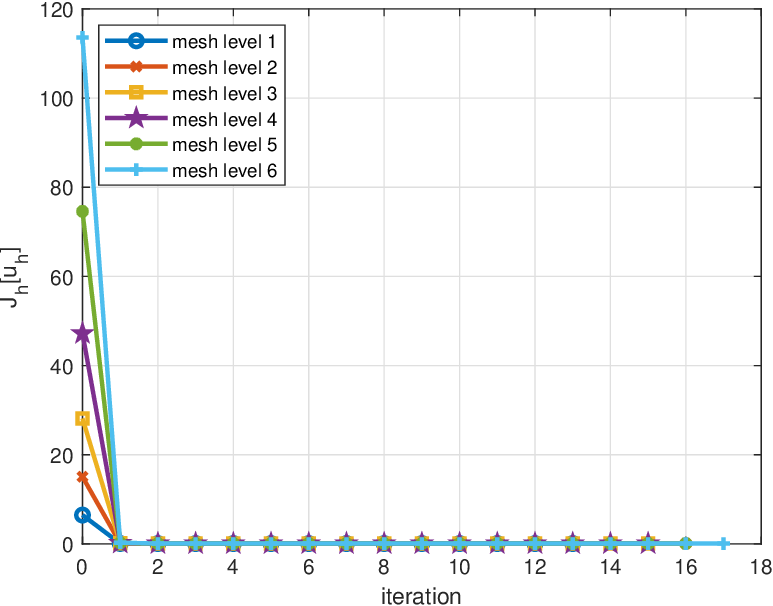}
		\caption{\(J_{h}(u_{h})\), \(k = 1\)}
	\end{subfigure}
	\begin{subfigure}{0.32\linewidth}
		\centering
		\includegraphics[width=\linewidth]{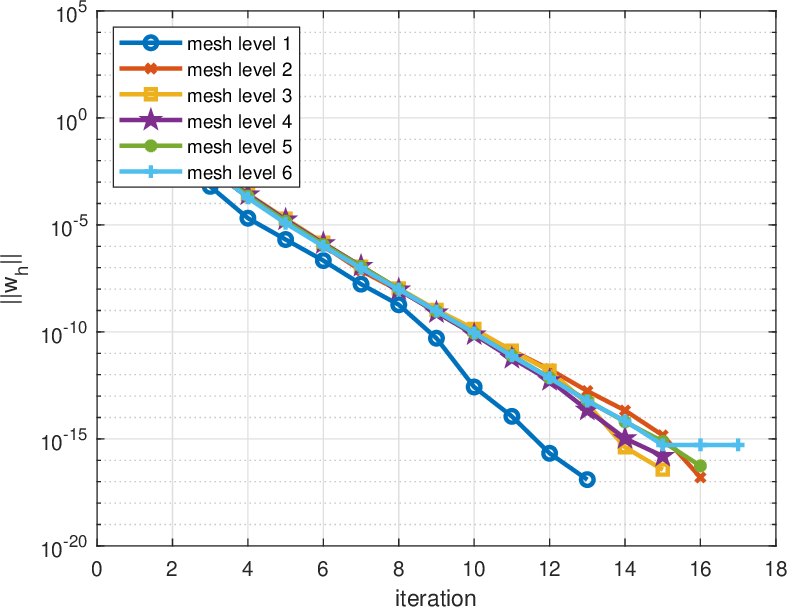}
		\caption{\(\norm{w_{h}}\), \(k = 1\)}
	\end{subfigure}
	\begin{subfigure}{0.32\linewidth}
		\centering
		\includegraphics[width=\linewidth]{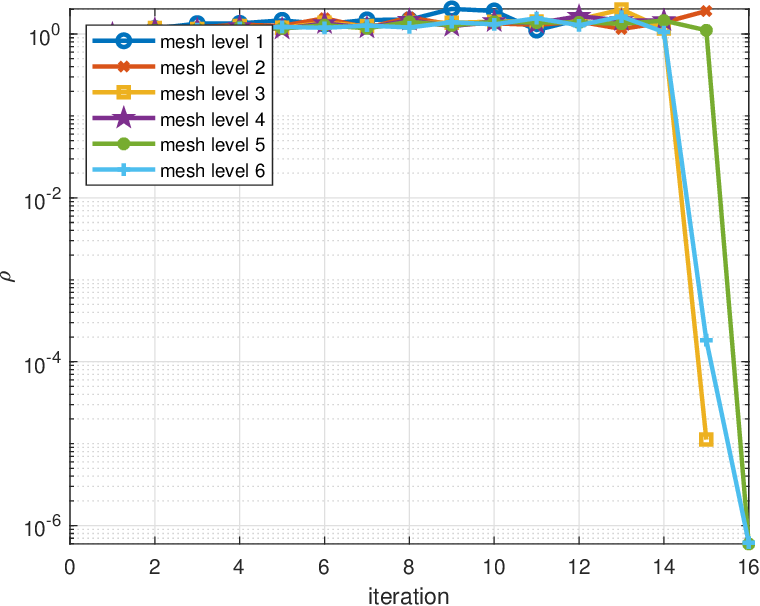}
		\caption{\(\rho\), \(k = 1\)}
	\end{subfigure}
	\\
	\begin{subfigure}{0.32\linewidth}
		\centering
		\includegraphics[width=\linewidth]{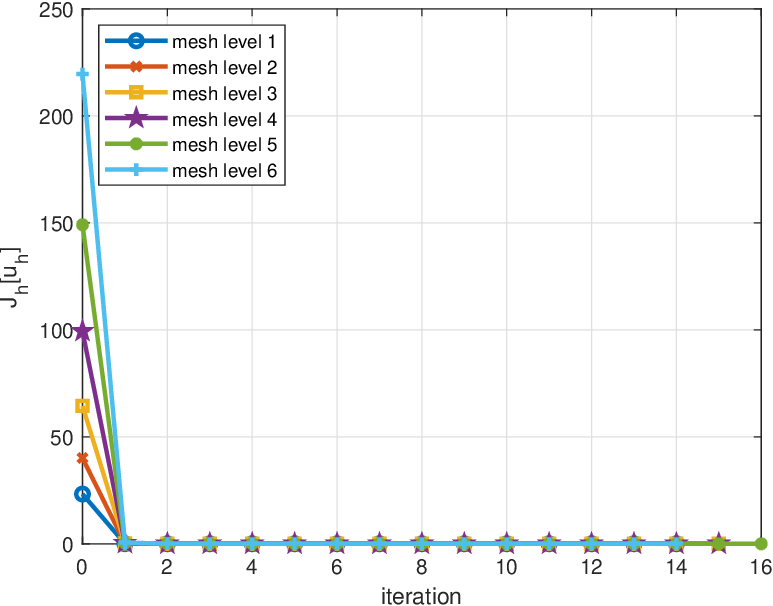}
		\caption{\(J_{h}(u_{h})\), \(k = 6\)}
	\end{subfigure}
	\begin{subfigure}{0.32\linewidth}
		\centering
		\includegraphics[width=\linewidth]{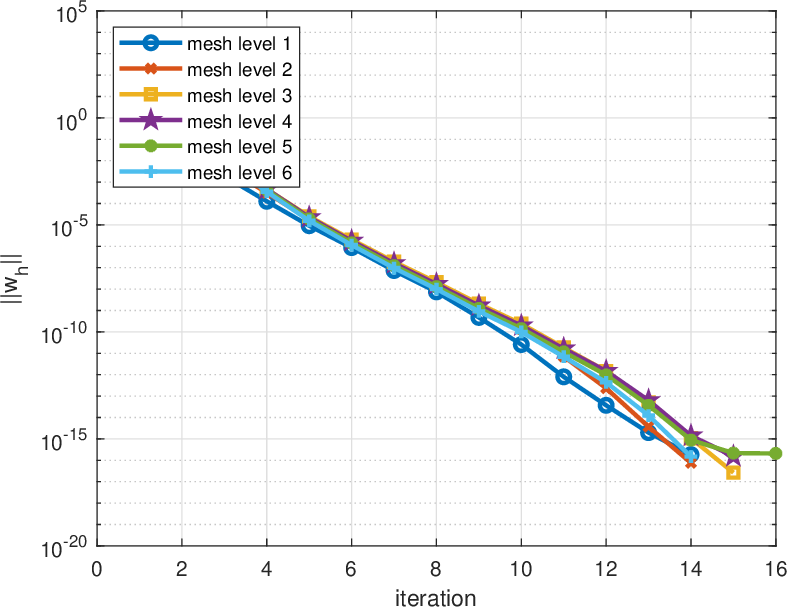}
		\caption{\(\norm{w_{h}}\), \(k = 6\)}
	\end{subfigure}
	\begin{subfigure}{0.32\linewidth}
		\centering
		\includegraphics[width=\linewidth]{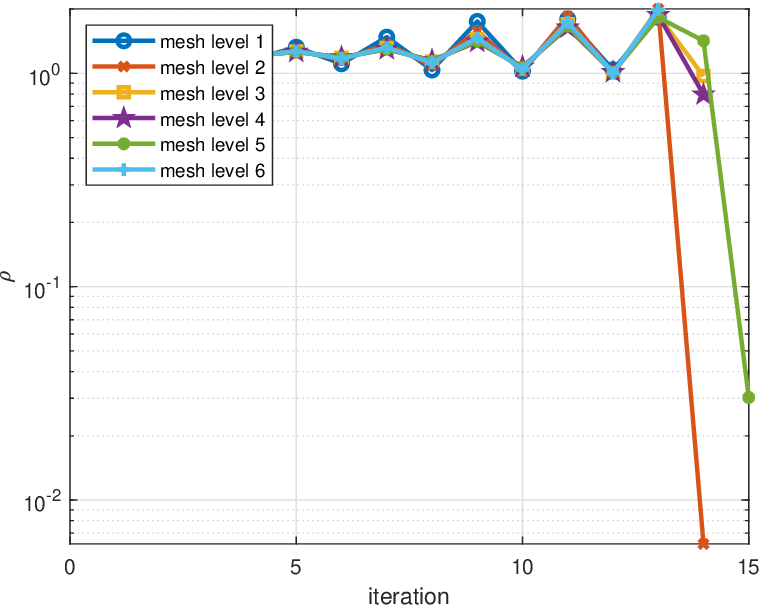}
		\caption{\(\rho\), \(k = 6\)}
	\end{subfigure}
	\caption{Convergence history for Example \ref{exm:sm}, with \(p = 1.5\). }
	\label{pic: 2D p-Laplace Dirichlet LDG, example 4, test1}
\end{figure}

\begin{table}[htbp]
	\centering
	\caption{Error table for Example \ref{exm:sm}. Parameters: \(p = 3\), \(\epsilon = 10^{-14}\), \(\delta_{w_{h}} = 10^{-16}\), \(\delta_{\rho} = 10^{-16}\).}
	\label{table: 2D p-Laplace Dirichlet LDG, example 4, test2}
	\begin{tabular}{| c | c c | c c | c c | c c |}
		\hline

		\(k\) & \(N_{e}\) & \(N_{\rm dof}\) & \(\norm{u - u_{h}}_{L^p}\) & order & \(\norm{\bm q - \bm q_{h}}_{L^p}\) & order & \(\norm{\bm\sigma - \bm\sigma_{h}}_{L^{p'}}\) & order \\

		\hline

		1	& 16 & 48 & 6.1661e-04 & - & 8.3824e-03 & - & 2.8982e-03 & - \\
			& 64 & 192 & 1.6146e-04 & 1.9332 & 4.3914e-03 & 0.9327 & 1.5995e-03 & 0.8576 \\
			& 256 & 768 & 4.1832e-05 & 1.9485 & 2.2942e-03 & 0.9367 & 8.3392e-04 & 0.9396 \\
			& 1024 & 3072 & 1.0638e-05 & 1.9753 & 1.1691e-03 & 0.9727 & 4.2445e-04 & 0.9743 \\
			& 4096 & 12288 & 2.6789e-06 & 1.9896 & 5.8994e-04 & 0.9867 & 2.1433e-04 & 0.9857 \\
			& 16384 & 49152 & 6.7208e-07 & 1.9949 & 2.9635e-04 & 0.9933 & 1.0787e-04 & 0.9905 \\

		\hline

		2	& 16 & 96 & 2.2017e-05 & - & 5.2367e-04 & - & 1.8345e-04 & - \\
			& 64 & 384 & 2.9965e-06 & 2.8773 & 1.4110e-04 & 1.8919 & 4.9821e-05 & 1.8805 \\
			& 256 & 1536 & 3.8253e-07 & 2.9696 & 3.6936e-05 & 1.9336 & 1.3101e-05 & 1.9271 \\
			& 1024 & 6144 & 4.8458e-08 & 2.9808 & 9.4465e-06 & 1.9672 & 3.3450e-06 & 1.9696 \\
			& 4096 & 24576 & 6.1051e-09 & 2.9886 & 2.3870e-06 & 1.9846 & 8.4330e-07 & 1.9879 \\
			& 16384 & 98304 & 1.0885e-09 & 2.4877 & 5.9969e-07 & 1.9929 & 2.1182e-07 & 1.9932 \\

		\hline

		3	& 16 & 160 & 1.1969e-06 & - & 3.4546e-05 & - & 1.1652e-05 & - \\
			& 64 & 640 & 7.9547e-08 & 3.9114 & 4.5249e-06 & 2.9326 & 1.5044e-06 & 2.9533 \\
			& 256 & 2560 & 5.0486e-09 & 3.9779 & 5.7407e-07 & 2.9786 & 1.9013e-07 & 2.9841 \\
			& 1024 & 10240 & 3.9896e-10 & 3.6616 & 7.3632e-08 & 2.9628 & 2.4486e-08 & 2.9569 \\
			& 4096 & 40960 & 2.5694e-10 & 0.6348 & 2.2007e-08 & 1.7424 & 4.9171e-09 & 2.3161 \\
			& 16384 & 163840 & 1.4915e-10 & 0.7847 & 5.5353e-09 & 1.9912 & 2.0881e-09 & 1.2356 \\

		\hline

		4	& 16 & 240 & 7.8745e-08 & - & 2.3821e-06 & - & 6.6297e-07 & - \\
			& 64 & 960 & 2.6198e-09 & 4.9097 & 1.5224e-07 & 3.9679 & 4.2029e-08 & 3.9795 \\
			& 256 & 3840 & 2.1988e-10 & 3.5747 & 2.6166e-08 & 2.5405 & 5.3147e-09 & 2.9833 \\
			& 1024 & 15360 & 3.0472e-10 & -0.4708 & 3.1129e-08 & -0.2506 & 6.5733e-09 & -0.3066 \\
			& 4096 & 61440 & 5.1319e-10 & -0.7520 & 4.1770e-08 & -0.4242 & 8.8976e-09 & -0.4368 \\
			& 16384 & 245760 & 2.5162e-10 & 1.0283 & 2.1292e-08 & 0.9722 & 3.6182e-09 & 1.2981 \\

		\hline
		
		5	& 16 & 336 & 4.8088e-09 & - & 1.7436e-07 & - & 4.0116e-08 & - \\
			& 64 & 1344 & 2.0183e-10 & 4.5745 & 1.4229e-08 & 3.6151 & 2.8580e-09 & 3.8111 \\
			& 256 & 5376 & 3.7226e-10 & -0.8832 & 3.1367e-08 & -1.1405 & 6.5744e-09 & -1.2019 \\
			& 1024 & 21504 & 3.2134e-10 & 0.2122 & 2.2226e-08 & 0.4970 & 5.7904e-09 & 0.1832 \\
			& 4096 & 86016 & 3.8469e-10 & -0.2596 & 1.4881e-08 & 0.5788 & 5.9937e-09 & -0.0498 \\
			& 16384 & 344064 & 9.6647e-10 & -1.3290 & 4.2430e-08 & -1.5117 & 1.3583e-08 & -1.1803 \\
		
		\hline

		6	& 16 & 448 & 5.0749e-10 & - & 1.9293e-08 & - & 7.7938e-09 & - \\
			& 64 & 1792 & 3.5543e-10 & 0.5138 & 1.5875e-08 & 0.2813 & 7.2192e-09 & 0.1105 \\
			& 256 & 7168 & 4.1596e-10 & -0.2269 & 1.9329e-08 & -0.2840 & 8.0760e-09 & -0.1618 \\
			& 1024 & 28672 & 3.9734e-10 & 0.0661 & 1.7754e-08 & 0.1226 & 6.7155e-09 & 0.2661 \\
			& 4096 & 114688 & 1.4872e-10 & 1.4177 & 6.3673e-09 & 1.4794 & 2.3250e-09 & 1.5303 \\
			& 16384 & 458752 & 9.0933e-10 & -2.6122 & 5.7407e-08 & -3.1725 & 7.6649e-09 & -1.7211 \\

		\hline
	\end{tabular}
\end{table}

\begin{figure}[htbp]
	\centering
	\begin{subfigure}{0.32\linewidth}
		\centering
		\includegraphics[width=\linewidth]{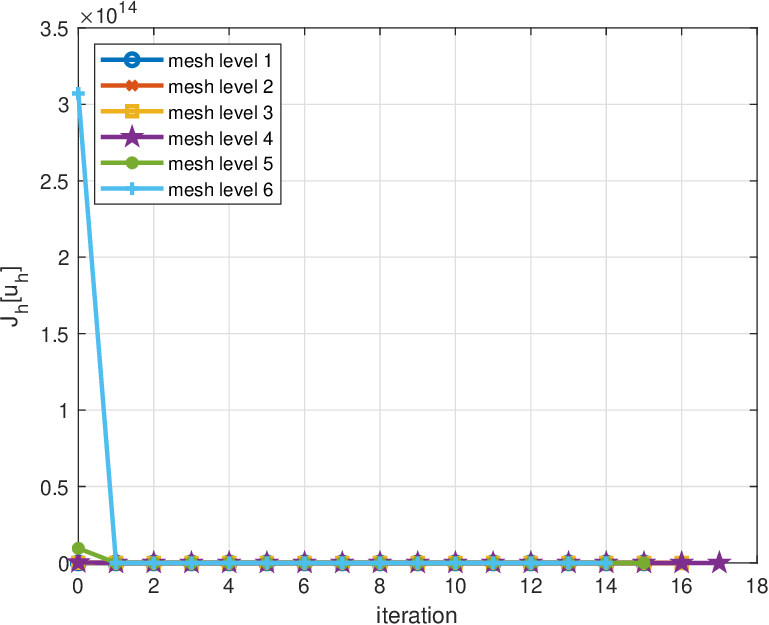}
		\caption{\(J_{h}(u_{h})\), \(k = 1\)}
	\end{subfigure}
	\begin{subfigure}{0.32\linewidth}
		\centering
		\includegraphics[width=\linewidth]{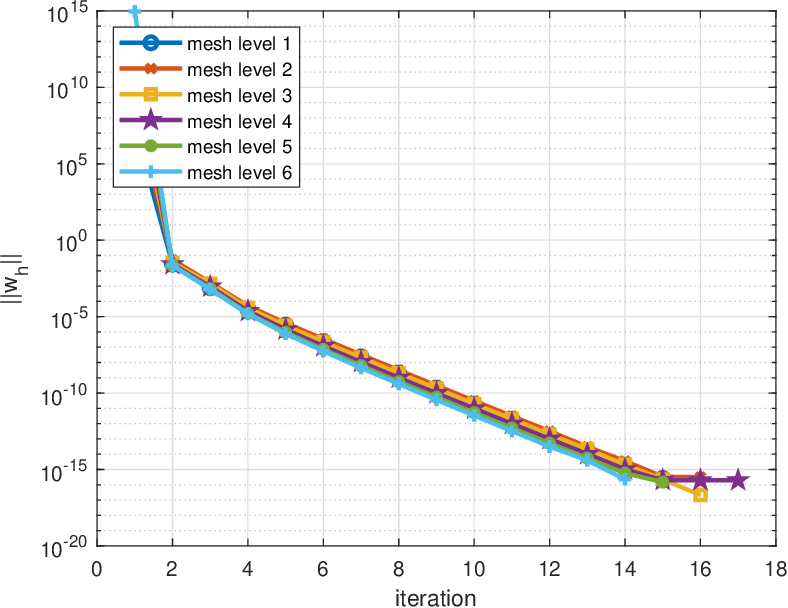}
		\caption{\(\norm{w_{h}}\), \(k = 1\)}
	\end{subfigure}
	\begin{subfigure}{0.32\linewidth}
		\centering
		\includegraphics[width=\linewidth]{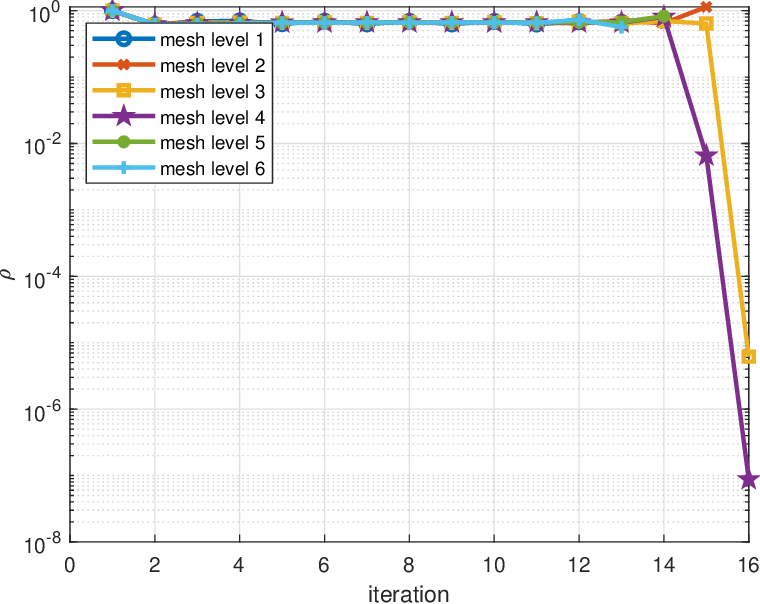}
		\caption{\(\rho\), \(k = 1\)}
	\end{subfigure}
	\\
	\begin{subfigure}{0.32\linewidth}
		\centering
		\includegraphics[width=\linewidth]{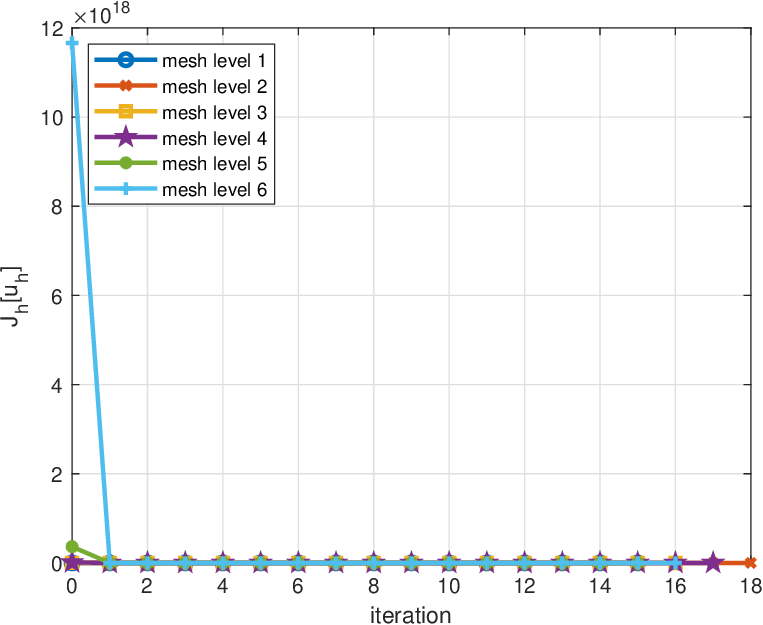}
		\caption{\(J_{h}(u_{h})\), \(k = 6\)}
	\end{subfigure}
	\begin{subfigure}{0.32\linewidth}
		\centering
		\includegraphics[width=\linewidth]{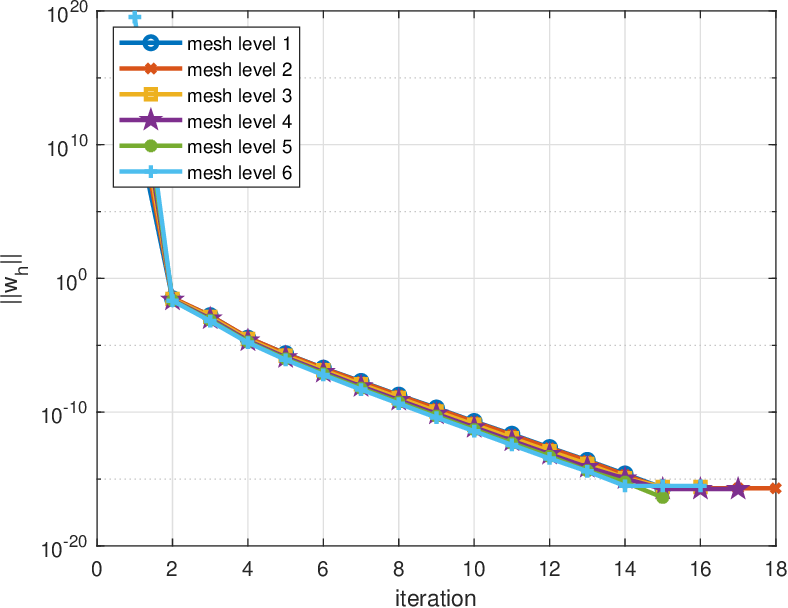}
		\caption{\(\norm{w_{h}}\), \(k = 6\)}
	\end{subfigure}
	\begin{subfigure}{0.32\linewidth}
		\centering
		\includegraphics[width=\linewidth]{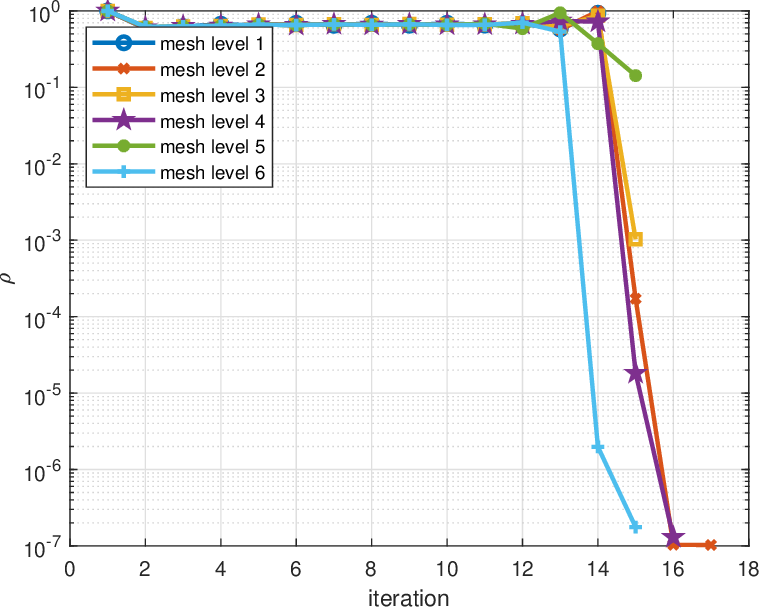}
		\caption{\(\rho\), \(k = 6\)}
	\end{subfigure}
	\caption{Convergence history for Example \ref{exm:sm}, with \(p = 3\). }
	\label{pic: 2D p-Laplace Dirichlet LDG, example 4, test2}
\end{figure}

\section{Conclusion}\label{se:con}

In this work, we study the high-order LDG method for the \(p\)-Laplace equation. Inspired by the variational minimization form of the PDE, we solve a discrete minimization problem instead of the original discrete weak problem to enhance stability, and a specially designed weighted preconditioner is employed to accelerate the nonlinear solver. In our analysis, we rigorously establish the solvability and equivalence between the discrete minimization problem and the original LDG scheme. For the first time, we provide (though non-optimal) \emph{a priori} error estimates for arbitrarily high-order polynomials under the assumption that the exact solution has sufficient regularity. Our error estimates demonstrate the expected potential for the LDG scheme to achieve high-order accuracy under reasonable assumptions, and the distance functional in our error estimates is non-equivalent to previous works. Moreover, under the same regularity assumption for the primal variable \(u\) and using sufficiently high-order polynomials, our estimated convergence rate is the same as that of HHO methods and similar to that of HDG methods. Our numerical results exhibit the best convergence rates for gradient variables in the LDG setting in most cases, while \emph{optimal} convergence rates are obtained for the primal variable when \(1 < p \leqslant 2\), which may be attributed to good local regimes. Additionally, the nonlinear solver achieves its desired \(hk\)-independent number of iterations, as shown in numerical examples. Possible future work includes improvements on the error estimates via dual arguments and via the local regime analysis.

\printbibliography[heading=bibintoc] 


\end{document}